\documentclass[12pt]{amsproc}

              [1996/10/25 v1.2b CONM-P Author Class]

\usepackage{amsmath,amsthm,amssymb, amscd, amsfonts}

\newcommand\SkPr{\mathcal P}

\newtheorem{theorem}{Theorem}[section]
\newtheorem{lema}[theorem]{Lemma}
\newtheorem{cor}[theorem]{Corollary}
\newtheorem{conjecture}[theorem]{Conjecture}
\newtheorem{prop}[theorem]{Proposition}

\theoremstyle{definition}
\newtheorem{definition}[theorem]{Definition}
\newtheorem{exa}[theorem]{Example}
\newtheorem{question}[theorem]{Question}

\newtheorem{obs}[theorem]{Remark}
\newtheorem{rmk}[theorem]{Remarks}

\numberwithin{equation}{section}

\def\pf{\begin{proof}}
\def\epf{\end{proof}}

\newcommand{\Ss}{{\mathcal S}}

\newcommand{\uno}{{\mathbb I}}

\renewcommand{\ker}{\mbox{\rm Ker\,}}

\newcommand{\Sym}{\mbox{\rm Sym\,}}

\newcommand{\ku}{{\Bbbk}}
\newcommand{\Z}{{\mathbb Z}}
\newcommand{\N}{{\mathbb N}}
\newcommand{\Q}{{\mathbb Q}}

\newcommand{\C}{{\mathbb C}}

\newcommand{\kv}{{\Q (v)}}

\newcommand{\ydh}{{}_{H}^{H}\mathcal{YD}}

\newcommand{\ydg}{{}_{\Gamma}^{\Gamma}\mathcal{YD}}
\newcommand{\ydw}{{}_{W}^{W}\mathcal{YD}}

\newcommand{\ydhsi}{{}_{H_{\sigma}}^{H_{\sigma}}\mathcal{YD}}

\newcommand{\ydzeta}{{}_{\Z}^{\Z}\mathcal{YD}}

\newcommand{\VGamma}{\widehat{\Gamma}}

\newcommand{\toba}{{\mathfrak B}}
\newcommand{\wtoba}{\widehat{\mathfrak B}}
\newcommand{\ftoba}{{\mathcal K}}
\newcommand{\Vdos }{{\mathcal V}(t, 2)}
\newcommand{\Vn}{{\mathcal V}(t, \theta)}
\newcommand{\Vuno }{{\mathcal V}(1, 2)}
\newcommand{\Vnuno}{{\mathcal V}(1, \theta)}

\newcommand{\End}{\mbox{\rm End\,}}
\newcommand{\Aut}{\mbox{\rm Aut\,}}

\newcommand{\im}{\mathop{\rm Im\,}}
\newcommand{\id}{\mathop{\rm id\,}}

\newcommand{\ord}{\mathop{\rm ord}}

\newcommand{\ad}{\mbox{\rm ad\,}}

\newcommand{\gr}{\mbox{\rm gr\,}}
\newcommand{\ev}{{\rm ev}}

\newcommand{\hecke}{{\mathcal H}_q(n)}

\theoremstyle{plain}

\begin{document}

\renewcommand{\baselinestretch}{1.2}
\renewcommand{\thefootnote}{}
\thispagestyle{empty}

\title{Pointed Hopf algebras}
\author[Andruskiewitsch and Schneider]{ Nicol\'as Andruskiewitsch and Hans-J\"urgen Schneider}
\address{Facultad de Matem\'atica, Astronom\'\i a y F\'\i sica\\
Universidad Nacional de C\'ordoba \\ (5000) Ciudad Universitaria
\\C\'ordoba \\Argentina}
\email{andrus@mate.uncor.edu}
\thanks{
\noindent Version of September 11, 2001. \\
This work was partially supported by ANPCyT, Agencia C\'ordoba Ciencia,
CONICET, DAAD, the Graduiertenkolleg of the Math. Institut (Universit\"at \ M\"unchen)  and Secyt (UNC)}
\address{Mathematisches Institut\\
Universit\"at \ M\"unchen\\
Theresienstra\ss e 39\\ \newline
D-80333 M\"unchen\\
Germany}
\email{hanssch@rz.mathematik.uni-muenchen.de}
\begin{abstract} This is a survey on pointed Hopf algebras over algebraically closed fields
of characteristic 0. We propose to classify pointed Hopf algebras $A$ by first determining
the graded Hopf algebra $\gr A$ associated to the coradical filtration of $A$.
The $A_{0}$-coinvariants elements form a braided Hopf algebra $R$ in the category of
Yetter-Drinfeld modules over the coradical $A_{0} = \ku \Gamma$, $\Gamma$ the group of
group-like elements of $A$, and $\gr A \simeq R \# A_{0}$. We call the braiding of the
primitive elements of $R$ the infinitesimal braiding of $A$. If this braiding is of
Cartan type \cite{AS2}, then it is often possible to determine $R$, to show that $R$
is generated as an algebra by its primitive elements and finally to compute all
deformations or liftings, that is pointed Hopf algebras such that
$\gr A \simeq R \# \ku \Gamma$.
In the last Chapter, as a concrete illustration of the method,
we describe explicitly all finite-dimensional pointed Hopf algebras
$A$ with abelian group of group-likes
$G(A)$ and infinitesimal braiding of type $A_{n}$ (up to some exceptional cases).
In other words, we compute all the liftings of type $A_n$; this result is our
main new contribution in this paper.

\end{abstract}

\maketitle

\tableofcontents



\section*{Introduction}

A Hopf algebra $A$ over a field $\ku$ is called {\em pointed} \cite{Sw}, \cite{M},
if all its simple left or right comodules are one-dimensional. The coradical $A_{0}$
of $A$ is the sum of all its simple subcoalgebras. Thus $A$ is pointed if and only
if $A_{0}$ is a group algebra.

We will always assume that the field $\ku$ is algebraically closed of characteristic
0 (although several results of the paper hold over arbitrary fields).

It is easy to see that $A$ is pointed if it is generated as an algebra by group-like
and skew-primitive elements. In particular, group algebras, universal enveloping
algebras of Lie algebras and the $q$-deformations of the universal enveloping
algebras of semisimple Lie algebras are all pointed.

\medbreak

An essential tool in the study of pointed Hopf algebras is the {\em coradical filtration}
$$A_{0} \subset A_{1} \subset \dots \subset A, \bigcup_{n\geq0} A_{n} = A$$
of $A$. It is dual to the filtration of an algebra by the powers of the Jacobson radical. For pointed Hopf algebras it is a Hopf algebra filtration, and the associated graded Hopf algebra $\gr A$ has a Hopf algebra projection onto $A_{0} = \ku \Gamma, \Gamma = G(A)$ the group of all group-like elements of $A$. By a theorem of Radford \cite{Ra}, $\gr A$ is a biproduct
$$\gr A \cong R \# \ku \Gamma,$$
where $R$ is a graded braided Hopf algebra in the category of left Yetter-Drinfeld modules over $\ku \Gamma$ \cite{AS2}.

\medbreak

This decomposition is an analog of the theorem of Cartier--Kostant--Milnor--Moore on the semidirect product decomposition of a cocommutative Hopf algebra into an infinitesimal and a group algebra part.

\medbreak

The vector space $V= P(R)$ of the primitive elements of $R$ is a Yetter-Drinfeld submodule. We call its braiding
$$c : V \otimes V \to V \otimes V$$
the {\bf infinitesimal braiding} of $A$. The infinitesimal braiding is the key to the structure of pointed Hopf algebras.

\medbreak

The subalgebra $\mathfrak{B}(V)$ of $R$ generated by $V$ is a braided Hopf subalgebra. As an algebra and coalgebra, $\mathfrak{B}(V)$ only depends on the infinitesimal braiding of $V$. In his thesis \cite{N} published in 1978, Nichols studied Hopf algebras of the form $\mathfrak{B}(V) \# \ku \Gamma$ under the name of bialgebras of type one. We call $\mathfrak{B}(V)$ the {\em Nichols algebra} of $V$. These Hopf algebras were found independently later by Woronowicz \cite{W} and other authors.

\medbreak

Important examples of Nichols algebras come from quantum groups \cite{Dr1}.
If $\mathfrak{g}$ is a semisimple Lie algebra, $U_{q}^{\ge 0}(\mathfrak{g})$,
$q$ not a root of unity, and the finite-dimensional Frobenius-Lusztig kernels $\mathfrak{u}_{q}^{\ge 0}(\mathfrak{g})$, $q$ a root of unity of order $N$,
are both of the form
$\mathfrak{B}(V) \# \ku \Gamma$ with $\Gamma = \mathbb{Z}^{\theta}$ resp. $(\mathbb{Z}/(N))^{\theta}, \theta \geq 1.$ (\cite{L3}, \cite{Ro1},
\cite{Sbg}, and \cite{L2}, \cite{Ro1}, \cite{Mu})
(assuming some technical conditions on $N$).

\medbreak

In general, the classification problem of pointed Hopf algebras has three parts:
\begin{itemize}
\item[(1)] Structure of the Nichols algebras $\mathfrak{B}(V)$.

\medbreak
\item[(2)] The lifting problem: Determine the structure of all pointed Hopf algebras $A$ with $G(A) = \Gamma$ such that $\gr A \cong  \mathfrak{B}(V) \# \ku \Gamma$.

\medbreak
\item[(3)] Generation in degree one: Decide which Hopf algebras $A$ are generated by group-like and skew-primitive elements, that is $\gr A$ is generated in degree one.
\end{itemize}
We conjecture that all finite-dimensional pointed Hopf algebras over an algebraically closed field of characteristic 0 are indeed generated by group-like and skew-primitive elements.

\medbreak

In this paper, we describe the steps of this program in detail and
explain the positive results obtained so far in this direction.
It is not our intention to give a complete survey on all aspects of pointed Hopf algebras.

\medbreak

We will mainly report on recent progress in the classification of pointed Hopf algebras  with {\em abelian} group of group-like elements.

If the group $\Gamma$ is abelian, and $V$ is a finite-dimensional Yetter-Drinfeld module,
then the braiding is given by a family of
non-zero scalars $q_{ij} \in \ku, 1 \leq i \leq \theta$, in the form
$$c(x_{i} \otimes x_{j}) = q_{ij} x_{j} \otimes x_{i}, \text{ where } x_{1},\dots,x_{\theta} \text{ is a basis of } V.$$
Moreover there are elements $g_{1}, \dots,g_{\theta} \in \Gamma$, and characters $\chi_{1},\dots,\chi_{\theta} \in \widehat{\Gamma}$ such that $q_{ij} = \chi_{j}(g_{i}).$ The group acts on $x_{i}$ via the character $\chi_{i}$, and $x_{i}$ is a $g_{i}$-homogeneous element with respect to the coaction of $\Gamma$.
We introduced braidings of Cartan type \cite{AS2} where
$$q_{ij} q_{ji} = q_{ii}^{a_{ij}},1\leq i,j \leq \theta, \text{ and } (a_{ij}) \text{ is a generalized Cartan matrix.}$$
If $(a_{ij})$ is a Cartan matrix of finite type, then the algebras $\toba(V)$ can be understood as twisting of the Frobenius-Lusztig kernels $\mathfrak{u}^{\ge 0}(\mathfrak{g})$, $\mathfrak{g}$ a semisimple Lie algebra.

By deforming the quantum Serre relations for simple roots which lie in two different
connected components of the Dynkin diagram, we define finite-dimensional pointed Hopf
algebras $\mathfrak{u}(\mathcal{D})$ in terms of a ''linking datum $\mathcal{D}$ of
finite Cartan type`` \cite{AS4}. They generalize the Frobenius-Lusztig kernels $\mathfrak{u}(\mathfrak{g})$ and are liftings of $\toba(V) \# \ku \Gamma$.

\medbreak
In some cases linking data of finite Cartan type are general enough to obtain
complete classification results.

For example, if $\Gamma= (\mathbb{Z}/(p))^{s}$, $p$ a prime $>17$ and $s\geq 1$,
we have determined the structure of all finite-dimensional
pointed Hopf algebras $A$ with
$G(A) \simeq \Gamma$. They are all of the form $\mathfrak{u}(\mathcal{D})$ \cite{AS4}.

Similar data allow a classification of infinite-dimensional pointed Hopf algebras
$A$ with abelian group $G(A)$, without zero divisors, with finite Gelfand-Kirillov
dimension and semisimple action of $G(A)$ on $A$,
in the case when the infinitesimal braiding is ''positive`` \cite{AS5}.

But the general case is more involved. We also have to deform the root vector
relations of the $\mathfrak{u}(\mathfrak{g})'s$.

\medbreak
The structure of pointed Hopf algebras $A$ with {\em non-abelian} group $G(A)$ is
largely unknown.  One basic open problem is to decide which finite groups appear as
groups of group-like elements of finite-dimensional pointed Hopf algebras which
are link-indecomposable in the sense of \cite{M1}.
In our formulation, this problem is the main part of the following question:
given a finite group $\Gamma$, determine all Yetter-Drinfeld modules $V$ over
$\ku\Gamma$  such that $\toba (V)$ is finite dimensional.
On the one hand, there are a number of severe constraints on $V$ \cite{Gn4}.
See also the exposition in  \cite[5.3.10]{bariloche}.
On the other hand, it is very hard to prove the finiteness of the dimension,
and in fact this has been done only for a few examples \cite{MiS}, \cite{FK},
\cite{FP} which are again related to root systems.
The examples over the symmetric groups in \cite{FK}
were introduced to describe the cohomology ring of the flag variety.
At this stage, the main difficulty is to decide
when certain Nichols algebras over non-abelian groups,
for example the symmetric  groups $\mathbb{S}_{n}$, are finite-dimensional.

\medbreak
The last Chapter provides a concrete illustration of the theory explained in this paper.
 We describe explicitly all finite-dimensional pointed Hopf algebras with abelian group
$G(A)$ and infinitesimal braiding of type $A_{n}$ (up to some exceptional cases). The
main results in this Chapter are new, and complete proofs are given. The only cases
which were known before are the easy case $A_{1}$ \cite{AS1}, and $A_{2}$ \cite{AS3}.

The new relations concern the root vectors $e_{i,j}, 1 \leq i < j \leq n+1$.
The relations $e_{i,j}^{N} = 0$ in $\mathfrak{u}_{q}^{\ge 0}(sl_{n + 1})$,
$q$ a root of unity of order $N$, are replaced by
$$e_{i,j}^{N} = u_{i,j} \text{ for a family } u_{i,j} \in \ku \Gamma,
1 \leq i < j \leq n+1,$$
depending on a family of free parameters in $\ku$.
See Theorem \ref{classification} for details.

\medbreak
Lifting of type $B_2$ was treated in \cite{BDR}.

\medbreak
To study the relations between a filtered object and its associated graded
object is a basic technique in modern algebra.
We would like to stress  that finite dimensional pointed Hopf algebras
enjoy a remarkable rigidity; it is seldom the case
that one is able to describe precisely
all the liftings of a graded object, as in this context.

\medbreak
\subsection*{Acknowledgements}  We would like to thank Jacques Alev, Mat\'\i as Gra\~na
and Eric M\"uller for suggestions and stimulating discussions.

\medbreak
\subsection*{Conventions} As said above, our ground field  $\ku$ is algebraically
closed field  of characteristic 0.   Throughout,
``Hopf algebra" means  ``Hopf algebra with bijective antipode".

We denote by $\tau: V \otimes W \to W \otimes V$ the usual transposition,
that is $\tau (v \otimes w) = w \otimes v$.

We use Sweedler's notation for the comultiplication and
coaction; but, to avoid confusions, we use the following variant
for the comultiplication of a braided Hopf algebra $R$:
$\Delta_R(r) = r^{(1)} \otimes r^{(2)}$.


\section{Braided Hopf algebras}\label{sectionbraided}

\subsection{Braided categories}\label{brcat}

\

Braided Hopf algebras play a central r\^ole in this paper. Although we have tried to minimize the use of categorical language, we briefly and informally recall the notion of a braided category which is the appropriate setting for braided Hopf algebras.

Braided categories were introduced in \cite{JS}.
We refer to \cite[Ch. XI, Ch. XIII]{Kassel} for a detailed exposition.
There is a hierarchy of categories with a tensor product functor:

\medbreak
a). A {\it monoidal} or {\it tensor} category is a collection $({\mathcal C},\otimes, a, \uno, \ell, r)$, where

\begin{itemize}
\item ${\mathcal C}$ is a category and $\otimes:{\mathcal C}\times
{\mathcal C}\to {\mathcal C}$ is a functor,

\medbreak
\item $\uno$ is an object of ${\mathcal C}$, and

\medbreak
\item  $a: V\otimes (W\otimes U) \to (V\otimes W) \otimes U$, $\ell: V\to V\otimes \uno$,
$r: V\to  \uno\otimes V$ are natural isomorphisms;

\end{itemize}

\noindent such that the so-called ``pentagon" and ``triangle" axioms are satisfied, see \cite[Ch. XI, (2.6) and (2.9)]{Kassel}. These axioms essentially express that the tensor product of a finite number of objects is well-defined, regardless of the place where parentheses are inserted; and that $\uno$ is a unit for the tensor product.

\medbreak

b). A  {\it braided (tensor)} category is a collection $({\mathcal C},\otimes, a, \uno, \ell, r, c)$, where

\begin{itemize}
\item  $({\mathcal C},\otimes, a, \uno, \ell, r)$ is a monoidal category and

\medbreak
\item  $c_{V, W}: V\otimes W \to W \otimes V$ is a natural isomorphism;
\end{itemize}

\noindent
such that the so-called ``hexagon" axioms are satisfied, see \cite[Ch. XIII, (1.3) and (1.4)]{Kassel}. A very important consequence of the axioms of a
braided category is the following equality for any objects $V, W, U$:
\begin{equation}\label{categbraid}
(c_{V,W}\otimes {\id}_{U}) ({\id}_{V}\otimes c_{U, W})(c_{U, V}\otimes {\id}_{W}) =
({\id}_{W}\otimes c_{U, V})(c_{U, W}\otimes {\id}_{V})({\id}_{U}\otimes c_{V, W}),
\end{equation}
see \cite[Ch. XIII, (1.8)]{Kassel}.
For simplicity we have omitted the associativity morphisms.

\medbreak

c). A  {\it symmetric} category is a braided category where $c_{V, W}c_{W, V} = {\id}_{W\otimes V}$ for all objects $V$, $W$.
Symmetric categories have been studied since the pioneering work of Mac Lane.

\medbreak

d). A left dual of an object $V$ of a monoidal category, is a triple $(V^*, \ev_V, b_V)$,
where $V^*$ is another object and $\ev_V: V^* \otimes V \to \uno$,
$b_V: \uno \to V\otimes V^*$ are morphisms such that the compositions
$$\begin{CD}
V @>>> \uno \otimes V @>{b_V \otimes \, \id_{V}}>> V\otimes V^*\otimes V
@>{\id_{V} \otimes \, \ev_{V}}>> V \otimes \uno @>>> V
\end{CD}$$
and
$$\begin{CD}
V^* @>>> V^* \otimes \uno @>{\id_{V^*} \otimes \, b_{V}}>> V^*\otimes V\otimes V^*
@>{\ev_{V} \otimes \, \id_{V^*}}>> \uno \otimes V^* @>>> V^*
\end{CD}$$
are, respectively, the identity of $V$ and $V^*$. A  braided category is {\it rigid} if any object $V$ admits a left dual \cite[Ch. XIV, Def. 2.1]{Kassel}.

\medbreak
\subsection{Braided vector spaces and Yetter-Drinfeld modules}\label{bvs-yd}

\

We begin with the fundamental
\begin{definition}\label{braidedvec}
Let $V$ be a vector space and $c: V\otimes V\to V\otimes V$ a linear isomorphism. Then $(V, c)$ is called a
{\it braided vector space}, if $c$ is a solution of the {\it braid equation},
that is
\begin{equation}\label{braideqn}
(c\otimes \id) (\id\otimes c) (c\otimes \id) = (\id\otimes c) (c\otimes \id) (\id\otimes c).
\end{equation}
\end{definition}
It is well-known that the braid equation is equivalent to the {\it quantum
Yang-Baxter equation}:
\begin{equation}\label{qybe}
R_{12}R_{13}R_{23} =R_{23} R_{13} R_{12}.
\end{equation}
Here we use the standard notation: $R_{13}: V\otimes V \otimes V \to V\otimes V \otimes V$ is the map given by $\sum_j r_j \otimes \id \otimes r^j$, if $R = \sum_j r_j \otimes r^j$.
Similarly for $R_{12}$, $R_{23}$.

The equivalence between solutions of \eqref{braideqn} and solutions of
\eqref{qybe} is given by the equality $c = \tau \circ R$. For this reason, some authors
call \eqref{braideqn} the quantum Yang-Baxter equation.

\medbreak

An easy and for this paper important example is given by a family of non-zero scalars $q_{ij} \in \ku, i,j \in I$, where $V$ is a vector space with basis $x_{i}, i \in I.$ Then
$$c(x_{i} \otimes x_{j}) = q_{ij}x_{j} \otimes x_{i}, \text{for all  } i,j \in I$$
is a solution of the braid equation.

Examples of braided vector spaces come from braided categories.
In this article, we are mainly concerned with examples related to
the notion of Yetter-Drinfeld modules.

\medbreak

\begin{definition}\label{ydr} Let $H$ be a Hopf algebra.
A (left) {\it Yetter-Drinfeld module} $V$
over $H$ is simultaneously a left $H$-module and a left $H$-comodule satisfying the compatibility condition
\begin{equation}\label{compcond}
\delta(h.v) = h_{(1)}v_{(-1)}\Ss h_{(3)} \otimes h_{(2)}.v_{(0)}, \qquad v\in V, \, h\in H. \end{equation}
\end{definition}
We denote by $\ydh$ the category  of Yetter-Drinfeld modules over $H$;
the morphisms in this category preserve both the action and the coaction of $H$.
The category $\ydh$ is a braided monoidal category;
indeed the tensor product of two Yetter-Drinfeld modules is again
a Yetter-Drinfeld module, with the usual tensor product module
and comodule structure. The compatibility
condition \eqref{compcond} is not difficult to verify.

\medbreak

For any two Yetter-Drinfeld-modules $M$ and $N$,
the braiding $c_{M,N}: M\otimes N \to N\otimes M$
is given by
\begin{equation}\label{braid}
c_{M,N}(m\otimes n) = m_{(-1)}.n \otimes m_{(0)}, \qquad m\in M, \, n\in N.
\end{equation}

\medbreak

The subcategory of $\ydh$ consisting of finite dimensional Yetter-Drinfeld
modules is rigid. Namely,  if
$V \in  \ydh$ is finite-dimensional, the dual $V^{*} = \text{Hom}(V,\ku)$ is in $\ydh$ with the following action and coaction:
\begin{itemize}
\item $(h \cdot f)(v) = f(\Ss(h)v)$ for all $h \in H$, $f \in V^{*}$, $v \in V$.

\medbreak
\item If $f \in V^{*}$, then $ \delta(f) = f_{(-1)} \otimes f_{(0)}$ is
determined by the equation $$ f_{(-1)}f_{(0)}(v) = \Ss^{-1}(v_{-1})f(v_{0}),
\qquad  v \in V.$$
\end{itemize}

Then the usual evaluation and coevaluation maps are morphisms in $\ydh$.

\medbreak

Let $V$, $W$ be two finite-dimensional Yetter-Drinfeld modules over $H$.
We shall consider the isomorphism $\Phi : W^{*} \otimes V^{*} \to
(V \otimes W)^{*}$ given by
\begin{equation}\label{isoduals}
\Phi(\varphi \otimes \psi)(v \otimes w) = \psi(v) \varphi(w), \qquad
\varphi \in W^{*}, \psi \in V^{*}, v \in V, w \in W.
\end{equation}

\medbreak
\begin{obs}
We see that a Yetter-Drinfeld module is a braided vector space.
Conversely, a braided vector space $(V, c)$ can be realized as a
Yetter-Drinfeld module
over some Hopf algebra $H$ if and only if  $c$ is {\it rigid} \cite{Tk1}.
If this is the case, it can be realized in many different ways.
\end{obs}

\medbreak
We recall that a Hopf bimodule over a Hopf algebra $H$
is simultaneously a bimodule
and a bicomodule satisfying all possible compatibility conditions.
The category $^H_H{\mathcal M}^H_H$ of all Hopf bimodules over $H$
is a braided category.
The category $\ydh$ is equivalent, as a braided category, to the
category of Hopf bimodules. This was essentially first observed in \cite{W}
and then independently in \cite[Appendix]{AD}, \cite{Sbg}, \cite{Ro1}.

\medbreak
If $H$ is a finite dimensional Hopf algebra,
then the category $\ydh$ is equivalent to the category of modules
over the double of $H$ \cite{Mj0}. The braiding in $\ydh$ corresponds to the braiding
given by the ``canonical" $R$-matrix of the double.
In particular, if $H$ is a semisimple Hopf algebra then $\ydh$ is a semisimple category.
Indeed, it is known that the double of a semisimple Hopf algebra is again semisimple.

\medbreak
The case of Yetter-Drinfeld modules over group algebras is especially important
for the applications to pointed Hopf algebras.
If $H = \ku \Gamma$, where $\Gamma$ is a group, then an $H$-comodule $V$
is just a $\Gamma$-graded vector space: $V = \oplus_{g\in \Gamma}V_{g}$,
where $V_{g} = \{v\in V \mid \delta(v) = g\otimes v\}$. We will write $\ydg$ for the category of Yetter-Drinfeld modules over $\ku \Gamma$, and say that $V \in \ydg$ is a Yetter-Drinfeld module over $\Gamma$ (when  the field is fixed).

\begin{obs}\label{ydgroup}
Let $\Gamma$ be a group, $V$ a left $\ku \Gamma$-module, and a left  $\ku \Gamma$-comodule with grading  $V = \oplus_{g\in \Gamma}V_{g}$. We define a linear isomorphism $c : V \otimes V \to V \otimes V $ by
\begin{equation}\label{grbraid}
c(x \otimes y) = gy \otimes x, \text{ for all } x \in V_{g},\  g \in \Gamma,\  y \in V.
\end{equation}
Then
\begin{itemize}
\item [a)] $V \in \ydg$ if and only if $gV_{h} \subset V_{ghg^{-1}} \text{ for all } g,h \in \Gamma.$

\medbreak
\item [b)] If $V \in \ydg$, then $(V,c)$ is a braided vector space.

\medbreak
\item [c)] Conversely, if V is a faithful $\Gamma$-module (that is, if for all $g \in \Gamma, gv = v$ for all $v \in V$, implies $g=1$), and if $(V,c)$ is a braided vector space, then $V \in \ydg$.
\end{itemize}
\pf
a) is clear from the definition.

By applying both sides of the braid equation to elements of the form $x \otimes y \otimes z, x \in V_{g}, y \in V_{h}, z \in V,$ it is easy to see that $(V,c)$ is a braided vector space if and only if
\begin{equation}\label{braidingsgr}
c(gy \otimes gz) = ghz \otimes gy, \text{ for all } g,h \in \Gamma, \  y \in V_{h},
\  z \in V.
\end{equation}
Let us write $gy = \sum_{a \in \Gamma} x_{a}$, where $x_{a} \in V_{a}$ for all $a \in \Gamma$. Then $c(gy \otimes gz) = \sum_{a \in \Gamma} agz \otimes x_{a}.$ Hence \eqref{braidingsgr} means that $agz = ghz$, for all $z \in V$ and $a \in \Gamma$ such that the homogeneous component $x_{a}$ is not zero. This proves b) and c).
\epf

\end{obs}
\begin{obs}\label{abelyd}
If $\Gamma$ is abelian,   a Yetter-Drinfeld
module over $H = \ku \Gamma$ is nothing but a
$\Gamma$-graded $\Gamma$-module.

Assume that $\Gamma$ is abelian and furthermore
that the action of $\Gamma$ is diagonalizable (this is always the
case if $\Gamma$ is finite). That is, $V = \oplus_{\chi\in \VGamma}V^{\chi}$,
where $V^{\chi} = \{v\in V \mid  gv = \chi(g)v \text{ for all } g \in \Gamma\}$.
Then
\begin{equation}\label{dobl-grad}
V = \oplus_{g\in \Gamma, \chi\in \VGamma}V_{g}^{\chi},
\end{equation}
\noindent where $V_{g}^{\chi} = V^{\chi} \cap V_{g}$. Conversely, any
vector space with a decomposition \eqref{dobl-grad} is a Yetter-Drinfeld module
over $\Gamma$. The braiding is given by
$$ c(x\otimes y) = \chi(g) y\otimes x,  \text{ for all } x\in V_{g},
\ g \in \Gamma, \  y\in V^{\chi}, \  \chi \in \VGamma.$$
\end{obs}

\medbreak
It is useful to characterize abstractly those braided vector spaces which come from Yetter-Drinfeld modules over groups or abelian groups. The first part of the following definition is due to M. Takeuchi.

\begin{definition}\label{grouptype} Let $(V, c)$ be a finite dimensional
braided vector space.

\begin{itemize}
\item  $(V, c)$   is of {\it group type}
if there exists a basis $x_{1}, \dots, x_{\theta}$ of $V$ and elements  $g_{i}(x_{j}) \in V$ for all $i,j$ such that
\begin{equation}
c(x_{i}\otimes x_{j}) = g_{i}(x_{j})\otimes x_{i}, \qquad 1 \le i, j \le \theta;
\end{equation}
necessarily $g_{i} \in GL(V)$.

\medbreak
\item  $(V, c)$ is of {\it finite group type} (resp. of {\it abelian group type}) if it is of group type and the subgroup of $GL(V)$ generated by $g_{1}, \dots, g_{\theta}$ is finite (resp. abelian).

\medbreak
\item  $(V, c)$ is of {\it diagonal type} if $V$ has
a basis $x_{1}, \dots, x_{\theta}$ such that
\begin{equation}\label{br-abfingrt}
c(x_{i}\otimes x_{j}) = q_{ij}x_{j}\otimes x_{i},\qquad 1 \le i, j \le \theta,
\end{equation}
for some $q_{ij}$ in $\ku$. The matrix $(q_{ij})$ is called the {\it matrix} of the braiding.

\medbreak
\item If $(V, c)$ is of  diagonal type, then we say that it
is {\it indecomposable}  if for all $i \neq j$,
there exists a sequence $i= i_{1}$, $i_{2}$, \dots, $i_{t} = j$
of elements of $\{1, \dots, \theta\}$ such that $q_{i_s, i_{s+1}}
q_{i_{s+1}, i_s} \neq 1$, $1\le s \le t-1$.
Otherwise, we say that the matrix is decomposable. We can also refer then
to the components of the matrix.
\end{itemize}
\end{definition}

\medbreak

If  $V \in \ydg$ is finite-dimensional with braiding $c$, then $(V,c)$
is of group type by \eqref{braid}. Conversely, assume that $(V,c)$ is a
finite-dimensional braided vector space of group type. Let $\Gamma$ be the
subgroup of $GL(V)$ generated by $g_{1}, \dots,g_{\theta}$. Define a coaction
by $\delta(x_{i}) = g_{i} \otimes x_{i}$ for all $i$. Then $V$ is a Yetter-Drinfeld
module over $\Gamma$ with braiding $c$ by Remark \ref{ydgroup}, c).

\medbreak
A braided vector space of diagonal type is clearly of abelian group
type; it is of finite group type if the $q_{ij}$'s are roots of one.

\medbreak
\subsection{Braided Hopf algebras}

\

The notion of ``braided Hopf algebra" is one of the basic features of braided categories.
We will deal in this paper only with braided Hopf algebras in categories of Yetter-Drinfeld modules, mainly over a group algebra.

\medbreak

Let $H$ be a Hopf algebra.
First, the tensor product in $\ydh$ allows us to define algebras and coalgebras
in $\ydh$. Namely,  an algebra in the category $\ydh$
is an associative algebra $(R, m)$, where $m: R\otimes R \to R$ is the product,
with unit $u: \ku \to R$, such that $R$ is a Yetter-Drinfeld module over $H$
and both $m$ and $u$ are morphisms in $\ydh$.

\medbreak
Similarly,  a coalgebra in the category $\ydh$
is a coassociative coalgebra $(R, \Delta)$,
where $\Delta: R\to R \otimes R$ is the coproduct,
with counit $\varepsilon: R \to \ku$,
such that $R$ is a Yetter-Drinfeld module over $H$
and both $\Delta$ and $\varepsilon$  are morphisms in $\ydh$.

\medbreak
Let now $R$, $S$ be two algebras in $\ydh$. Then the braiding
$c:S \otimes R \to R\otimes S$ allows us to provide the
Yetter-Drinfeld module $R\otimes S$ with a "twisted"
algebra structure in $\ydh$. Namely, the product in $R\otimes S$ is
$m_{R\otimes S} := (m_{R}\otimes m_{S}) (\id \otimes c\otimes \id)$:
$$\begin{CD}
R\otimes S \otimes R\otimes S @>>> R\otimes S\\
@V\id \otimes c\otimes \id VV @VV = V \\
R\otimes R \otimes S\otimes S  @>m_R \otimes m_S>> R\otimes S.
\end{CD}$$

We shall denote this algebra by $R \underline{\otimes} S$.
The difference with the usual tensor product algebra is the presence
of the braiding $c$ instead of the usual transposition $\tau$.

\begin{definition} A {\it braided   bialgebra in $\ydh$} is a collection
$(R, m, u, \Delta, \varepsilon)$, where

\begin{itemize}
\item $(R, m,u)$ is an algebra in $\ydh$.

\medbreak
\item  $(R, \Delta,\varepsilon)$  is a coalgebra in $\ydh$.

\medbreak
\item $\Delta: R\to R \underline{\otimes} R$ is a morphism of algebras.

\medbreak
\item $u: \ku \to R$ and $\varepsilon: R \to \ku$ are morphisms of algebras.
\end{itemize}

We say that it is a {\it braided  Hopf algebra in $\ydh$} if in addition:

\begin{itemize}
\item The identity is convolution invertible in $\End(R)$; its inverse
is the antipode of $R$.
\end{itemize}

A  {\it graded} braided Hopf algebra in $\ydh$ is a braided Hopf algebra
$R$ in $\ydh$ provided with a grading $R= \oplus_{n\ge 0} R(n)$ of Yetter-Drinfeld modules, such that $R$ is a graded algebra and a
graded coalgebra.
\end{definition}

\begin{obs}
There is a non-categorical version of braided Hopf algebras, see \cite{Tk1}.
Any braided Hopf algebra in $\ydh$ gives rise to a braided Hopf algebra
in the sense of  \cite{Tk1} by forgetting the action and coaction, and preserving the multiplication, comultiplication and braiding.
For the converse see  \cite[Th. 5.7]{Tk1}. Analogously, one can define
graded braided Hopf algebras in the spirit of  \cite{Tk1}.
\end{obs}

Let $R$ be a finite-dimensional Hopf algebra
in $\ydh$. The dual $S = R^{*}$ is a braided Hopf algebra in $\ydh$
with multiplication $\Delta_R^{*} \Phi$ and comultiplication $\Phi^{-1} m_R^{*}$, {\it cf.} \eqref{isoduals};
this is $R^{*bop}$ in the notation of \cite[Section 2]{AG}.

\medbreak
In the same way, if $R = \oplus_{n\ge 0} R(n)$ is a graded braided Hopf algebra
in $\ydh$ with finite-dimensional homogeneous components,
then the graded dual $S = R^{*} =  \oplus_{n\ge 0} R(n)^{*}$
is a graded braided Hopf algebra in $\ydh$.

\medbreak
 \subsection{Examples. The quantum binomial formula}\label{quantbinomial}

\

We shall provide many examples of braided Hopf algebras in Chapter
\ref{nicholsalgs}. Here we discuss a very simple class
of braided Hopf algebras.

We first recall the well-known quantum binomial formula. Let $U$ and $V$ be elements of an associative algebra over $\ku[q]$, $q$ an indeterminate, such that $VU =qUV$. Then
\begin{equation}\label{qbf}
(U + V)^n = \sum_{1\le i \le n} \binom{n}{i}_q U^{i} V^{n - i}, \qquad \text{ if } n \ge 1.
\end{equation}
Here
$$
\binom{n}{i}_q = \frac{(n)_q!}{(i)_q!(n - i)_q!}, \quad \text{ where } (n)_q! =
 \prod_{1\le i \le n} (i)_q, \quad \text{ and } (i)_q =  \sum_{0\le j \le i-1} q^j.
$$
By specialization, \eqref{qbf} holds for $q\in \ku$. In particular,
if $U$ and $V$ are elements of an associative algebra over $\ku$,
and $q$ is a primitive  $n$-th root of 1, such that $VU =qUV$ then
\begin{equation}\label{qbfpwr}
(U + V)^n = U^n + V^n .
\end{equation}

\begin{exa}\label{qlslemma} Let $(q_{ij})_{1\le i, j\le \theta}$
be a matrix such that
\begin{align}
q_{ij}q_{ji} &= 1,\quad  1\le i, j\le \theta, \ i \neq j.
\end{align}
Let $N_{i}$ be the order of $q_{ii}$, when this is finite.

Let $R$ be the algebra presented by generators $x_{1}, \dots, x_{\theta}$ with relations
\begin{flalign}\label{qls2bis}  &x_{i}^{N_i} = 0,
\quad
\text{if }\ord q_{ii} < \infty.& \\
\label{qls2tris} &x_{i}x_{j} = q_{ij}x_{j}x_{i}, \quad 1\le i < j\le \theta.&
\end{flalign}
Given a group $\Gamma$ and elements $g_{1}, \dots, g_{\theta}$ in the center of $\Gamma$, and characters $\chi_{1}, \dots, \chi_{\theta}$ of $\Gamma$, there exists a unique structure of Yetter-Drinfeld module over $\Gamma$ on $R$, such that
$$
x_i \in R^{\chi_{i}}_{g_{i}}, \qquad 1\le i \le \theta.
$$
Note that the braiding is determined by
$$
c(x_{i} \otimes x_{j}) = q_{ij} \,x_{j} \otimes x_{i}, \text{ where } q_{ij} = \chi_{j}(g_{i}), \quad  1\le i, j\le \theta.
$$
Furthermore, $R$ is a braided Hopf algebra with the comultiplication given
by $\Delta(x_{i}) = x_{i} \otimes 1 + 1\otimes x_{i}$. To check that the comultiplication
preserves \eqref{qls2bis} one uses \eqref{qbfpwr}; the verification for \eqref{qls2tris} is easy. We know \cite{AS1} that $\dim R$ is infinite unless all the
orders of $q_{ii}$'s are finite; in this last case, $\dim R = \prod_{1\le i\le \theta}N_{i}$.
We also have $P(R) = \oplus_{1\le i\le \theta}\ku x_{i}$.
\end{exa}

\medbreak
\subsection{Biproducts, or bosonizations}\label{bosoniz}

\

Let $A$, $H$ be Hopf algebras and $\pi: A \to H$ and $\iota: H \to A$ Hopf algebra homomorphisms.
Assume that $\pi\iota = \id_{H}$, so that $\pi$ is surjective, and $\iota$
is injective. By analogy with elementary group theory, one would like
to reconstruct $A$ from $H$ and the kernel of $\pi$ as a semidirect product. However,
the natural candidate for the kernel of $\pi$ is the algebra of
coinvariants
$$R:= A^{\text{co } \pi} = \{a\in A: (\id\otimes \pi)\Delta (a) = a\otimes 1\}$$
which is {\it not}, in general, a Hopf algebra. Instead, $R$
is a braided Hopf algebra in $_{H}^{H}\mathcal{YD}$ with the following structure:

\begin{itemize}
\item The action  $\cdot$ of $H$ on $R$
is the restriction of the adjoint action
(composed with $\iota$).

\medbreak
\item The coaction is $(\pi \otimes \id)\Delta$.

\medbreak
\item $R$ is a subalgebra of $A$.

\medbreak
\item The comultiplication is $\Delta_R(r) =
r_{(1)} \iota\pi\mathcal S(r_{(2)}) \otimes r_{(3)}$, for all $r \in R$.
\end{itemize}

Given a braided Hopf algebra $R$ in $\ydh$,
one can consider  the  {\it bosonization} or {\it biproduct} of $R$ by $H$ \cite{Ra}, \cite{Mj1}.
This is a usual Hopf algebra $R\# H$, with underlying vector space
$R\otimes H$,
whose multiplication and comultiplication are given by

\begin{equation}\label{smash1}
\begin{aligned}(r\# h)(s\# f) &= r (h_{(1)} \cdot s)\# h_{(2)}f, \\ \Delta(r\# h) &=
r^{(1)} \# (r^{(2)})_{(-1)} h_{(1)} \otimes (r^{(2)})_{(0)}
 \#  h_{(2)}.\end{aligned}\end{equation}

The maps $\pi: R\# H \to H$ and $\iota: H \to R\# H$,
$\pi(r\# h) = \epsilon(r)h$,
$\iota(h) = 1\# h$, are Hopf algebra homomorphisms; we have
$R = \{a\in R\# H: (\id\otimes \pi)\Delta (a) = a\otimes 1\}$.

\medbreak

Conversely, if $A$ and $H$ are Hopf algebras as above
and $R = A^{\text{co } \pi}$, then $A \simeq R\# H$.

\medbreak

Let $\vartheta: A \to R$ be the map given by
$\vartheta(a) = a_{(1)} \iota \pi\Ss(a_{(2)})$.
Then
\begin{equation}\label{varthetadjoint}
\vartheta(ab) = a_{(1)}\vartheta(b)  \iota\pi\Ss(a_{(2)}),\end{equation} for all
$a,b \in A$,
 and $\vartheta(\iota(h)) = \varepsilon(h)$ for all $h\in H$; therefore,
for all $a\in A$, $h\in H$, we have
$\vartheta(a \iota(h)) = \vartheta(a) \varepsilon(h)$ and
\begin{equation}\label{varthetamodule}
\vartheta(\iota(h)a) = h \cdot \vartheta(a).
\end{equation}
Notice also that $\vartheta$ induces a coalgebra isomorphism $A/A \iota(H)^{+}
\simeq R$. In fact, the isomorphism $A \to R\# H$ can be expressed explicitly as
$$
a \mapsto \vartheta(a_{(1)}) \# \pi(a_{(2)}), \qquad a\in A.
$$

If $A$ is a Hopf algebra, the  adjoint representation
$\ad$
of $A$ on itself is given by
$$\ad x(y)= x_{(1)} y \Ss(x_{(2)}).$$
If $R$ is a braided Hopf algebra in $_{H}^{H}\mathcal{YD}$,
then there is also a braided adjoint
representation $\ad_c$ of $R$ on itself defined by
$$\ad_c x(y)=\mu(\mu\otimes\Ss_{R})(\id\otimes c)(\Delta_{R}\otimes\id)(x\otimes y),$$
where $\mu$ is the multiplication and
$c\in\End(R\otimes R)$ is the braiding.
Note that if $x\in {\mathcal P} (R)$ then the braided adjoint representation
of $x$ is just
\begin{equation}\label{br-adj}
\ad_cx(y)=\mu(\id-c)(x\otimes y) =: [x, y]_{c}.
\end{equation}
For any $x,y \in R$, we call $[x, y]_{c} := \mu(\id-c)(x\otimes y)$
 a {\it braided commutator}.

\medbreak
When $A = R\# H$, then for all $b,d\in R$,
\begin{equation}\label{bradj}
\ad_{(b\# 1)}(d\# 1)=(\ad_cb(d))\# 1.\end{equation}

\medbreak
\subsection{Some properties of braided Hopf algebras}\label{propbraided-hopf}

\

In this Section, we first collect several useful facts about braided Hopf algebras in the category of Yetter-Drinfeld modules over an abelian group $\Gamma$.
We begin with some identities on braided commutators.

\medbreak
In the following two Lemmas, $R$ denotes a
braided Hopf algebra in $\ydg$. Let $a_{1}, a_{2}, \dots \in R$ be elements
such that $a_{i} \in R^{\chi_i}_{g_i}$, for some  $\chi_i\in \VGamma$,
$g_i\in \Gamma$.

\begin{lema}\label{brcomm} (a).
\begin{equation}
\label{brcomm2}
[[a_{1}, a_{2}]_{c}, a_{3}]_{c} + \chi_{2}(g_{1}) a_{2}[a_{1}, a_{3}]_{c}
= [a_{1}, [a_{2}, a_{3}]_{c}]_{c} + \chi_{3}(g_{2}) [a_{1}, a_{3}]_{c} a_{2}.
\end{equation}

\medbreak
(b). If $[a_{1}, a_{2}]_{c} = 0$ and $[a_{1}, a_{3}]_{c} = 0$
then $[a_{1}, [a_{2}, a_{3}]_{c}]_{c} = 0$.

\medbreak
(c). If $[a_{1}, a_{3}]_{c} = 0$ and $[a_{2}, a_{3}]_{c} = 0$
then $[[a_{1}, a_{2}]_{c}, a_{3}]_{c} = 0$.

\medbreak
(d). Assume that $\chi_{1}(g_{2})\chi_{2}(g_{1})\chi_{2}(g_{2}) = 1$.
Then

\begin{equation}\label{brcomm1}
[[a_{1}, a_{2}]_{c}, a_{2}]_{c} = \chi_{2}(g_{1}) \chi_{1}(g_{2})^{-1}
[a_{2}, [a_{2}, a_{1}]_{c}]_{c}
\end{equation}  \end{lema}

\pf Left to the reader. \epf

\medbreak
The following technical Lemma will be used at a crucial point in Section \ref{nich-A_n}.

\begin{lema}\label{brcomm-an} Assume that $\chi_{2}(g_{2}) \neq -1$ and
\begin{align}
\label{techni2}
\chi_{1}(g_{2})\chi_{2}(g_{1}) \chi_{2}(g_{2}) &= 1, \\
\label{techni3}
\chi_{2}(g_{3})\chi_{3}(g_{2}) \chi_{2}(g_{2}) &= 1.
\end{align}
If
\begin{align}
\label{techni5}
[a_{2}, [a_{2}, a_{1}]_{c}]_{c} &= 0, \\
\label{techni6}
[a_{2}, [a_{2}, a_{3}]_{c}]_{c} &= 0, \\
\label{techni7}
[a_{1}, a_{3}]_{c} &= 0, \end{align}
then
\begin{equation}\label{brcomm3}
[[[a_{1}, a_{2}]_{c}, a_{3}]_{c}, a_{2}]_{c} = 0.\end{equation}  \end{lema}

\pf We compute:
\begin{align*}
[[[a_{1}, a_{2}]_{c}, a_{3}]_{c}, a_{2}]_{c}
&= a_{1} a_{2}a_{3}a_{2}
- \chi_{2}(g_{1}) \, a_{2} a_{1} a_{3}a_{2}
  - \chi_{3}(g_{1})\chi_{3}(g_{2}) \, a_{3} a_{1} a_{2}^{2}\\
& + \chi_{3}(g_{1})\chi_{3}(g_{2}) \chi_{2}(g_{1}) \,  a_{3}a_{2} a_{1} a_{2}
- \chi_{2}(g_{1})\chi_{2}(g_{2}) \chi_{2}(g_{3}) \, a_{2} a_{1} a_{2}a_{3} \\
& + \chi_{2}(g_{1})^2\chi_{2}(g_{2}) \chi_{2}(g_{3}) \, a_{2}^{2} a_{1} a_{3}
+ \chi_{2}(g_{1})\chi_{2}(g_{2}) \chi_{2}(g_{3}) \chi_{3}(g_{1})\chi_{3}(g_{2}) \,
a_{2} a_{3} a_{1} a_{2} \\
&- \chi_{2}(g_{1})^2\chi_{2}(g_{2}) \chi_{2}(g_{3}) \chi_{3}(g_{1})\chi_{3}(g_{2})
 \, a_{2} a_{3}a_{2} a_{1}.
\end{align*}
We index consecutively the terms in the right-hand side by roman numbers:
$(I), \dots, (VIII)$. Then $(II) + (VII) = 0$, by \eqref{techni3} and \eqref{techni7}.
Now,
\begin{align*}
(I) &= \frac{1}{\chi_{3}(g_{2}) (1 + \chi_{2}(g_{2}))} a_{1} a_{2}^{2} a_{3}
+ \frac{\chi_{2}(g_{2})\chi_{3}(g_{2})}{1 + \chi_{2}(g_{2})} a_{1}  a_{3} a_{2}^{2}
\\ \\&= \frac{1}{\chi_{3}(g_{2}) (1 + \chi_{2}(g_{2}))} a_{1} a_{2}^{2} a_{3}
+ \frac{\chi_{2}(g_{2})\chi_{3}(g_{2})\chi_{3}(g_{1})}{1 + \chi_{2}(g_{2})}
a_{3} a_{1} a_{2}^{2}
\\ \\&= (I\, a) + (I\, b),
\end{align*}
by \eqref{techni6} and \eqref{techni7}.
By the same equations \eqref{techni6} and \eqref{techni7}, we also have
\begin{align*}
(VIII) &=  - \frac{\chi_{2}(g_{1})^2\chi_{2}(g_{2}) \chi_{2}(g_{3}) \chi_{3}(g_{1})}{1 + \chi_{2}(g_{2})} a_{2}^{2} a_{3} a_{1}
- \frac{\chi_{2}(g_{1})^2\chi_{2}(g_{2})^2 \chi_{2}(g_{3})
\chi_{3}(g_{1})\chi_{3}(g_{2})^2}{1 + \chi_{2}(g_{2})}   a_{3} a_{2}^{2} a_{1}
\\ \\ &= - \frac{\chi_{2}(g_{1})^2\chi_{2}(g_{2}) \chi_{2}(g_{3}) }{1 + \chi_{2}(g_{2})} a_{2}^{2} a_{1} a_{3}
- \frac{\chi_{2}(g_{1})^2\chi_{2}(g_{2})^2
\chi_{2}(g_{3}) \chi_{3}(g_{1})\chi_{3}(g_{2})^2}{1 + \chi_{2}(g_{2})}   a_{3} a_{2}^{2} a_{1}
\\ \\ &= (VIII\, a) + (VIII\, b).
\end{align*}
We next use \eqref{techni5} to show that
\begin{align*}
(I\, a) + (V)+ (VI) + (VIII\, a) &=  0,\\ \\ (I\, b) + (III)+ (IV) + (VIII\, b)&= 0.
\end{align*}
In the course of the proof of these equalities, we need \eqref{techni2} and \eqref{techni3}.
This finishes the proof of \eqref{brcomm3}.  \epf

Let $H$ be a Hopf algebra.
Then the existence of an integral for finite-dimensional braided Hopf algebras implies
\begin{lema}\label{Poincare}
Let $R = \bigoplus_{n = 0}^{N}R(n)$ be a finite-dimensional graded
braided Hopf algebra in $\ydh$ with
$R(N) \neq 0$. There exists $\lambda \in R(N)$ which is a left integral
on $R$ and such that
$$R(i) \otimes R(N - i) \to
\ku, \quad x \otimes y \mapsto \lambda(xy),$$
is a non-degenerate pairing,
for all $0 \leq i \leq N$.   In particular, $$\text{dim} R(i) = \text{dim} R(N - i).$$
\end{lema}

\pf This is essentially due to Nichols \cite[1.5]{N}. In this formulation,
one needs the existence of non-zero integrals on $R$; this follows from \cite{FMS}.
See \cite[Prop. 3.2.2]{AG} for details. \epf

\medbreak
\subsection{The infinitesimal braiding of Hopf algebras whose coradical is a Hopf subalgebra}\label{infin}

\

For the convenience of the reader,
we first recall in this Section some basic definitions from coalgebra theory.

\begin{definition}
Let $C$ be a coalgebra.
\begin{itemize}
\item  $G(C) := \{x\in C \setminus \{0\} \mid \Delta(x) = x \otimes x\}$
is the set of all group-like elements of $C$.

\medbreak
\item  If $g, h \in G(C)$, then $x \in C$ is {\it $(g,h)$-skew primitive}
if $\Delta (x) = x \otimes h + g \otimes x$. The space of all $(g,h)$-skew primitive
elements of $C$ is denoted by ${\mathcal P}(C)_{g, h}$. If $C$ is a bialgebra or a braided bialgebra,
and $g = h  =1$, then  $P(C) = {\mathcal P}(C)_{1, 1}$ is the space of {\it primitive} elements.

\medbreak
\item The {\it coradical} of $C$ is $C_0 := \sum D$,  where
$D$ runs through  all the  simple subcoalgebras of $C$;
it is the largest cosemisimple subcoalgebra of $C$.
In particular, $\ku G(C) \subseteq C_{0}$.

\medbreak
\item  $C$ is {\it pointed} if $\ku G(C) = C_{0}$.

\medbreak
\item The {\it coradical  filtration} of $C$ is  the ascending filtration
$C_{0} \subseteq C_{1} \subseteq \dots \subseteq C_{j} \subseteq C_{j + 1} \subseteq \dots,
$ defined by
$C_{j +1}:= \{x\in C \mid \Delta(x) \in C_{j}\otimes C + C\otimes C_{0}\}$.
This is a  coalgebra  filtration:
$\Delta C_j \subseteq \sum _{0\le i\le j} C_{i} \otimes C_{j-i}$; and it is exhaustive:
$C = \bigcup_{n\ge 0} C_n$.

\medbreak
\item A {\it  graded coalgebra} is a coalgebra $G$ provided with a grading
$G = \oplus_{n\ge 0} G(n)$ such that $\Delta G(j) \subseteq \sum _{0\le i\le j} G(i) \otimes G(j-i)$ for all $j \geq 0$.

\medbreak
\item A {\it coradically graded} coalgebra \cite{CM} is a graded coalgebra $G
= \oplus_{n\ge 0} G(n)$ such that its coradical filtration coincides with the standard
ascending filtration arising from the grading:
$G_{n} = \oplus_{m\le n} G(m).$ A {\it strictly graded} coalgebra \cite{Sw} is a
coradically graded coalgebra $G$ such that $G(0)$ is one-dimensional.

\medbreak
\item The graded coalgebra  associated to the coalgebra filtration of $C$ is
$\gr C = \oplus_{n\ge 0} \gr C(n)$, where
$\gr C(n) := C_{n}/C_{n-1}$, $n>0$, $\gr C(0) := C_{0}$. It
is a coradically graded coalgebra.
\end{itemize}
\end{definition}

\medbreak
We shall need  a basic technical fact on pointed coalgebras.

\begin{lema} \label{injcoalgmap} \cite[5.3.3]{M}.
A morphism of pointed coalgebras which is
injective in the first term of the coalgebra filtration, is injective.
  \qed
\end{lema}

\medbreak
Let now  $A$ be a Hopf algebra. We shall  assume in what follows
that the coradical $A_0$ is not only
a subcoalgebra but a Hopf subalgebra of $A$; this is the case if
$A$ is pointed.

\medbreak
To study the structure of $A$, we consider its coradical filtration;
because of our assumption on $A$, it is also an algebra filtration \cite{M}.
Therefore, the associated graded
coalgebra $\gr A$ is a graded Hopf algebra. Furthermore,
$H := A_{0} \simeq  \gr A(0)$
is a Hopf subalgebra of $\gr A$; and the projection
$\pi: \gr A \to \gr A(0)$ with kernel $\oplus_{n > 0} \gr A(n)$,
is a Hopf algebra map and a retraction of the inclusion. We can
then apply the general remarks of Section \ref{bosoniz}.
 Let $R$ be the algebra of coinvariants of $\pi$; $R$ is a braided Hopf algebra in
 $\ydh$ and $\gr A$ can be reconstructed  from
$R$ and $H$ as a bosonization $ \gr A \simeq R \# H$.

\medbreak
The  braided Hopf algebra $R$ is graded, since
it inherits the gradation from $\gr A$:
$R = \oplus_{n\ge 0} R(n)$, where $R(n) = \gr A(n) \cap R$.
Furthermore, $R$ is strictly graded; this means,
\begin{enumerate}
\item[{\bf (a)}.] $R(0) = \ku 1$ (hence the  coradical is trivial, {\it cf.}
\cite[Chapter 11]{Sw}).

\medbreak
\item[{\bf (b)}.] $R(1) = P(R)$ (the space of primitive elements of $R$).  \end{enumerate}

It is in general not true that a braided Hopf algebra $R$ satisfying {\bf (a)} and {\bf (b)},
also satisfies

\begin{enumerate}
\item[{\bf (c)}.] $R$ is generated as an algebra over $\ku$ by $R(1)$.  \end{enumerate}

A braided graded Hopf algebra satisfying {\bf (a)}, {\bf (b)} and {\bf (c)} is called a Nichols algebra. In the next chapter we will discuss this notion in detail.
Notice that the subalgebra $R'$ of $R$ generated by $R(1)$, a Hopf subalgebra of $R$,
is indeed a Nichols algebra.

\begin{definition}
The braiding
$$c: V\otimes V\to V\otimes V,$$
of $V := R(1) = P(R)$ is called the {\it infinitesimal braiding} of $A$.

The graded braided Hopf algebra $R$ is called the {\it diagram} of $A$.

The dimension of  $V = P(R)$ is called the {\it rank} of $A$.
\end{definition}

\section{Nichols algebras}\label{nicholsalgs}

Let $H$ be a Hopf algebra. In this Chapter, we discuss a functor
$\toba$ from the  category $\ydh$ to the category of braided Hopf
algebras in $\ydh$; given a Yetter-Drinfeld module $V$, the
braided Hopf algebra $\toba(V)$ is called the {\it Nichols
algebra} of $V$.

\medbreak The structure of a Nichols algebra appeared first in the
paper ''Bialgebras of  type one`` \cite{N} of Nichols and was
rediscovered later by several authors. In our language, a
bialgebra of type one  is just a bosonization $\toba(V) \# H$.
Hence Nichols algebras are the $H$-coinvariant elements of
bialgebras of type one, also called quantum symmetric algebras in
\cite{Ro2}. Several years after \cite{N}, Woronowicz defined
Nichols algebras in his approach to "quantum differential
calculus" \cite{W}; again, they appeared as the invariant part of
his "algebra of quantum differential forms".
 Lusztig's algebras   $\mathfrak f$ \cite{L3}, defined
by the non-degeneracy of a certain invariant bilinear form, are
Nichols algebras. In fact Nichols algebras can always be defined
by the non-degeneracy of an invariant bilinear form \cite{AG},
when $H$ is the group algebra of a finite group. The algebras
$\toba(V)$ are called bitensor algebras in \cite{Sbg}. See also
\cite{Kh1, Gr1, Gr2}.

\medbreak In a sense, Nichols algebras are similar to symmetric
algebras;  indeed, both notions coincide in the trivial braided
category of vector spaces, or more generally in any symmetric
category ({\it e. g.} in the category of super vector spaces). But
when the braiding is not a symmetry, a Nichols algebra could have
a much richer structure. We hope that this will be clarified in
the examples. On the other hand, Nichols algebras are also similar
to universal enveloping algebras. However, in spite of the efforts
of several authors, it is not clear to us how to achieve a
compact, functorial definition of a "braided Lie algebra" from a
Nichols algebra.

\medbreak We believe that Nichols algebras are very interesting
objects of  an essentially new nature.

\medbreak
\subsection{Definition of Nichols algebras}\label{DefNichols}

\

We now present one of the main notions of this survey.

\begin{definition} Let $V$ be a Yetter-Drinfeld module over $H$. A braided graded Hopf algebra $R= \oplus_{n\ge 0} R(n)$ in $\ydh$
is called a  {\it Nichols algebra} of $V$ if $\ku\simeq R(0)$ and
$V\simeq R(1)$ in $\ydh$,
and
\begin{flalign}
\label{primdegone} &P(R) = R(1),& \\
\label{gendegone} &R  \text{ is generated as an algebra by } R(1).&
\end{flalign}
The dimension of $V$ will be called the {\it rank} of $R$.
\end{definition}

\medbreak
We need some preliminaries to
show the existence and uniqueness of the Nichols algebra
of $V$ in $\ydh$.

\medbreak

Let $V$ be a Yetter-Drinfeld module over $H$.
Then the tensor algebra $T(V)=\bigoplus_{n\geq
0}T(V)(n)$  of the vector space $V$ admits a natural structure
of a Yetter-Drinfeld module, since $\ydh$ is a braided category.
It is then an algebra in $\ydh$.
There exists a unique algebra map $\Delta: T(V) \to T(V)
\underline{\otimes} T(V)$ such that $\Delta(v) = v\otimes 1 + 1\otimes
v$, for all $v\in V$. For example, if $x,y \in V$, then
$$\Delta(xy) = 1 \otimes xy + x \otimes y + x_{(-1)} \cdot y \otimes x_{(0)} + yx \otimes 1.$$
With this structure, $T(V)$ is a graded
braided Hopf algebra in
$\ydh$ with counit $\varepsilon : T(V) \to \ku$,
$\varepsilon(v) = 0$, if $v \in V$. To show the existence of the antipode, one notes that
the coradical of the coalgebra $T(V)$
is $\ku$, and uses a result of Takeuchi \cite[5.2.10]{M}.
Hence  all the braided bialgebra quotients of $T(V)$ in $\ydh$
are braided Hopf algebras in $\ydh$.

 \medbreak Let us consider the class $\mathfrak{S}$
of all $I \subset T(V)$ such that
\begin{itemize}

\item $I$ is a homogeneous ideal
generated by homogeneous elements of degree $\geq 2$,

\medbreak
\item $I$ is also a coideal, {\it i. e.} $\Delta(I) \subset I\otimes T(V) +
T(V) \otimes I$.
\end{itemize}

Note that we do {\it not} require that the ideals $I$
are Yetter-Drinfeld submodules of $T(V)$.
Let then $\widetilde{\mathfrak{S}}$ be the subset of $\mathfrak{S}$
consisting of all $I \in
\mathfrak{S}$ which are Yetter-Drinfeld submodules of $T(V)$.
The ideals
$$I(V) = \sum_{I \in \mathfrak{S}} I, \qquad
\widetilde{I}(V) = \sum_{J \in \widetilde{\mathfrak{S}}} J$$ are the largest
elements in $\mathfrak{S}$,  respectively $\widetilde{\mathfrak{S}}$.

If $I \in \mathfrak{S}$ then
$R := T(V)/I = \oplus_{n\ge 0} R(n)$
is a graded  algebra and a graded coalgebra with
$$R(0) = \ku, \qquad V\simeq R(1)  \subset P(R).$$
If actually $I \in \widetilde{\mathfrak{S}}$, then $R$
is a graded braided Hopf algebra in $\ydh$.

\medbreak

We can show now existence and uniqueness of Nichols algebras.

\begin{prop} \label{Nichols} Let $\toba(V) := T(V)/\widetilde{I}(V)$. Then the
following hold:

\begin{enumerate}
\item $V = P(\toba(V))$, hence $\toba(V)$ is a Nichols algebra of $V$.

\medbreak
\item $I(V) = \widetilde{I}(V)$.

\medbreak
\item Let $R = \oplus_{n\ge 0} R(n)$ be a graded Hopf algebra in $\ydh$
such that $R(0) = \ku 1$ and $R$ is generated as an algebra by $V:= R(1)$.
Then there exists a surjective
map of graded Hopf algebras $R \to \toba(V)$, which is
an isomorphism of Yetter-Drinfeld modules in degree 1.

\medbreak
\item  Let $R = \oplus_{n\ge 0} R(n)$ be a Nichols algebra of $V$.
Then $R\simeq \toba(V)$ as braided Hopf algebras in $\ydh$.

\medbreak
\item Let $R = \bigoplus_{n\geq 0}R(n)$ be a graded
braided Hopf algebra in $\ydh$ with
$R(0) = \ku 1$ and $R(1) = P(R)$ $= V$. Then
$\toba(V)$ is isomorphic to the subalgebra $\ku \langle V\rangle$ of $R$ generated by $V$.
\end{enumerate}
\end{prop}

\pf 1. We have to show the equality $V = P(\toba(V))$. Let us
consider the inverse image $X$ in $T(V)$ of all homogeneous
primitive elements of $\toba(V)$ in degree $n\geq 2$. Then $X$ is a graded
Yetter-Drinfeld submodule of
$T(V)$, and for all $x \in X$, $\Delta(x) \in x \otimes 1 + 1 \otimes x +
T(V) \otimes \widetilde{I}(V) + \widetilde{I}(V) \otimes T(V).$ Hence the ideal generated by $\widetilde{I}(V)$
and $X$ is in $\widetilde{\mathfrak{S}}$, and $X \subset \widetilde{I}(V)$ by the maximality of
$\widetilde{I}(V)$. Hence the image of $X$ in $\toba(V)$ is zero. This proves our claim
since the primitive elements form a graded submodule.

2. We have to show that the surjective map
$\toba(V) \to T(V)/I(V)$ is bijective.
This follows from 1. and Lemma \ref{injcoalgmap}.

3. The kernel $I$ of the canonical projection $T(V) \to R$ belongs
to $\widetilde{\mathfrak{S}}$; hence $I \subseteq \widetilde{I}(V)$.

4. follows again from  Lemma \ref{injcoalgmap}, as in 2.

5. follows from 4.
\epf

If $U$ is a braided subspace of $V \in
\ydh$, that is a subspace such that $c(U \otimes U) \subset U \otimes U$,
where $c$ is the braiding of $V$, we can define $\toba(U) := T(U)/
I(U)$ with the obvious meaning of $I(U)$. Then the description
in Proposition \ref{Nichols} also applies to $\toba(U)$.

\begin{cor} The assignment $V \mapsto \toba(V)$ is a functor
from $\ydh$ to the category of braided Hopf algebras in $\ydh$.

If $U$ is a Yetter-Drinfeld submodule of $V$, or more generally
if $U$ is a braided subspace of $V$, then  the canonical
map $\toba(U) \to \toba(V)$ is injective.
\end{cor}
\pf
If $\phi: U \to V$ is a morphism in $\ydh$, then $T(\phi) : T(U) \to T(V)$ is a morphism of braided Hopf algebras. Since $T(\phi)(I
(U))$ is a coideal
and a Yetter-Drinfeld submodule of $T(V)$, the ideal generated by $T(\phi)(\widetilde{I}(U))$ is contained in $\widetilde{I}(V)$. Hence by Proposition \ref{Nichols}, $\mathfrak{B}$ is a functor.

The second part of the claim follows from Proposition \ref{Nichols}, 5.
\epf

The duality between conditions \eqref{primdegone} and \eqref{gendegone}
in the definition of Nichols algebra, emphasized by Parts 3 and 5 of
Proposition \ref{Nichols}, is explicitly stated in the following

\begin{lema}\label{dualityprinciple} Let  $R = \oplus_{n\ge 0} R(n)$
be a graded braided Hopf algebra
in $\ydh$; suppose the homogeneous components
are finite-dimensional   and  $R(0) = \ku 1$.
Let $S = \oplus_{n \ge 0} R(n)^*$ be the graded dual of $R$.
Then $R(1) = P(R)$ if and only if $S$ is generated as an algebra by $S(1).$
\end{lema}
\pf See for instance \cite[Lemma 5.5]{AS2}. \epf

\begin{exa}\label{5.11}
Let $F$ be a field of positive characteristic $p$. Let $S$ be the
(usual) Hopf algebra $F[x]/\langle x^{p^{2}}\rangle$ with $x\in P(S)$.
Then $x^p \in P(S)$.
Hence $S$  satisfies \eqref{gendegone} but not \eqref{primdegone}.
\end{exa}

\medbreak
\begin{exa}\label{5.12}
Let $S = \ku[X] = \oplus_{n\ge 0} S(n)$
be a polynomial algebra in one variable. We consider $S$ as a
braided Hopf  algebra in $\ydh$, where $H = \ku \Gamma$, $\Gamma$
an infinite cyclic group with generator $g$, with action, coaction
and comultiplication given by
$$
\delta(X^n)= g^n \otimes X^n, \quad g \cdot X = qX, \quad\Delta(X) = X\otimes 1
+ 1 \otimes X.
$$
Here $q \in \ku$ is a root of 1 of order $N$.
That is, $S$ is a so-called quantum line. Then $S$ satisfies \eqref{gendegone}
but not \eqref{primdegone} since $X^N$ is also primitive. Hence the
graded dual $R = S^d = \oplus_{n\ge 0} S(n)^*$ is a braided Hopf algebra
satisfying \eqref{primdegone} but not \eqref{gendegone}. \end{exa}

However, in characteristic 0 we do not know any finite dimensional example
of a braided Hopf algebra satisfying \eqref{primdegone} but not \eqref{gendegone}.

\begin{conjecture}\label{conjecturedegree1} \cite[Conjecture 1.4]{AS2} Any {\it finite dimensional}
braided Hopf
algebra in $\ydh$ satisfying \eqref{primdegone} also satisfies \eqref{gendegone}.
(Recall that the base field $\ku$ has characteristic zero.)
\end{conjecture}

\medbreak
The compact description of $\toba(V)$ in Lemma \ref{Nichols}
shows that it depends only on the algebra and coalgebra structure of $T(V)$.
Since the comultiplication of the tensor algebra was defined using the "twisted" multiplication of
$T(V)\underline{\otimes} T(V)$, we see that
$\toba(V)$ {\it depends as an algebra and coalgebra only
on the braiding of $V$}. The explicit formula for the comultiplication
of $T(V)$ leads to the following alternative description of $\toba(V)$.

\medbreak
\subsection{Skew-derivations and bilinear forms} \label{skew-der-bil}

\

We want to describe two important techniques to prove identities in Nichols algebras
even without knowing the defining relations explicitly.

The first technique was
introduced by Nichols \cite[3.3]{N} to deal
with $\toba(V)$ over group algebras $\ku \Gamma$ using
skew-derivations.
Let $V \in \ydg$ be of finite dimension $\theta$.
We choose a
basis $x_{i} \in V_{g_{i}}$ with $g_{i} \in \Gamma, 1 \leq i \leq \theta,$ of
$\Gamma$-homogeneous elements. Let $I \in \mathfrak{S}$ and $R = T(V)/I$ (see Section \ref{DefNichols}). Then $R$
is a graded Hopf algebra in $\ydg$ with $R(0) = \ku 1$ and $R(1) = V$. For
all $1 \leq i \leq \theta$ let $\sigma_{i} : R \to R$ be the algebra
automorphism given by the action of $g_{i}$.

Recall that if $\sigma: R \to R$ is an algebra automorphism, an
$(\id, \sigma)$-derivation $D: R\to R$ is a $\ku$-linear map such that
$$
D(xy) = xD(y) + D(x) \sigma(y), \qquad \text{for all } x,y \in R.
$$

\begin{prop}\label{derivation}
1) For all $1 \leq i \leq \theta$, there exists a uniquely determined
$(id,\sigma_{i})$-derivation $D_{i} : R \to R$ with $D_{i}(x_{j}) =
\delta_{i,j}$ (Kronecker $\delta$) for all $j$.

2) $I = I(V)$, that is $R = \toba(V)$, if and only if $\bigcap_{i =
1}^{\theta} \text{ker}(D_{i}) = \ku 1.$
\end{prop}

\pf See for example \cite[2.4]{MiS}.
\epf

Let us illustrate this Proposition in a very simple case.

\begin{exa}\label{ejemplo} Let $V$ be as above and assume that $g_i \cdot x_i = q_i x_i$,
for some $q_i\in \ku^{\times}$, $1\le i \le \theta$. Then for any $n\in \N$,

(a). $D_i(x_i^n) = (n)_{q_i}x_i^{n-1}$.

(b). $x_i^n \neq 0$ if and only if $(n)_{q_i}! \neq 0$.
\end{exa}

\pf (a) follows by induction on $n$ since $D_i$ is a skew-derivation;
(b) follows from (a) and Proposition \ref{derivation}, since $D_j$
vanishes on any power of $x_i$, for $j\neq i$.
\epf

The second technique was used by Lusztig \cite{L3} to prove very deep results
about quantum enveloping algebras using a canonical bilinear form.

Let $(V, c)$ be a braided vector space of diagonal type
as in \eqref{br-abfingrt} and assume that $q_{ij} = q_{ji}$
for all $i, j$. Let $\Gamma$ be the free abelian group of rank $\theta$
with basis $g_{1}, \dots, g_{\theta}$.
We define characters $\chi_{1}, \dots, \chi_{\theta}$ of $\Gamma$ by
$$
\chi_{i} (g_{j}) = q_{ji}, \qquad 1\le i,j\le \theta.
$$
We consider $V$ as a Yetter-Drinfeld module over $\ku \Gamma$ by defining
$x_{i} \in V^{\chi_{i}}_{g_{i}}$, for all $i$.

\begin{prop}\label{radical-nichols}
 Let $B_1, \dots, B_{\theta}$ be non-zero elements in
$\ku$. There is a unique  bilinear form
$(\, \vert \, ): T(V) \times T(V) \to \ku$ such that $(1\vert 1) = 1$ and
\begin{flalign}\label{bil1} &(x_j \vert x_j) = \delta_{ij} B_{i}, \qquad \text{for all }
i, j;&\\
\label{bil2} &(x \vert yy') = (x_{(1)} \vert y)(x_{(2)} \vert y'), \qquad \text{for all }
x, y, y'\in T(V);& \\
\label{bil3}
&(xx' \vert y) = (x \vert y_{(1)})(x' \vert y_{(2)}), \qquad \text{for all }
x, x', y\in T(V).&
\end{flalign}
This form is symmetric and also satisfies
\begin{equation}
(x \vert y) = 0, \qquad \text{for all }x\in T(V)_g, y\in T(V)_h, g\neq h \in \Gamma.
\end{equation}
The homogeneous components of $T(V)$ with respect to its usual $\N$-grading are
also orthogonal with respect to $(\, \vert \, )$.

The quotient $T(V) / I(V)$, where  $I(V) =
\{x\in T(V): (x \vert y) = 0 \forall y\in T(V)\}$ is the radical of the form,
is canonically isomorphic to the Nichols algebra of $V$. Thus, $(\, \vert \, )$
induces a non-degenerate bilinear form on $\toba(V)$, denoted by the same name.
\end{prop}

\pf The existence and uniqueness of the form, and the claims about symmetry and
orthogonality, are proved exactly as in \cite[1.2.3]{L3}. It follows from
the properties of the form that $I(V)$ is a Hopf ideal. We now check that
$T(V) / I(V)$ is the Nichols algebra of $V$; it is enough to verify that
the primitive elements of $T(V) / I(V)$ are in $V$. Let $x$ be a primitive element
in $T(V) / I(V)$, homogeneous of degree $n\ge 2$. Then $(x\vert yy') =0$
for all $y$, $y'$ homogeneous of degrees $m, m'\ge 1$ with $m + m' = n$; thus $x = 0$.
\epf

A generalization of the preceding result, valid for any
finite dimensional Yetter-Drinfeld module
over any group, can be found in \cite[3.2.17]{AG}.

\medbreak
The Proposition shows that Lusztig's algebra $\bf f$ \cite[Chapter 1]{L3}
is the Nichols algebra of $V$ over the field of rational
functions $\kv$, with $q_{ij} = v^{i\cdot j}$ if
$I =\{1, \dots, \theta\}$ and $(I, \cdot)$ a Cartan datum.
In particular, we can take a generalized symmetrizable
Cartan matrix $(a_{ij})$, $0 < d_i \in \N$ for all $i$
with $d_ia_{ij} = d_ja_{ji}$ for all $i,j$ and define $i\cdot j := d_ia_{ij}$.

\medbreak
\subsection{The braid group} \label{brgroup}

\

Let us recall that the braid group
$\mathbb B_{n}$ is presented by generators $\sigma_{1}, \dots,
\sigma_{n-1}$ with relations
\begin{align*}
\sigma_{i}\sigma_{i+1}\sigma_{i} &= \sigma_{i+1}\sigma_{i}\sigma_{i+1},
\quad 1\le i \le n-2,\\
\sigma_{i}\sigma_{j} &= \sigma_{j}\sigma_{i}, ,
\quad 1\le i, j \le n-2, \ \vert i-j \vert > 1.
\end{align*}
Here are some basic well-known facts about the braid group.

\medbreak
There is a natural projection $\pi: \mathbb B_{n} \to \mathbb S_{n}$ sending
$\sigma_{i}$ to the transposition $\tau_{i} :=(i, i+1)$ for all $i$.
 The projection $\pi$ admits a set-theoretical section
$s:\mathbb S_{n} \to \mathbb B_{n}$ determined by
\begin{align*}
s(\tau_{i}) &= \sigma_{i}, \quad 1\le i \le n-1,\\
s(\tau\omega) &= s(\tau)s(\omega), \quad \text{if } \ell(\tau\omega) = \ell(\tau) + \ell(\omega).
\end{align*}
Here $\ell$ denotes the length of an element of $\mathbb S_{n}$
with respect to the set of generators $\tau_{1}, \dots,
\tau_{n-1}$. The map $s$ is called the Matsumoto section.
In other words, if $\omega = \tau_{i_1} \dots \tau_{i_M}$
is a reduced expression of  $\omega \in \mathbb S_{n}$,
then $s(\omega) = \sigma_{i_1} \dots \sigma_{i_M}$.

\medbreak Let $q\in \ku$, $q\neq 0$. The quotient of the group
algebra $\ku(\mathbb B_{n})$ by the two-sided ideal generated by
the relations
$$ (\sigma_{i} - q)(\sigma_{i} + 1), \quad 1\le i \le n-1,$$
is the so-called {\it Hecke algebra} of type $A_{n}$, denoted by
$\hecke$.

\medbreak
Using the section $s$, the following
distinguished elements of the group algebra $\ku\mathbb B_{n}$ are defined:

$$\mathfrak S_{n} := \sum_{\sigma \in \mathbb S_{n}} s(\sigma),
\qquad \mathfrak S_{i,j} := \sum_{\sigma \in X_{i,j}} s(\sigma);$$
here $X_{i,j} \subset \mathbb S_{n}$ is the set of all $(i,j)$-shuffles.
The element $\mathfrak S_{n}$ is called the quantum symmetrizer.

\medbreak

Given a braided vector space $(V, c)$, there are representations
of the braid groups $\rho_n: \mathbb B_{n} \to \Aut (V^{\otimes n})$
for any $n \ge 0$, given  by
$$
\rho_n(\sigma_{i}) = \id \otimes \dots \otimes  \id \otimes
c \otimes  \id \otimes \dots  \otimes \id,
$$
where $c$ acts in the tensor product of the $i$ and $i+1$ copies of $V$.
By abuse of notation, we shall denote by $\mathfrak S_{n}$, $\mathfrak S_{i,j}$
also the corresponding endomorphisms $\rho(\mathfrak S_{n})$, $\rho(\mathfrak S_{i,j})$ of $V^{\otimes n} = T^{n}(V)$.

\medbreak

If $C = \bigoplus_{n\geq0} C(n)$ is a graded coalgebra with comultiplication $\Delta$,
we denote by
$\Delta_{i,j}: C(i+j) \to C(i) \otimes C(j)$, $i,j \geq0$, the $(i,j)$-graded
component of the map $\Delta$.

\begin{prop}\label{nic-alg}
Let $V \in \ydh$. Then
\begin{align}
  \Delta_{i,j} =  \mathfrak S_{i,j},\\
\toba(V) = \bigoplus_{n\geq0} T^{n}(V)/ker(\mathfrak{S}_{n}). \label{generalB}
\end{align}
\end{prop}
\pf See for instance \cite{Sbg}. \epf

This description of the relation of $\toba(V)$ does not mean that the relations are known. In general it is very hard to compute the kernels of the maps $\mathfrak{S}_{n}$ in concrete terms.
For any braided vector space $(V,c)$, we may define $\toba(V)$ by \eqref{generalB}.

\medbreak
Using the action of the braid group, $\toba(V)$ can also be described as
a subalgebra of the quantum shuffle algebra \cite{N, Ro1, Ro2, Sbg}.

\medbreak
\subsection{Invariance under twisting}\label{Twisting}

\

Twisting  is a method to construct new Hopf algebras by "deforming" the comultiplication; originally due to Drinfeld \cite{Dr2}, it was adapted to Hopf algebras in \cite{Re}.

\medbreak
Let $A$ be a Hopf algebra and $F
\in A\otimes A$ be an invertible element. Let $\Delta_{F} := F \Delta F^{-1}: A
\to A\otimes A$; it is again an algebra map. If
\begin{align} \label{twconduno}
(1\otimes F) (\id \otimes \Delta)(F) &=
(F\otimes 1) (\Delta \otimes \id)(F),
\\\label{twconddos} (\id \otimes \varepsilon)(F) &=
(\varepsilon \otimes \id)(F) = 1,
\end{align}
then $A_{F}$ (the same algebra,
but with comultiplication $\Delta_{F}$) is again a Hopf algebra.
We shall say  that $A_{F}$ is obtained
from $A$ via twisting by $F$;  $F$ is a cocycle in a suitable sense.

\medbreak
There is a dual version of the twisting operation, which amounts to a
twist of the multiplication \cite{DT}.
Let $A$ be a Hopf algebra and let $\sigma: A\times A \to \ku$ be an invertible
2-cocycle$^1$\footnote{$^1$ Here "invertible" means that the associated linear map
$\sigma: A\otimes A \to \ku$ is invertible with respect to the convolution product.},
that is
\begin{align*} \sigma(x_{(1)}, y_{(1)}) \sigma(x_{(2)} y_{(2)}, z) &=
\sigma(y_{(1)}, z_{(1)}) \sigma(x, y_{(2)}z_{(2)}), \\
\sigma(x, 1) &= \sigma(1, x) = \varepsilon(x),
\end{align*}
for all $x,y,z\in A$. Then $A_{\sigma}$ -- the same $A$ but with the multiplication
$\cdot_{\sigma}$ below -- is again a Hopf algebra, where
$$
x\cdot_{\sigma}y = \sigma(x_{(1)}, y_{(1)}) x_{(2)} y_{(2)} \sigma^{-1}(x_{(3)}, y_{(3)}).
$$

For details, see for instance \cite[10.2.3 and 10.2.4]{KS}.

\bigbreak
Assume now that $H$ is a Hopf algebra, $R$ is a braided Hopf algebra in $\ydh$, and $A = R\# H$.
Let $\pi: A \to H$ and $\iota: H \to A$ be the canonical projection and injection.
Let $\sigma: H\times H \to \ku$ be an invertible 2-cocycle, and define
$\sigma_{\pi}: A\times A \to \ku$ by
$$
\sigma_{\pi} := \sigma (\pi \otimes \pi);
$$
$\sigma_{\pi}$ is an invertible 2-cocycle, with inverse $(\sigma^{-1})_{\pi}$.
The maps $\pi: A_{\sigma_{\pi}} \to H_{\sigma}$, $\iota: H_{\sigma} \to A_{\sigma_{\pi}}$
are still Hopf algebra maps. Because the comultiplication is not changed, the
space of coinvariants of $\pi$ is $R$; this is a subalgebra of $A_{\sigma_{\pi}}$ that we denote $R_{\sigma}$; the multiplication in $R_{\sigma}$ is given by
\begin{equation}
\label{twistingalgebra}
x._{\sigma}y = \sigma(x_{(-1)}, y_{(-1)}) x_{(0)} y_{(0)}, \qquad x, y \in R = R_{\sigma}.
\end{equation}
Equation \eqref{twistingalgebra} follows easily using \eqref{smash1}. Clearly, $R_{\sigma}$ is a Yetter-Drinfeld Hopf algebra in $\ydhsi$. The  coaction of $H_{\sigma}$ on $R_{\sigma}$ is the same as the coaction of $H$ on $R$, since the comultiplication was not altered.
The explicit formula for the action of $H_{\sigma}$ on $R_{\sigma}$ can be written down; we shall do this only in the setting we are interested in.

\medbreak
Let $H = \ku \Gamma$ be a group algebra; an invertible 2-cocycle
$\sigma: H \times H \to \ku$ is uniquely determined by its restriction
$\sigma: \Gamma \times \Gamma \to \ku^{\times}$, a group
2-cocycle with respect to the trivial action.

\medbreak
\begin{lema} Let $\Gamma$ be an abelian group and let $R$
be a braided Hopf algebra in $\ydg$. Let
$\sigma: \Gamma \times \Gamma \to \ku^{\times}$ be a 2-cocycle.
Let $S$ be the subalgebra of $R$ generated by $P(R)$.
In the case $y\in S_h^{\eta}$, for some $h \in \Gamma$ and $\eta \in \VGamma$,
the action of $H = H_{\sigma}$ on $R_{\sigma}$ is
\begin{equation}
\label{twistingaction}
g\rightharpoonup_{\sigma}y = \sigma(g,h)\sigma^{-1}(h, g)\eta(g) y, \qquad g \in \Gamma.
\end{equation}
Hence, the braiding $c_{\sigma}$ in  $R_{\sigma}$ is given in this case by
\begin{equation} \label{twistingbraiding}
c_{\sigma} (x \otimes y) = \sigma(g,h)\sigma^{-1}(h, g)\eta(g) \,  y \otimes x, \qquad
x\in R_{g}, g \in \Gamma.
\end{equation}
Therefore, for such $x$ and $y$, we have
\begin{equation} \label{twistingad}
 [x, y]_{c_{\sigma}} = \sigma(g,h) [x, y]_{c}.
\end{equation} \end{lema}

\pf To prove \eqref{twistingaction}, it is enough to assume $y\in P(R)_h^{\eta}$.

Let $A = R\# H$; in $A_{\sigma_{\pi}}$, we have
\begin{align*} g._{\sigma} y &=  \sigma(g, \pi(y))\, g \, \sigma^{-1}(g, 1)
+ \sigma(g, h) \, gy\, \sigma^{-1}(g, 1) + \sigma(g, h)\, gh \, \sigma^{-1}(g, \pi(y))
\\&=   \sigma(g, h) gy;\\
y._{\sigma} g &=  \sigma( \pi(y), g)\, g \, \sigma^{-1}(1, g)
+ \sigma(h, g)\, yg \, \sigma^{-1}(1,g) + \sigma(h, g) \,hg \, \sigma^{-1}(\pi(y), g)
\\&=   \sigma(h, g) yg;
\end{align*}
hence
$$ g._{\sigma} y  =   \sigma(g, h) gy = \sigma(g, h)  \eta(g) yg = \sigma(g, h)
\sigma^{-1}(h, g)\eta(g) y._{\sigma} g, $$
which is equivalent to \eqref{twistingaction}. Now \eqref{twistingbraiding} follows
at once, and \eqref{twistingad} follows from \eqref{twistingalgebra}
and \eqref{twistingbraiding}:
\begin{equation*}
 [x, y]_{c_{\sigma}} =  x._{\sigma}y - ._{\sigma} c_{\sigma} (x \otimes y)
 = \sigma(g, h) xy  -\sigma(g,h)\sigma^{-1}(h, g)\eta(g) \,\sigma(h, g)yx = \sigma(g,h) [x, y]_{c}.  \end{equation*}
\epf

The proof of the following Lemma is clear, since the comultiplication of a
Hopf algebra is not changed by twisting.

\begin{lema}\label{twistnichols} Let $H$ be a  Hopf algebra and let
$R$ be a braided Hopf algebra in $\ydh$. Let $\sigma: H\times H \to \ku$ be an
invertible 2-cocycle. If $R =\oplus_{n\ge 0} R(n)$ is a braided graded Hopf
algebra in $\ydh$, then
$R_{\sigma}$ is a braided graded Hopf algebra in $\ydhsi$
with $R(n) = R_{\sigma}(n)$ as vector spaces for all $n\ge 0$.
Also $R$ is a Nichols algebra if and only if $R_{\sigma}$
is a Nichols algebra in $\ydhsi$. \qed
\end{lema}

\section{Types of Nichols algebras}\label{examnichols}

We now discuss several  examples
of Nichols algebras.  We are interested in explicit presentations,
{\it e. g.} by generators and relations, of $\toba(V)$, for braided
vector spaces in suitable classes, for instance, those of group type.
We would also like to determine when $\toba(V)$ has finite dimension, or polynomial growth.

\subsection{Symmetries and braidings of Hecke type}

\medbreak
We begin with the simplest class of braided vector spaces.

\begin{exa} Let $\tau: V\otimes V \to V\otimes V$ be the usual transposition;
the braided vector space $(V, \tau)$ can be realized as a Yetter-Drinfeld
module over any Hopf algebra $H$, with trivial action and coaction. Then $\toba(V) = \Sym(V)$, the symmetric algebra of $V$.

\medbreak
The braided vector space $(V, -\tau)$, which can be realized {\it e. g.}
in $\ydzeta$, has $\toba(V) = \Lambda(V)$, the exterior algebra of $V$.
\end{exa}

\begin{exa} Let $V = \oplus_{i\in \Z/2} V(i)$
be a super vector space and let
$c: V\otimes V \to V\otimes V$ be the supersymmetry:
$$
c(v\otimes w) = (-1)^{i.j} w\otimes v\qquad v \in V(i),  \ w\in V(j).
$$
Clearly, $V$ can  be realized as a Yetter-Drinfeld  module over $\Z/2$.
Then $\toba(V) \simeq \Sym(V(0)) \otimes \Lambda(V(1))$,
the super-symmetric algebra of $V$.
\end{exa}

\medbreak
The simple form of $\toba(V)$ in these examples
can be explained in the following context.

\begin{definition} We say that a braided vector space $(V, c)$ is of
{\it Hecke-type} with label $q\in \ku$, $q\neq 0$,  if
$$
(c-q)(c+ 1) = 0.
$$
In this case, the representation of the braid group
$\rho_n: \mathbb B_{n} \to \Aut (V^{\otimes n})$ factorizes through
the Hecke algebra $\hecke$, for all $n\ge 0$; {\it cf.} Section \ref{brgroup}.

\medbreak
If $q = 1$, one says that $c$ is a {\it symmetry}. Then
$\rho_n$ factorizes through  the symmetric group $\mathbb S_{n}$,
for all $n\ge 0$.
The categorical version of
symmetries is that of symmetric categories, see Section \ref{brcat}.
\end{definition}

\begin{prop} \label{hecke} Let $(V, c)$ be a braided vector space  of
Hecke-type with label $q$, which is either 1 or not a root of 1.
Then $\toba(V)$ is a quadratic algebra; that is, the ideal $I(V)$
is generated by $I(V)(2) = \ker \mathfrak S_{2}$.

Moreover, $\toba(V)$ is a Koszul algebra and its Koszul dual is
the Nichols algebra $\toba(V^*)$ corresponding to the braided vector
space $(V^*, q^{-1}c^{t})$.
\end{prop}

A nice exposition on Koszul algebras is \cite[Chapter 2]{BGS}.

\pf The argument for the first claim is taken from \cite[Prop. 3.3.1]{AA}.
The image of the quantum symmetrizer $\mathfrak S_{n}$ in the
Hecke algebra $\hecke$ is $[n]_q! M_{\varepsilon}$, where $M_{\varepsilon}$
satisfies the following properties:
$$
M_{\varepsilon}^2 = M_{\varepsilon}, \quad M_{\varepsilon} c_{i} =
c_{i} M_{\varepsilon} = q M_{\varepsilon}, \qquad 1\le i \le n-1.
$$
See for instance \cite{HKW}. Now, we have to show that
$\ker \mathfrak S_{n} = T^n(V) \cap I$,
where $I$ is the ideal generated by $\ker \mathfrak S_{2}
= \ker (c+1) = \im (c-q)$; but clearly
$T^n(V) \cap I = \sum_{i} I^{n,i}$,  where
$$I^{n,i} =
T^{i -1}(V) \otimes \im (c-q) \otimes T^{n -i -1}(V) = \im (c_{i} - q).
$$
It follows  that $T^n(V) \cap I \subseteq \ker \mathfrak S_{n}$,
a fact that we already know from the general theory. But moreover,
$T^n(V) \cap I$ is a $\hecke$-submodule of $T^n(V)$ since
$$
c_{j}(c_{i} - q) = (c_{j} - q)(c_{i} - q) + q (c_{i} - q).
$$
This computation also shows that the action of $\hecke$ on the quotient
module $T^n(V) / T^n(V) \cap I$ is via the character that sends $\sigma_{i}$
to $q$; hence $M_{\varepsilon}$ acts on $T^n(V) / T^n(V) \cap I$
by an automorphism, and {\it a fortiori} $T^n(V) \cap I
\supseteq \ker \mathfrak S_{n}$.
Having shown the first claim, the second claim is essentially a result
from \cite{Gu, Wa}; see also the exposition in \cite[Sections 3.3 and 3.4]{AA}.
\epf

\begin{exa} Let $q\in \ku^{\times}$, $q$ is not a root of 1.
The braided vector space $(V, q\tau)$ can be realized in $\ydzeta$. It can be shown that
$\toba(V) = T(V)$, the tensor algebra of $V$, for all $q$ in an open set. Problem: Determine this open set.
\end{exa}

It would be interesting to know whether other conditions on the minimal polynomial of a braiding
have consequences on the structure of the corresponding Nichols algebra. The first candidate should be a braiding of BMW-type.

\medbreak
\subsection{Braidings of diagonal type}

\

In this Section,  $(V, c)$ denotes a finite dimensional braided
vector space of diagonal type; that is, $V$ has a basis $x_{1},
\dots, x_{\theta}$ such that \eqref{br-abfingrt} holds for some
non-zero $q_{ij}$ in $\ku$. Our first goal is to determine
polynomial relations on the generators $x_{1}, \dots, x_{\theta}$
that should hold in $\toba(V)$. We look at polynomial expressions
in these generators which are homogeneous of degree $\ge 2$,  and
give rise to primitive elements in any  braided Hopf algebra
containing $V$ inside its primitive elements. For related
material, see \cite{Kh1}.

\begin{lema}\label{primitivos}
Let $R$ be a braided Hopf algebra in $\ydh$, for some Hopf algebra $H$, such that
$V \subseteq P(R)$ as braided vector spaces.
\begin{itemize}
\item [(a).] If $q_{ii}$ is a root of 1 of order $N > 1$ for some $i\in \{1,
\dots, \theta\}$, then $x_{i}^N \in P(R)$.

\medbreak
\item [(b).] Let $1 \leq i,j \leq \theta, i \neq j,$ such that $q_{ij}q_{ji} = q_{ii}^r$, where $0 \leq - r < \ord q_{ii}$
(which could be infinite). Then $(\ad_{c} x_{i})^{1-r} (x_{j})$ is primitive in $R$.
\end{itemize}
\end{lema}

\pf (a) and (b) are consequences of the quantum binomial formula, see {\it e. g.}
\cite[Appendix]{AS2} for (b).  \epf

We apply these first remarks to $\toba(V)$ and see how conditions
on the Nichols algebra induce conditions on the braiding.

\begin{lema}\label{rossito}  Let $R = \toba(V)$.
\begin{itemize}
\item [(a).] If $q_{ii}$ is a root of 1 of order $N > 1$ then
$x_{i}^N = 0$. In particular, if $\toba (V)$ is an integral domain,
then $q_{hh} = 1$ or it is not a root of 1, for all $h$.

\medbreak
\item [(b).] If $i\neq j$, then
$(\ad_{c}x_{i})^{r} (x_{j}) = 0$ if and only if
$(r)!_{q_{ii}} \prod_{0\le k \le r - 1}
\left(1 - q_{ii}^{k} q_{ij}q_{ji} \right) = 0$.

\medbreak
\item [(c).] If $i\neq j$ and $q_{ij}q_{ji} = q_{ii}^r$, 
where $0 \leq  - r < \ord q_{ii}$
(which could be infinite), then $$(\ad_{c} x_{i})^{1-r} (x_{j}) = 0.$$

\medbreak
\item [(d).] If $\toba(V)$ has finite Gelfand-Kirillov dimension,
then for all $i\neq j$, there exists $r_{ij} > 0$ such
that $(\ad_{c}x_{i})^{r_{ij}} (x_{j}) = 0$.

\end{itemize}
\end{lema}

\pf  Parts (a) and (c) follow from Lemma \ref{primitivos};
part (a) is also a special case of Example \ref{ejemplo};
and part (c) also follows from (b).
Part (b) is stated in \cite[Lemma 14]{Ro2}. It can be shown
using the skew-derivations $D_j$ of Section \ref{skew-der-bil}.
Indeed, we first claim that $D_j\left((\ad_{c}x_{i})^{r} (x_{j}) \right)
= \prod_{0\le k \le r - 1} \left(1 - q_{ii}^{k} q_{ij}q_{ji} \right)\, x_{i}^{r}$.
We set $z_r = (\ad_{c}x_{i})^{r} (x_{j})$ and compute
\begin{align*} D_j\left(\ad_{c}x_{i} (z_r) \right) &= D_j(x_{i}z_r - (g_i\cdot z_r) x_i)
\\
&= D_j(x_{i}z_r - q_{ii}^{r} q_{ij} z_r x_i)
\\
& = x_{i} D_j(z_r) - q_{ii}^{r} q_{ij}q_{ji} D_j(z_r) x_i
\end{align*}
and the claim follows by induction. Thus, by Example \ref{ejemplo},
$D_j\left((\ad_{c}x_{i})^{r} (x_{j}) \right)
= 0$ if and only if
$(r)!_{q_{ii}} \prod_{0\le k \le r - 1}
\left(1 - q_{ii}^{k} q_{ij}q_{ji} \right) = 0$.
We next claim that
$D_i\left((\ad_{c}x_{i})^{r} (x_{j}) \right)
= 0$.
We  compute
\begin{align*} D_i\left(\ad_{c}x_{i} (z_r) \right) &= D_i(x_{i}z_r - (g_i\cdot z_r) x_i)
\\
& = x_{i} D_i(z_r) + g_i\cdot z_r - g_i\cdot z_r - D_i(g_i\cdot z_r) \, g_i\cdot x_i
\end{align*}
and the claim follows by induction. Finally, it is clear that
$D_{\ell}\left((\ad_{c}x_{i})^{r} (x_{j}) \right) = 0$, for all $\ell \neq i,j$.
Part (b) follows then from Proposition \ref{derivation}.

Part(d) is an important result of Rosso \cite[Lemma 20]{Ro2}. \epf

We now discuss how the twisting operation, {\it cf.} Section \ref{Twisting},
affects Nichols algebras of diagonal type.

\begin{definition} We shall say that two braided vector spaces $(V, c)$ and $(W, d)$ of diagonal type, with matrices $(q_{ij})$ and $({\widehat q}_{ij})$, are {\it twist-equivalent}
if $\dim V = \dim W$ and, for all  $i, j$, $q_{ii} = {\widehat q}_{ii}$  and
\begin{equation}\label{twistequiv}
q_{ij}q_{ji} = {\widehat q}_{ij}{\widehat q}_{ji}.  \end{equation} \end{definition}

\begin{prop} Let $(V, c)$ and $(W, d)$ be two twist-equivalent braided vector spaces  of diagonal type, with matrices $(q_{ij})$ and $({\widehat q}_{ij})$; say with respect to basis $x_{1}, \dots x_{\theta}$, resp. ${\widehat x}_{1}, \dots {\widehat x}_{\theta}$. Then there exists a linear isomorphism $\psi: \toba(V) \to \toba (W)$ such that
\begin{equation}
\psi(x_{i}) = {\widehat x}_{i}, \qquad 1\le i \le \theta.
  \end{equation} \end{prop}

\pf Let $\Gamma$ be the free abelian group of rank $\theta$, with basis $g_{1}, \dots, g_{\theta}$. We define characters $\chi_{1}, \dots, \chi_{\theta}$, ${\widehat \chi}_{1}, \dots, {\widehat \chi}_{\theta}$ of $\Gamma$ by
$$
\chi_{i} (g_{j}) = q_{ji}, \qquad {\widehat \chi}_{i} (g_{j}) = {\widehat q}_{ji}, \qquad 1\le i,j\le \theta.
$$
We consider $V$, $W$ as Yetter-Drinfeld modules over $\Gamma$ by declaring
$x_{i} \in V^{\chi_{i}}_{g_{i}}$, ${\widehat x}_{i} \in V^{{\widehat \chi}_{i}}_{g_{i}}$.
Hence,
$\toba(V), \toba (W)$ are braided Hopf algebras in $\ydg$.

Let $\sigma: \Gamma \times \Gamma \to \ku^{\times}$ be the unique bilinear
form such that
\begin{equation}
\sigma(g_{i}, g_{j}) = \begin{cases} &{\widehat q}_{ij}  q^{-1}_{ij}, \qquad i\le j, \\
&1, \qquad \qquad i > j;\end{cases}  \end{equation}
it is a group cocycle. We claim that $\varphi: W \to \toba(V)_{\sigma} (1)$,
$\varphi({\widehat x}_{i}) = x_{i}$, $1\le i\le \theta$,
is an isomorphism in $\ydg$. It clearly preserves the coaction; for the action, we
assume $i\le j$ and compute
\begin{align*}
g_{j}\cdot_{\sigma}x_{i} &= \sigma(g_{j}, g_{i})\sigma^{-1}(g_{i}, g_{j})\chi_{i}(g_{j}) x_{i}
\\
&= ({\widehat q}_{ij})^{-1}  q_{ij} q_{ji} x_{i} = {\widehat q}_{ji} x_{i},
\end{align*}
and also
\begin{align*}
g_{i}\cdot_{\sigma}x_{j} &= \sigma(g_{i}, g_{j})\sigma^{-1}(g_{j}, g_{i})\chi_{j}(g_{i}) x_{j}
\\
&= {\widehat q}_{ij}  q^{-1}_{ij}   q_{ij} x_{j} = {\widehat q}_{ij} x_{j},
\end{align*}
where we have used \eqref{twistingaction} and the hypothesis \eqref{twistequiv}.
This proves the claim. By Proposition \ref{Nichols}, $\varphi$ extends to an isomorphism
$\varphi: \toba (W) \to \toba(V)_{\sigma}$; $\psi = \varphi^{-1}$  is the map
we are looking for. \epf

\begin{rmk} (i). The map $\psi$ defined in the proof is much more than just linear; by \eqref{twistingalgebra} and \eqref{twistingad}, we have for all $g,h \in \Gamma$,

\begin{align}
\label{psimult}
\psi(xy) &= \sigma^{-1}(g, h)\psi(x)\psi(y), \qquad x \in \toba (V)_{g},
\quad y \in \toba (V)_{h}; \\ \label{psiad}
 \psi([x, y]_{c}) &=  \sigma^{-1}(g, h) [\psi(x), \psi(y)]_{d},
\qquad x \in \toba(V)^{\chi}_{g}, \quad y \in \toba (V)_{h}^{\eta}.
\end{align}

(ii). A braided vector space $(V, c)$ of diagonal type, with matrix $(q_{ij})$,
is twist-equivalent to  $(W, d)$, with a {\it symmetric} matrix $({\widehat q}_{ij})$.
\end{rmk}

Twisting is a very important tool. For many problems, twisting allows to reduce to
the case when the diagonal braiding is symmetric; then the theory of
quantum groups can be applied.

\medbreak
\subsection{Braidings of diagonal type but not Cartan}

\

In the next Chapter, we shall concentrate on braidings of Cartan type.
There are a few examples  of Nichols algebras $\toba(V)$ of
finite group type and rank 2, which are {\it not}
of Cartan type, but where we know that the dimension
is finite.
We now list the examples we know, following \cite{N, Gn4}.
The braided vector space is necessarily of diagonal
type; we shall give the matrix $Q$ of the braiding, the constraints
on their entries and the dimension $d$ of $\toba(V)$.
Below, $\omega$, resp. $ \zeta$, denotes an arbitrary primitive third root
of 1, resp. different from 1.

{\allowdisplaybreaks
\begin{flalign}\label{mg1}  &\begin{pmatrix}
q_{11} & q_{12} \\ q_{21} & -1
\end{pmatrix}; \quad q_{11}^{-1} = q_{12}q_{21} \neq  1;
 & d &= 4\ord (q_{12}q_{21}).& \\ \notag \\
\label{mg2} &\begin{pmatrix}
q_{11} & q_{12} \\ q_{21} & \omega
\end{pmatrix}; \quad q_{11}^{-1} = q_{12}q_{21} \neq  \pm 1, \omega^{-1};
&d &= 9\ord (q_{11})\ord (q_{12}q_{21}\omega).&
\\ \notag \\
\label{mg3} &\begin{pmatrix}
-1 & q_{12} \\ q_{21} & \omega
\end{pmatrix}; \quad q_{12}q_{21} = - 1;
&d &= 108.&
\\ \notag \\
\label{mg4} &\begin{pmatrix}
-1 & q_{12} \\ q_{21} & \omega
\end{pmatrix}; \quad q_{12}q_{21} = \omega;
&d &= 72.&
 \\ \notag \\
\label{mg5} &\begin{pmatrix}
-1 & q_{12} \\ q_{21} & \omega
\end{pmatrix}; \quad q_{12}q_{21} = - \omega;
&d &= 36.&
\\ \notag \\
\label{mg6} &\begin{pmatrix}
-1 & q_{12} \\ q_{21} & \zeta
\end{pmatrix}; \quad q_{12}q_{21} = \zeta^{-2};
&d &=  4\ord ( \zeta) \ord (- \zeta^{-1}).&
\end{flalign}}

\medbreak
\subsection{Braidings of finite non-abelian group type}

\

We begin with a class of examples studied in \cite{MiS}.

Let $\Gamma$ be a group and $T \subset \Gamma$ a subset such that for all
$g \in \Gamma, t
\in T, gtg^{-1} \in T$. Thus $T$ is a union of conjugacy classes of $\Gamma$.
Let $\phi : \Gamma \times T \to \ku \setminus \{0\}$ be a
function such that for all $g,h \in \Gamma$ and $t \in T$,
\begin{flalign}
\phi(1,t) &= 1,& \label{phi1}\\
\phi(gh,t) &= \phi(g,hth^{-1}) \phi(h,t). &\label{phi2}
\end{flalign}
We can then define a Yetter-Drinfeld module $V = V(\Gamma,T, \phi)$ over $\Gamma$ with
$\ku$-basis $x_{t}$, $t \in T,$ and action and coaction of $\Gamma$ given by
\begin{flalign}
gx_{t} &= \phi(g,t)x_{gtg^{-1}},&\label{action}\\
\delta(x_{t}) &= t \otimes x_{t} &\label{coaction}
\end{flalign}
for all $g \in \Gamma, t \in T$.

Conversely, if the function $\phi$ defines a Yetter-Drinfeld module on the vector space $V$ by \eqref{action}, \eqref{coaction}, then $\phi$ satisfies \eqref{phi1}, \eqref{phi2}.

Note that the braiding $c$ of $V(\Gamma,T, \phi)$ is determined by
$$c(x_{s} \otimes x_{t}) = \phi(s,t) x_{sts^{-1}} \otimes x_{t} \text{ for all } s,t \in T,$$
hence by the values of $\phi$ on $T \times T$.

\medbreak
The main examples come from the theory of Coxeter groups (\cite[Chapitre
IV]{BL}). Let $S$ be a subset of a group $W$ of elements of order $2$. For
all $s,s' \in S$ let $m(s,s')$ be the order of $ss'$. $(W,S)$ is called a
{\em Coxeter system} and $W$ a {\em Coxeter group} if $W$ is generated by
$S$ with defining relations $(ss')^{m(s,s')} = 1$ for all $s,s' \in S$ such
that $m(s,s')$ is finite.

\medbreak
Let $(W,S)$ be a Coxeter system. For any $g \in W$ there is a sequence
$(s_{1}, \dots, s_{q})$ of elements in $S$ with $g = s_{1} \cdot \dots \cdot
s_{q}$. If $q$ is minimal among all such representations, then $q = l(g)$ is
called the {\em length} of $g$, and $(s_{1}, \dots, s_{q})$  is a {\em
reduced representation} of $g$.

\begin{definition}\label{Coxeterexample}
Let $(W,S)$ be a Coxeter system, and $T = \{gsg^{-1} \mid g \in W, s \in S
\}$. Define $\phi : W \times T \to \ku \setminus \{0\}$ by
\begin{equation}\label{Coxeterphi}
\phi(g,t) = (-1)^{l(g)} \text{ for all } g \in W, t \in T.
\end{equation}
\end{definition}
This $\phi$ satisfies \eqref{phi1} and \eqref{phi2}. Thus we have associated to each Coxeter group the Yetter-Drinfeld module $V(W, T,\phi) \in \ydw$.

\medbreak

The functions $\phi$ satisfying \eqref{phi1}, \eqref{phi2} can be constructed up to a diagonal change of the basis from characters of the centralizers of elements in the conjugacy classes. This is a special case of the description of the simple modules in $\ydg$ (see \cite{Wi} and also \cite{L4}); the equivalent classification of the simple Hopf bimodules over $\Gamma$ was obtained in \cite{DPR}
(over $\ku$) and then in \cite{Ci} (over any field).

\medbreak

Let $t$ be an element in $\Gamma$. We denote by $\mathcal O_{t}$ and $\Gamma^{t}$ the
conjugacy class and the centralizer of $t$ in $\Gamma$. Let $U$ be any
left $\ku \Gamma^{t}$-module. It is easy to see that the induced representation
$V = \ku \Gamma \otimes _{\ku \Gamma^{t}} U$ is a Yetter-Drinfeld module over
$\Gamma$ with the induced action of $\Gamma$ and the coaction
$$\delta: V \to \ku \Gamma \otimes V, \quad  \delta(g \otimes u) =
gtg^{-1} \otimes g \otimes u \text{ for all } g \in \Gamma, u \in U.$$
We will denote this Yetter-Drinfeld module over $\Gamma$ by $M(t,U)$.

Assume that $\Gamma$ is finite. Then $V = M(t,U)$ is a simple Yetter-Drinfeld module if $U$
is a simple representation of $\Gamma^{t}$, and each simple module in $\ydg$ has this form.
If we take from each conjugacy class one element $t$ and non-isomorphic simple
$\Gamma^{t}$-modules, any two of these simple Yetter-Drinfeld modules are non-isomorphic.

Let $s_{i}, 1 \leq i \leq \theta$, be a complete system of representatives of the
residue classes of $\Gamma^{t}$. We define $t_{i} = s_{i}ts_{i}^{-1}$ for all
$1 \leq i \leq \theta$. Thus
$$\Gamma/\Gamma^{t} \to \mathcal O_{t}, \quad s_{i}\Gamma^{t} \mapsto t_{i}, 1 \leq i \leq \theta,$$
is bijective, and as a vector space, $V = \bigoplus_{1 \leq i \leq \theta} s_{i} \otimes U$.
For all $g \in \Gamma$ and $1 \leq i \leq \theta$, there is a uniquely determined
$1 \leq j \leq \theta$ with $s_{j}^{-1}gs_{i} \in \Gamma^{t}$, and the action of
$g$ on $s_{i} \otimes u, u \in U$, is given by
$$g(s_{i} \otimes u) = s_{j} \otimes (s_{j}^{-1}gs_{i})u.$$
In particular, if $U$ is a one-dimensional $\Gamma^{t}$-module with basis $u$ and action
$hu = \chi(h)u$ for all $h \in \Gamma^{t}$ defined by the character $\chi :
\Gamma^{t} \to \ku \setminus \{0\}$, then $V$ has a basis $x_{i} = s_{i} \otimes u,
1 \leq i \leq \theta$, and the action and coaction of $\Gamma$  are given by
$$
gx_{i} = \chi(s_{j}^{-1}gs_{i})x_{j} \quad \text{ and }  \delta(x_{i}) = t_{i}\otimes x_{i},
$$
 if $s_{j}^{-1}gs_{i} \in \Gamma^{t}$. Note that $gt_{i}g^{-1} = t_{j}$. Hence the
module we have constructed is $V(\Gamma,T,\phi)$, where $T$ is the conjugacy class
of $t$, and $\phi$ is given by $\phi(g,t_{i}) = \chi(s_{j}^{-1}gs_{i})$.

\medbreak

We now construct another example of a function $\phi$ satisfying \eqref{phi1}, \eqref{phi2}.

\begin{definition}
Let $T$ be the set of all transpositions in the symmetric group $\mathbb{S}_{n}$.
Define $\phi : \mathbb{S}_{n} \times T \to \ku \setminus \{0\}$ for all $g \in
\mathbb{S}_{n}, 1 \leq i < j \leq n,$ by
\begin{equation}\label{symmetric2}
\phi(g, (ij)) =
\begin{cases}
1   &, \text{ if } g(i ) < g(j),\\
-1   &, \text{ if } g(i) > g(j).
\end{cases}
\end{equation}
\end{definition}

Let $t = (12)$. The centralizer of $t$ in $\mathbb{S}_{n}$ is $\langle (34), (45),
\dots, (n-1,n)\rangle \cup
\langle (34), (45), \dots, (n-1,n)\rangle (12)$. Let $\chi$ be the character of $(\mathbb{S}_{n})^{t}$ with $\chi((ij)) = 1$ for all $3 \leq i < j \leq n$, and
$\chi((12)) = -1$. Then the function $\phi$ defined by \eqref{symmetric2} is
given by the character $\chi$ as described above.

\medbreak
Up to base change we have found all functions $\phi$ satisfying \eqref{phi1}, \eqref{phi2}
for $\Gamma = \mathbb{S}_{n}$, where $T$ is the conjugacy class of all transpositions,
and $\phi(t,t) = -1$ for all $t  \in T$. The case $\phi(t,t)= 1$ for some $t \in T$
would lead to a Nichols algebra $\mathfrak{B}(V)$ of infinite dimension.

\medbreak
To determine the structure of $\toba(V)$ for the Yetter-Drinfeld modules defined
by the functions $\phi$ in \eqref{Coxeterphi} and \eqref{symmetric2} seems to be
a fundamental and very hard combinatorial problem. Only a few partial results are
known \cite {MiS}, \cite{FK}, \cite{FP}.

\medbreak
We consider some special cases; here the method of skew-derivations
is applied, see Proposition \ref{derivation}.

\begin{exa}\label{Sn}
Let $W = \mathbb{S}_{n}, n \geq 2$, and $T = \{(ij) \mid 1 \leq i < j \leq n\}$ the
set of all transpositions. Define $\phi$ by \eqref{Coxeterphi} and let $V =
V(W,T,\phi)$. Then the following relations hold in $\toba(V)$ for all $1
\leq i < j \leq n, 1 \leq k < l \leq n$:
\begin{align}
x_{(ij)}^{2}& = 0. \label{qrn}\\
\text{If } \{i,j\} \cap \{k,l\} = \emptyset, \text{ then }\qquad
x_{(ij)}x_{(kl)} + x_{(kl)}x_{(ij)}& = 0. \label{order2} \\
\text{If } i < j < k, \text{ then }\qquad   x_{(ij)}x_{(jk)} +
x_{(jk)}x_{(ik)} + x_{(ik)}x_{(ij)}& = 0, \label{order3} \\
 x_{(jk)}x_{(ij)} + x_{(ik)}x_{(jk)} + x_{(ij)}x_{(ik)}& = 0. \notag
\end{align}

\end{exa}

\begin{exa}\label{Sn'}
Let $W = \mathbb{S}_{n}, n \geq 2$, and $T = \{(ij) \mid 1 \leq i < j \leq n\}$ the
set of all transpositions. Define $\phi$ by \eqref{symmetric2} and let $V =
V(W,T,\phi)$. Then the following relations hold in $\toba(V)$ for all $1
\leq i < j \leq n, 1 \leq k < l \leq n$:
\begin{align}
x_{(ij)}^{2}& = 0. \label{qrn'}\\
\text{If } \{i,j\} \cap \{k,l\} = \emptyset, { then }\qquad
x_{(ij)}x_{(kl)} - x_{(kl)}x_{(ij)}& = 0. \label{order2'} \\
\text{If } i < j < k, { then } \qquad  x_{(ij)}x_{(jk)} - x_{(jk)}x_{(ik)} -
x_{(ik)}x_{(ij)}& = 0, \label{order3'}\\
 x_{(jk)}x_{(ij)} - x_{(ik)}x_{(jk)} - x_{(ij)}x_{(ik)}& = 0. \notag
\end{align}

\end{exa}

The algebras $\widetilde{\toba}(V)$ generated by all $x_{(ij)}, 1 \leq i <  j \leq n,$ with the
quadratic relations in the examples \ref{Sn} resp. \ref{Sn'} are braided Hopf algebras in the category of Yetter-Drinfeld modules over $\mathbb{S}_{n}$. $\widetilde{\toba}(V)$ in example \ref{Sn'} is the
algebra
$\mathcal{E}_{n}$ introduced by Fomin and Kirillov in \cite{FK} to describe the cohomology
ring of the flag variety. We believe that indeed the quadratic relations in
the examples \ref{Sn} and \ref{Sn'} are defining relations for $\toba(V)$,
that is $\widetilde{\toba}(V) = \toba(V)$ in these cases.

\medbreak
It was noted in \cite{MiS} that the conjecture in \cite{FK} about the ''Poincar\'{e}-duality`` of the dimensions of the homogeneous components of the algebras $\mathcal{E}_{n}$ (in case they are finite-dimensional) follows from the braided Hopf algebra structure as a special case of Lemma \ref{Poincare}.

\medbreak
Another result about the algebras $\mathcal{E}_{n}$ by Fomin and Procesi \cite{FP} says that  $\mathcal{E}_{n+1}$ is a free module over $\mathcal{E}_{n}$, and
$ P_{\mathcal{E}_{n}}$ divides $P_{\mathcal{E}_{n+1}}$, where $P_{A}$ denotes the Hilbert series of a graded algebra $A$. The proof in \cite{FP} used the relations in Example \ref{Sn'}.

\medbreak
This result is in fact a special case of a very general splitting theorem for braided Hopf algebras in \cite[Theorem3.2]{MiS} which is an application of the  fundamental theorem for Hopf modules in the braided situation. This splitting theorem generalizes the main result of \cite{Gn3}.

\medbreak

In \cite{MiS} some partial results are obtained about the structure of the Nichols algebras of Coxeter groups. In particular

\begin{theorem}\cite[Corollary 5.9]{MiS}\label{MiS}
Let $(W,S)$ be a Coxeter system, $T$ the set of all $W$-conjugates of
elements in $S$, $\phi$ defined by \eqref{Coxeterphi}, $V = V(W,T,\phi)$ and $R = \toba(V)$. For all $g \in W$, choose a reduced representation $g = s_{1} \cdots
s_{q}$, $s_{1}, \cdots ,s_{q} \in S$, of $g$, and define
$$x_{g} = x_{s_{1}} \cdots x_{s_{q}}.$$

Then the subalgebra of $R$
generated by all $x_{s}, s \in S$ has the $\ku$-basis $x_{g}, g \in W$.
For all $g \in W$, the $g$-homogeneous component $R_{g}$ of $R$ is isomorphic to $R_{1}$.

If $R$ is finite-dimensional, then $W$ is finite and $dim(R) = ord(W)
dim(R_{1})$. \qed
\end{theorem}

\medbreak
This theorem holds for more general functions $\phi$, in particular for $S_{n}$ and $\phi$ defined in \eqref{symmetric2}.

\medbreak
Let $(W,S)$ be a Coxeter system and $V = V(W,T,\phi)$ as in Theorem \cite{MiS}. Then $\toba(V)$ was computed in \cite{MiS} in the following cases:
\begin{itemize}
\item $W= \mathbb{S}_{3}, S=\{(12),(23)\}$: The relations of $\toba(V)$ are the quadratic relations in Example \ref{Sn}, and $\text{dim}\toba(V) = 12$.

\medbreak
\item $W= \mathbb{S}_{4}, S=\{(12),(23),(34)\}$: The relations of $\toba(V)$ are the quadratic relations in Example \ref{Sn}, and $\text{dim}\toba(V) = 24 \cdot 24$.

\medbreak
\item $W=D_{4}$, the dihedral group of order 8, $S =\{t,t'\}$, where $t,t'$ are generators of $D_{4}$ of order 2 such that $tt'$ is of order 4. There are quadratic relations and relations of order 4 defining $\toba(V)$, and  $\text{dim}\toba(V) = 64$.
\end{itemize}

In all three cases the integral, which is the longest non-zero word in the generators $x_{t}$, can be described in terms of the longest element in the Coxeter group.
In all the other cases it is not known whether $\toba(V)$ is finite-dimensional.

\medbreak

In \cite{FK} it is shown that
\begin{itemize}
\item dim$(\mathcal{E}_{3}) = 12.$

\medbreak
\item dim$(\mathcal{E}_{4}) = 24 \cdot 24.$

\medbreak
\item dim$(\mathcal{E}_{5})$ is finite by using a computer program.
\end{itemize}
Again, for the other cases $n > 5$ it is not known whether $\mathcal{E}_{n}$ is finite-dimensional.

\medbreak
In \cite[5.3.2]{Gn4} another example of a finite-dimensional Nichols algebra of a braided vector space $(V,c)$ of finite group type is given with dim$(V)=4$ and dim$(\toba(V)) = 72.$ The defining relations of $\toba(V)$ are quadratic and of order 6.

\medbreak

By a result of Montgomery \cite{M1}, any pointed Hopf algebra $B$ can be decomposed as a crossed
product
$$B \simeq A  \#_{\sigma} \ku G, \quad \sigma \text{ a 2-cocycle}$$
of A, its link-indecomposable component containing 1 (a
Hopf subalgebra) and a group algebra $\ku G$. However, the structure of such link-indecomposable Hopf algebras A, in particular in the case when $A$ is
finite-dimensional and the group of its group-like elements $G(A)$ is
non-abelian, is not known. To define {\it link-indecomposable pointed Hopf algebras}, we recall the definition of the {\it quiver}  of $A$ in \cite{M1}.
The vertices of the quiver of $A$ are the elements of the group $G(A)$; for $g,h \in G(A)$, there
exists an arrow from $h$ to $g$ if $P_{g,h}(A)$ is non-trivial, that is if $\ku (g - h) \subsetneqq P_{g,h}(A)$. The Hopf algebra $A$ is called link-indecomposable, if its quiver is connected as an undirected graph.

\begin{definition}
Let $\Gamma$ be a finite group and $V \in \ydg$. V is called {\it link-indecomposable} if  the group $\Gamma$ is generated by the elements $g$ with $V_{g} \neq 0$.
\end{definition}
By \cite[4.2]{MiS}, $V \in \ydg$ is link-indecomposable if and only if the Hopf algebra $\toba(V) \# \ku \Gamma$ is link-indecomposable.

Thus by the examples constructed above, there are link-indecomposable, finite-dimensional pointed Hopf algebras $A$ with $G(A)$ isomorphic to $S_{n}, 3 \leq n \leq 5$, or to $D_{4}$.

\begin{question}
Which finite groups are isomorphic to $G(A)$ for some finite-dimensional, link-indecomposable pointed Hopf algebra $A$? Are there finite groups which do not occur in this form?
\end{question}

Finally, let us come back to the simple Yetter-Drinfeld modules $V = M(t,U) \in \ydg$, where $t \in \Gamma$ and $U$ is a simple left $\Gamma^{t}$-module of dimension $>1$. In this case, strong restrictions are known for $\toba(V)$ to be finite-dimensional.
By Schur's lemma, $t$ acts as a scalar $q$ on $U$.

\begin{prop}\label{conddegree} \cite[3.1]{Gn4} Assume that $\dim \toba (V)$ is finite.
If $\dim U \ge 3$, then  $q=-1$;
and if $\dim U =2$, then  $q=-1$ or $q$ is a  root of unity of order three. \qed
\end{prop}

\medbreak
In the proof of Proposition \ref{conddegree}, a result of Lusztig on braidings of Cartan type (see \cite[Theorem 3.1]{AS2}) is used. In a similar way Gra{\~n}a showed

\begin{prop}\label{finiteB} \cite[3.2]{Gn4}
Let $\Gamma$ be a finite group of odd order, and $V \in \ydg$. Assume that $\toba(V)$ is finite-dimensional. Then the multiplicity of any simple Yetter-Drinfeld module over $\Gamma$ as a direct summand in $V$ is at most 2.

In particular, up to isomorphism there are only finitely many Yetter-Drinfeld modules $V \in \ydg$ such that $\toba(V)$ is finite-dimensional.
\qed \end{prop}

\medbreak
The second statement in Proposition \ref{finiteB} was a conjecture in a preliminary version of \cite{AS2}.

\medbreak
\subsection{Braidings of (infinite) group type}

\

We briefly mention Nichols algebras over a free abelian
group of finite rank with a braiding which is not diagonal.

\begin{exa} Let $\Gamma = \langle g\rangle$ be a free group in one generator.
Let $\Vdos $ be the Yetter-Drinfeld module of dimension 2 such that $\Vdos  = \Vdos _{g}$
and the action of $g$ on $\Vdos$ is given, in a basis $x_{1}, x_{2}$, by
$$ g\cdot  x_{1} = t x_{1}, \qquad g \cdot x_{2} = t x_{2} + x_{1}. $$
Here $t\in \ku^{\times}$. Then:

(a). If $t$ is not a root of 1, then $\toba(\Vdos) = T(\Vdos )$.

(b). If $t= 1$, then $\toba(\Vuno) = \ku<x_{1},x_{2} \vert x_{1}x_{2} = x_{2}x_{1} + x_{1}^2>$;
this is the well-known Jordanian quantum plane.
\end{exa}

\begin{exa} More generally, if $t\in \ku^{\times}$,
let $\Vn$ be the Yetter-Drinfeld module of dimension $\theta \ge 2$ such that $\Vn = \Vn_{g}$
and the action of $g$ on $\Vn$ is given, in a basis $x_{1}, \dots, x_{\theta}$, by
$$ g \cdot x_{1} = t x_{1}, \qquad g \cdot x_{j} = t x_{j} + x_{j -1}, \quad 2\le j \le \theta. $$
Note there is an inclusion of Yetter-Drinfeld modules $\Vdos \hookrightarrow \Vn$;
hence, if $t$ is not a root of 1,  $\toba(\Vn)$ has exponential growth. \end{exa}

\begin{question} Compute $\toba(\Vnuno)$; does it have finite growth?
\end{question}

\section{Nichols algebras of Cartan type}\label{examnichols-cartan}

We now discuss fundamental examples of Nichols algebras
of diagonal type that come from the theory of quantum groups.

\medbreak We first need to fix some notation.
Let $A = (a_{ij})_{1\le i, j\le \theta}$ be a generalized symmetrizable Cartan matrix \cite{K};
let $(d_{1}, \dots, d_{\theta})$ be positive integers such that $d_{i}a_{ij} = d_{j}a_{ji}$.
Let $g$ be the Kac-Moody algebra  corresponding to the Cartan matrix $A$.
Let $\mathcal  X$ be the set of connected components of the Dynkin diagram
corresponding to it.
For each $I \in \mathcal  X$, we let
$\mathfrak g_{I}$
be the Kac-Moody Lie algebra corresponding to the generalized Cartan matrix
$(a_{ij})_{i,j \in I}$ and $\mathfrak n_{I}$ be the Lie subalgebra of
$\mathfrak g_{I}$ spanned by all its  positive roots. We omit the subindex $I$ when $I = \{1, \dots, \theta\}$.
We assume that for each $I\in \mathcal X$, there exist $c_I, d_I$
such that $I = \{j: c_I\le j \le d_I\}$; that is, after reordering
the Cartan matrix is
 a matrix of blocks corresponding to the connected components.
Let $I\in {\mathcal X}$ and $i \sim j$ in $I$;
then $N_{i} = N_{j}$, hence $N_{I} := N_{i}$ is well defined.
Let $\Phi_{I}$, resp. $\Phi_{I}^{+}$, be the root system, resp. the subset of positive roots,
corresponding to the Cartan matrix
$(a_{ij})_{i, j\in I}$; then $\Phi = \bigcup_{I\in \mathcal X}\Phi_{I}$, resp. $\Phi^+ =
\bigcup_{I\in \mathcal X}\Phi_{I}^{+}$ is the root system, resp. the subset of positive roots,
corresponding to the Cartan matrix
$(a_{ij})_{1\le i, j\le \theta}.$
 Let $\alpha_{1}, \dots,\alpha_{\theta}$  be the set of simple roots.

\medbreak
Let $\mathcal W_I$ be the Weyl group  corresponding to the Cartan matrix
$(a_{ij})_{i, j\in I}$; we identify it with a subgroup of the Weyl group
$\mathcal W$ corresponding to the Cartan matrix
$(a_{ij})$.

If $(a_{ij})$ is of finite type, we fix
a reduced decomposition of the longest element $\omega_{0, I}$ of $\mathcal W_{I}$
in terms of simple reflections.
Then we obtain a reduced decomposition of the longest element $\omega_{0}
= s_{i_1} \dots s_{i_P}$ of $\mathcal W$
from the expression of $\omega_{0}$ as product of the $\omega_{0, I}$'s in some
fixed order of the components, say the order arising from the order of the vertices.
Therefore
$\beta_{j} := s_{i_1} \dots s_{i_{j-1}}(\alpha_{i_j})$
is a numeration of $\Phi^+$.

\begin{exa}
Let $q\in \ku$, $q\neq 0$, and consider the braided vector space
$({\mathbb V}, c)$,
where ${\mathbb V}$
is a vector space  with a basis $x_{1}, \dots, x_{\theta}$ and
the braiding $c$ is given by
\begin{equation}\label{braiding-cartan}
c(x_{i}\otimes x_{j}) = q^{d_{i}a_{ij}} x_{j} \otimes x_{i},\end{equation}
\end{exa}

\begin{theorem} \label{quantum+}  \cite[33.1.5]{L3} Let $({\mathbb V}, c)$ be a
braided vector space with braiding matrix  \eqref{braiding-cartan}.
If $q$ is not algebraic over $\Q$, then
$$\toba({\mathbb V}) = \ku\langle x_{1},\dots, x_{\theta}
\vert \ad_{c}(x_{i})^{1 - a_{ij}} x_{j} = 0, \quad i\neq j\rangle.$$
\qed
\end{theorem}

The Theorem says that $\toba({\mathbb V})$ is the well-known
"positive part" $U_{q}^{+}(g)$ of the Drinfeld-Jimbo quantum enveloping algebra of $g$.

\medbreak
To state the following important Theorem, we recall the definition of braided commutators \eqref{br-adj}. Lusztig defined root vectors $X_{\alpha} \in {\mathfrak B}({\mathbb V})$,
$\alpha \in \Phi^{+}$ \cite{L2}. One can see from \cite{L1, L2} that,
up to a non-zero scalar, each root vector can be written as an
iterated braided commutator in some sequence
$X_{\ell_{1}}, \dots, X_{\ell_{a}}$ of simple root vectors
such as $[[X_{\ell_{1}}, [X_{\ell_{2}}, X_{\ell_{3}}]_{\mathfrak c}]_{\mathfrak c},
[X_{\ell_{4}}, X_{\ell_{5}}]_{\mathfrak c}]_{\mathfrak c}$.
See also  \cite{Ri}.

\begin{theorem}\label{flkernels} \cite{L1, L2, L3, Ro1, Mu}.
 Let $({\mathbb V}, c)$ be a braided vector space
with braiding matrix  \eqref{braiding-cartan}. Assume that
$q$ is a root of 1 of odd order $N$; and that 3 does not divide
$N$ if there exists $I \in \mathcal  X$ of type $G_2$.

The algebra $\toba(\mathbb V)$ is finite dimensional if and only if $(a_{ij})$ is a finite Cartan matrix.

If this happens, then  ${\mathfrak B}({\mathbb V})$ can be presented by generators $X_{i}$, $1\le i \le \theta$, and relations
\begin{align}
\ad _{c}(X_{i})^{1-a_{ij}}(X_{j}) & = 0, \qquad i\neq j, \\
X_{\alpha}^{N}&=0, \qquad \alpha \in \Phi^+.
\end{align}
Moreover, the following elements constitute a basis of
${\mathfrak B}({\mathbb V})$:
$$X_{\beta_{1}}^{h_{1}} X_{\beta_{2}}^{h_{2}} \dots X_{\beta_{P}}^{h_{P}},
\qquad 0 \le h_{j} \le N - 1, \quad  1\le
j \le P.$$
\qed\end{theorem}

 The Theorem says that $\toba({\mathbb V})$
is the well-known
"positive part" $\mathfrak u_{q}^+(g)$ of the so-called Frobenius-Lusztig kernel of $g$.

Motivated by the preceding Theorems and results, we introduce
the following notion, generalizing \cite{AS2} (see also \cite{FG}).

\begin{definition} Let $(V,c)$ a braided vector space of diagonal type
with basis $x_{1}, \dots, x_{\theta}$, and matrix $(q_{ij})$, that is
$$c(x_{i} \otimes x_{j}) = q_{i,j} x_{j} \otimes x_{i},
\text{ for all } 1 \leq  i,j \leq \theta.$$
We shall say that $(V, c)$  is {\it of Cartan type}
if  $q_{ii} \neq 1$ for all $i$, and  there are integers $a_{ij}$ with $a_{ii} = 2$,
$1\le i \le \theta$, and $0 \le -a_{ij} < \ord q_{ii}$ (which could be infinite),
$1\le i \neq j \le \theta$, such that
$$q_{ij}q_{ji} = q_{ii}^{a_{ij}}, \qquad 1\le i, j \le \theta.$$
\end{definition}

Since clearly
$a_{ij} =0$ implies that $a_{ji} =0$ for all $i\neq j$,
$(a_{ij})$ is a generalized Cartan matrix
in the sense of the book \cite{K}.
We shall adapt the terminology from generalized Cartan matrices and Dynkin diagrams to
braidings of Cartan type. For instance, we shall say that $(V, c)$ is of
{\it finite Cartan type} if it is of Cartan type and the corresponding GCM
is actually of finite type, {\it i. e.} a Cartan matrix associated to a
finite dimensional semisimple Lie algebra.
We shall say that a Yetter-Drinfeld module $V$ is {\it of Cartan type}
 if the matrix  $(q_{ij})$ as above is of Cartan type.

\begin{definition}
Let $(V,c)$ be a braided vector space of Cartan type with Cartan matrix $(a_{ij})$.
We say that $(V,c)$ is of {\it FL-type} (or Frobenius-Lusztig type)
if there exist positive integers $d_{1}, \dots, d_{\theta}$ such that
\begin{align}
 &\text{For all } i,j,\  d_{i} a_{ij} = d_{j} a_{ji} \text{ (thus } (a_{ij})
\text{ is symmetrizable).}\\
&\text{There exists a root of unity } q \in \ku \text{  such that }
q_{ij} = q^{d_{i}a_{ij}} \text{ for all } i,j.
\end{align}
We call $(V,c)$ {\it locally of FL-type} if any principal $2 \times 2$
submatrix of $(q_{ij})$ defines a braiding of FL-type.
\end{definition}

\medbreak
We now fix for each $\alpha \in \Phi^{+}$ such a representation
of $X_{\alpha}$ as an iterated braided commutator. For a general braided vector space
$(V, c)$ of finite Cartan type,
we define root vectors $x_{\alpha}$ in the tensor algebra $T(V)$,
$\alpha \in \Phi^{+}$, as the same formal iteration of braided commutators
in the elements $x_{1}, \dots, x_{\theta}$ instead of
$X_{1}, \dots, X_{\theta}$ but with respect to the braiding $c$ given by the
general matrix $\left(q_{ij}\right)$.

\begin{theorem}\label{tobas-cartantype} \cite[Th. 1.1]{AS2}, \cite[Th. 4.5]{AS4}.
Let $(V, c)$ be a braided vector space of Cartan type.  We also assume
that $q_{ij}$ has odd order for all $i,j$.

\medbreak

(i). Assume that $(V, c)$ is locally of
FL-type and that, for all $i$, the order of $q_{ii}$    is relatively
prime to 3 whenever $a_{ij} = - 3$ for some $j$, and
 is different from 3, 5, 7, 11, 13, 17.
If $\toba(V)$ is  finite dimensional, then $(V, c)$ is of finite Cartan type.

\medbreak

(ii). If $(V, c)$ is of finite Cartan type, then $\toba(V)$ is  finite
dimensional, and if moreover 3 does not divide the order of $q_{ii}$ for all $i$ in a connected component of the Dynkin diagram of type $G_2$, then
$$\dim \toba(V)
 = \prod_{I\in \mathcal  X}N_{I}^{ \dim \mathfrak n_{I}},$$
where $N_I = \ord(q_{ii})$ for all $i\in I$ and $I\in \mathcal  X$.
The Nichols algebra ${\mathfrak B}(V)$ is presented by generators
$x_{i}$, $1\le i \le \theta$, and relations
\begin{align} \label{serrebis}
\ad _{c}(x_{i})^{1-a_{ij}}(x_{j}) & = 0, \qquad i\neq j, \\ \label{rootbis}
x_{\alpha}^{N_{I}}&=0, \qquad \alpha \in \Phi_{I}^+,  \ I\in \mathcal X.
\end{align}
Moreover, the following elements constitute a basis of
${\mathfrak B}(V)$:
$$x_{\beta_{1}}^{h_{1}} x_{\beta_{2}}^{h_{2}} \dots x_{\beta_{P}}^{h_{P}}, \qquad 0 \le h_{j} \le N_{I} - 1, \text{ if }\, \beta_j \in I, \quad  1\le
j \le P.$$
 \qed \end{theorem}

Let $\wtoba(V)$ be the braided Hopf algebra in $\ydg$ generated by
$x_1, \dots, x_{\theta}$ with relations \eqref{serrebis}, where the $x_{i}$'s are
primitive. Let $\ftoba(V)$ be the subalgebra of $\wtoba(V)$ generated by
$x_{\alpha}^{N_{I}}$, $\alpha\in \Phi_I^+$, $I\in \mathcal X$; it is a Yetter-Drinfeld
submodule of $\wtoba(V)$.

\begin
{theorem}\label{powrootvec-subalg} \cite[Th. 4.8]{AS4} $\ftoba(V)$  is a
braided Hopf subalgebra in $\ydg$ of $\wtoba(V)$. \qed \end{theorem}

\section{Classification of pointed Hopf algebras by the lifting method}
\label{meth-class}

\subsection{Lifting of Cartan type}\label{method-strategy}

\

We propose  subdividing the classification problem for finite-dimensional
pointed Hopf algebras into the following problems:

\begin{itemize}
\item [(a).] Determine all braided vector spaces $V$ of group type
such that $\mathfrak  B(V)$  is finite dimensional.

\medbreak
\item [(b).] Given a finite group $\Gamma$, determine all realizations of braided
vector spaces $V$ as
in (a) as Yetter-Drinfeld modules over $\Gamma$.

\medbreak
\item [(c).]  The lifting problem: For $\mathfrak  B(V)$ as in (a),
compute all Hopf algebras $A$ such that $\gr A \simeq \mathfrak  B(V) \# H$ .

\medbreak
\item [(d).] Investigate whether any  finite dimensional pointed Hopf algebra   is generated as an algebra
by its group-like and skew-primitive elements.
\end{itemize}
\medbreak

Problem (a) was discussed in Chapters \ref{examnichols} and \ref{examnichols-cartan}.
We have seen the very important class of braidings of finite Cartan type and some
isolated examples where the Nichols algebra is finite-dimensional. But the general
case of problem (a) seems to require completely new ideas.

\medbreak
Problem (b) is of a computational nature. For braidings of finite Cartan type with
Cartan matrix $(a_{ij})_{1 \leq i,j \leq \theta}$ and an abelian group $\Gamma$
we have to compute elements $g_{1}, \dots, g_{\theta} \in \Gamma$ and characters
$\chi_{1}, \dots, \chi_{\theta} \in \widehat{\Gamma}$ such that
\begin{equation}\label{Cartantype}
\chi_{i}(g_{j}) \chi_{j}(g_{i}) = \chi_{i}(g_{i})^{a_{ij}}, \text{ for all } 1 \leq i,j \leq \theta.
\end{equation}
To find these elements one has to solve a system of quadratic congruences in several unknowns. In many cases they do not exist. In particular, if $\theta > 2(\text{ord}\Gamma)^{2}$, then the braiding cannot be realized over the group $\Gamma$. We refer to \cite[Section 8]{AS2} for details.

\medbreak
Problem (d) is the subject of Section \ref{degone}.

\medbreak We will now discuss the lifting problem (c).

\medbreak

 The coradical filtration $\ku \Gamma = A_{0} \subset A_{1} \subset \dots$ of a pointed Hopf algebra $A$ is stable under the adjoint action of the group.  For abelian groups $\Gamma$ and finite-dimensional Hopf algebras, the following stronger result holds.
It is the starting point of the lifting procedure, and we will use it several times.

 If $M$  is a $\ku \Gamma$-module, we denote by $M^{\chi} = \{ m \in M \mid gm
= \chi(g)m \text{ for all }g \in \Gamma\}$, $\chi \in \widehat{\Gamma}$,
the isotypic component of type $\chi$.

\begin{lema}\label{start}
Let $A$ be a finite-dimensional Hopf algebra with abelian group
$G(A) = \Gamma$ and diagram $R$. Let $V = R(1) \in \ydg$ with basis $x_{i} \in V_{g_{i}}^{\chi_{i}}, g_{i} \in \Gamma, \chi_{i} \in \widehat{\Gamma},
1 \leq i \leq \theta.$

(a). The isotypic component of trivial type of $A_1$ is $A_0$. Therefore,
$A_1 = A_0 \oplus (\oplus_{\chi \neq \varepsilon} (A_{1})^{\chi})$ and
\begin{equation}\label{start1}
\oplus_{\chi \neq \varepsilon} (A_{1})^{\chi} \xrightarrow{\simeq} A_{1}/A_{0}\xleftarrow{\simeq} V \# \ku \Gamma.
\end{equation}

(b). For all $g\in \Gamma$, $\chi\in \VGamma$
with $\chi\neq \varepsilon$,
\begin{align}\label{start2}
{\mathcal P}_{g,1}(A)^{\chi} &\neq 0 \iff \text{ there is some } 1\le \ell\le \theta:
g=g_{\ell}, \chi= \chi_{\ell}; \\ \label{start3}
{\mathcal P}_{g,1}(A)^{\varepsilon} &= \ku(1-g).
\end{align}
\end{lema}
\pf
(a) follows from \cite[Lemma 3.1]{AS1} and implies (b). See \cite[Lemma 5.4]{AS1}. \epf

We assume that $A$ is a finite-dimensional pointed Hopf algebra with abelian group $G(A) = \Gamma$, and that
$$\text{gr}A \simeq \toba(V) \# \ku \Gamma,$$
where $V \in \ydg$ is a given Yetter-Drinfeld module with basis $x_{i} \in V_{g_{i}}^{\chi_{i}}, g_{1} \dots, g_{\theta} \in \Gamma,\chi_{1}, \dots, \chi_{\theta} \in \widehat{\Gamma}, 1 \leq i \leq \theta.$

We first lift the basis elements $x_{i}$. Using \eqref{start1}, we choose $a_{i}\in \SkPr(A)_{g_{i}, 1}^{\chi_{i}}$ such that the canonical image of $a_{i}$ in $A_{1}/A_{0}$ is $x_{i}$ (which we identify with $x_{i} \# 1$), $1\leq i \leq \theta$. Since the elements $x_{i}$ together with $\Gamma$ generate $\toba(V) \# \ku \Gamma$, it follows from a standard argument that $a_{1}, \dots, a_{\theta}$ and the elements in $\Gamma$ generate $A$ as an algebra.

\medbreak

Our aim is to find relations between the $a_{i}'s$ and the elements in $\Gamma$ which define a quotient Hopf algebra of the correct dimension $\text{dim}\toba(V) \cdot \text{ord}(\Gamma)$. The idea is to ''lift`` the relations between the $x_{i}'s$ and the elements in $\Gamma$ in $\toba(V) \# \ku \Gamma$.

\medbreak

We now assume moreover that $V$ is of finite Cartan type with Cartan matrix $(a_{ij})$ with respect to the basis $x_{1},\dots,x_{\theta}$, that is \eqref{Cartantype} holds.
 We also assume
\begin{flalign}
&\text{ord}(\chi_{j}(g_{i})) \text{ is odd for all } i,j, & \label{cond1}\\
 &N_{i}= \text{ord}(\chi_{i}(g_{i})) \text{ is prime to 3} \text{ for all }i \in I, I \in \mathcal{X} \text{ of type }  G_2.\label{cond2}&
\end{flalign}

 We fix a presentation $\Gamma = \langle  y_{1}\rangle\oplus \dots\oplus  \langle  y_{\sigma}\rangle$, and denote by $M_{\ell}$ the order of $y_{\ell}$, $1\leq \ell \leq \sigma$.
Then Theorem \ref{tobas-cartantype} and formulas \eqref{smash1} imply that $\toba(V) \# \ku \Gamma$
can be presented by generators $h_{\ell}$, $1\leq
\ell \leq \sigma$, and $x_{i}$, $1\leq i \leq \theta$ with defining relations

\begin{flalign}
\label{qls5}
& h_{\ell}^{M_{\ell}} = 1, \quad 1\leq \ell \leq \sigma;&\\
\label{qls6}& h_{\ell}h_{t} = h_{t}h_{\ell}, \quad 1\leq t <\ell \leq \sigma;&\\
\label{qls7}
& h_{\ell} x_{i} = \chi_{i}(y_{\ell}) x_{i} h_{\ell}, \quad 1\leq
\ell \leq \sigma,\quad 1\leq i \leq \theta;&\\
\label{qls8}
& x_{\alpha}^{N_{I}}=0, \qquad \alpha \in \Phi_{I}^+, I\in \mathcal X ;&\\
\label{qls9}
& \ad (x_{i})^{1-a_{ij}}(x_{j})  = 0, \qquad i\neq j,&
\end{flalign}

and where the Hopf algebra structure is determined by
\begin{flalign}
\label{qls10}
& \Delta(h_{\ell}) = h_{\ell} \otimes h_{\ell}, \quad 1\leq \ell \leq \sigma;&\\
\label{qls11}
& \Delta(x_{i}) = x_{i}  \otimes 1 +  g_{i} \otimes x_{i}, \quad 1\leq i  \leq \theta.&
\end{flalign}

Thus $A$ is generated by the elements $a_{i}, 1 \leq i \leq \theta$, and $h_{l}, 1 \leq l \leq \sigma$. By our previous choice, relations \eqref{qls5}, \eqref{qls6}, \eqref{qls7} and \eqref{qls10}, \eqref{qls11} all hold in $A$ with the $x_{i}'s$ replaced by the $a_{i}'s$.

\medbreak

The remaining problem is to lift the quantum Serre relations \eqref{qls9} and the root vector relations \eqref{qls8}. We will do this in the next two Sections.

\medbreak
\subsection{Lifting the quantum Serre relations }\label{Serre}

\

We divide the problem into two cases.
\begin{itemize}
\item Lifting of the ``quantum Serre relations"
$ x_{i}x_{j} - \chi_{j}(g_{i}) x_{j}x_{i} = 0$, when $i\neq j$ are in different components of the Dynkin diagram.

\medbreak
\item Lifting of the ``quantum Serre relations"
$\ad _{c}(x_{i})^{1-a_{ij}}(x_{j})  = 0$, when $i\neq j$ are in the same component
of the Dynkin diagram.
\end{itemize}
The first case is settled in the next result from \cite[Theorem 6.8 (a)]{AS4}.

\begin{lema}\label{Serre1}
Assume that $1\leq i,j \leq \theta, i < j$ and $ i \not\sim j$. Then
\begin{equation}\label{linkSerre}
a_{i}a_{j} -\chi_{j}(g_{i}) a_{j}a_{i} = \lambda_{ij} (1 -
g_{i} g_{j}),
\end{equation}
where $\lambda_{ij}$ is a scalar in $\ku$ which can be chosen such that
\begin{flalign}
\label{qls4}
&\lambda_{ij} \text{ is arbitrary if } g_{i}g_{j} \neq 1 \text{ and }
\chi_{i}\chi_{j} = \varepsilon, \text{ but 0 otherwise}.&
\end{flalign}
\end{lema}
\pf
It is easy to check that $a_{i} a_{j} - \chi_{j}(g_{i}) a_{j} a_{i}
\in \SkPr(A)_{g_{i}g_{j}, 1}^{\chi_{i}\chi_{j}}$.
Suppose that $\chi_{i} \chi_{j} \neq \varepsilon$ and  $a_{i} a_{j}
- \chi_{j}(g_{i}) a_{j} a_{i} \neq 0$. Then by \eqref{start2},
$\chi_{i} \chi_{j} =\chi_{l}$ and $g_{i} g_{j} = g_{l}$ for some $1 \leq l \leq \theta$.

Substituting $g_{l}$ and $\chi_{l}$ in $\chi_{i}(g_{l}) \chi_{l}(g_{i}) = \chi_{i}(g_{i})^{a_{il}}$ and using $\chi_{i}(g_{j}) \chi_{j}(g_{i}) = 1$ (since $a_{ij} = 0$, since $i$ and $j$ lie in different components), we get $\chi_{i}(g_{i})^{2} =\chi_{i}(g_{i})^{a_{il}}.$

Thus we have shown that $a_{il} \equiv  2 \mod N_{i}$, and in the same way $a_{jl} \equiv  2 \mod N_{j}$. Since $i \not\sim j$, $a_{il}$ or $a_{jl}$ must be $0$, and we obtain the contradiction $N_{i} = 2$ or $N_{j} = 2$.

Therefore  $\chi_{i} \chi_{j} = \varepsilon$, and the claim follows from \eqref{start3}, or $a_{i} a_{j} - \chi_{j}(g_{i}) a_{j} a_{i} = 0$, and the claim is trivial.
\epf

Lemma \ref{Serre1} motivates the following notion.

\begin{definition}\label{link-dat}\cite[Definition 5.1]{AS4}
We say that two vertices $i$ and $j$
{\em are linkable} (or that $i$ {\em is linkable to} $j$) if
\begin{flalign}\label{link0} &i\not\sim j,&\\
\label{link1} &g_{i}g_{j} \neq 1 \text {  and}& \\
\label{link2}
&\chi_{i}\chi_{j} = \varepsilon.&
\end{flalign}
\end{definition}
The  following elementary properties are easily verified:
\begin{flalign}\label{link20} &\text{If $i$ is linkable to $j$, then } \chi_{i}(g_{j})\chi_{j}(g_{i}) = 1,
\quad \chi_{j}(g_{j}) = \chi_{i}(g_{i})^{-1}.& \\
\label{sndproplink}
&\text{If $i$ and $k$, resp. $j$ and $\ell$, are linkable,
then $a_{ij} = a_{k\ell}$, $a_{ji} = a_{\ell k}$.}& \\
\label{sndproplink2}
&\text{A vertex $i$ can not be linkable to two different vertices $j$ and $h$.}&
\end{flalign}
A {\em linking datum} is
a collection $(\lambda_{ij})_{1 \le i < j \le \theta, \, i\not\sim j}$ of elements in $\ku$
such that $\lambda_{ij}$ is arbitrary if $i$ and $j$
are linkable but 0 otherwise. Given a linking datum,
we say that two vertices $i$ and $j$
{\em are linked} if $\lambda_{ij}\neq 0$.

\medbreak

The notion of a linking datum encodes the information about lifting of relations in the first case.

\begin{definition}
The collection $\mathcal{D}$ formed by a finite Cartan matrix $(a_{ij})$, and $g_{1},\dots,g_{\theta} \in \Gamma, \chi_{1},\dots,\chi_{\theta} \in \widehat{\Gamma}$ satisfying \eqref{Cartantype}, \eqref{cond1} and \eqref{cond2}, and a linking datum $(\lambda_{ij})_{1 \le i < j \le \theta, \, i\not\sim j}$ will be called a {\it linking datum of finite Cartan type} for $\Gamma$.
We define the Yetter-Drinfeld module $V \in \ydg$ of $\mathcal{D}$ as the vector space with basis $x_{1} ,\dots,x_{\theta}$ with $x_{i} \in V_{g_{i}}^{\chi_{i}}$ for all $i$.

\medbreak

If $\mathcal{D}$ is a linking datum of finite Cartan type for $\Gamma$, we define the Hopf algebra $\mathfrak{u}(\mathcal{D})$ by generators $a_{i}, 1 \leq i \leq \theta$, and
$h_{l}, 1 \leq l \leq \sigma$ and the relations \eqref{qls5},\eqref{qls6},\eqref{qls7},\eqref{qls8}, the quantum Serre relations
\eqref{qls9} for $i \neq j$ and $i \sim j$,
\eqref{qls10},\eqref{qls11} with the $x_{i}$'s replaced by the  $a_{i}$'s, and the
lifted quantum Serre relations \eqref{linkSerre}.
\end{definition}

\medbreak
In the definition of $\mathfrak{u}(\mathcal{D})$ we could always assume that the
linking datum contains only elements $\lambda_{ij} \in \{0,1\}$ (by multiplying
the generators $a_{i}$ with non-zero scalars).

\begin{theorem}\label{littlelink}\cite[Th. 5.17]{AS4}
Let $\Gamma$ be a finite abelian group and $\mathcal{D}$ a linking datum of finite
Cartan type for $\Gamma$ with Yetter-Drinfeld module $V$. Then $\mathfrak{u}(\mathcal{D})$
is a finite-dimensional pointed Hopf algebra  with
$\text{gr}\mathfrak{u}(\mathcal{D})\simeq \toba(V) \# \ku \Gamma$. \qed
\end{theorem}

The proof of the Theorem is by induction on the number of irreducible components of
the Dynkin diagram. In the induction step a new Hopf algebra is constructed by
twisting the multiplication of the tensor product of two Hopf algebras by a 2-cocycle.
The 2-cocycle is defined in terms of the linking datum.

\medbreak
Note that the Frobenius-Lusztig kernel $\mathfrak{u}_{q}(\mathfrak{g})$ of a semisimple Lie algebra $\mathfrak{g}$ is a special case of $\mathfrak{u}(\mathcal{D})$. Here the Dynkin diagram of $\mathcal{D}$ is the disjoint union of two copies of the Dynkin diagram of $\mathfrak{g}$, and corresponding points are linked pairwise. But many other linkings are possible, for example 4 copies of $A_{3}$ linked in a circle \cite[5.13]{AS4}. See \cite{D} for a combinatorial description of all linkings of Dynkin diagrams.

\medbreak
Let us now turn to the second case. Luckily it turns out that (up to some small order exceptions) in the second case the Serre relations simply hold in the lifted situation without any change.

\begin{theorem}\label{QSR} \cite[Theorem 6.8]{AS4}.
Let $I \in \mathcal X$. Assume that $N_I\neq 3$.
If  $I$ is of type $B_n$, $C_n$ or $F_4$, resp. $G_2$, assume further that $N_I\neq 5$,
resp. $N_I\neq 7$. Then the quantum Serre relations  hold for all
$i, j \in I, i \neq j, i \sim j$.
\qed \end{theorem}

\medbreak
\subsection{Lifting the root vector relations}\label{lfpwr}

\

Assume first that the root $\alpha$ is simple
and corresponds to a vertex $i$. It is not difficult to see, using the quantum binomial formula,
that $a_{i}^{N_i}$ is a $(g_{i}^{N_{i}}, 1)$-skew-primitive.
By Lemma \ref{start}, we have
\begin{equation}\label{simpleroots}
a_{i}^{N_{i}} = \mu_{i}\left(1 - g_{i}^{N_{i}}\right),\end{equation}
for some scalar $\mu_{i}$; this scalar can be chosen so that
\begin{flalign}
\label{qls3}
& \mu_{i} \text{ is arbitrary if } g_{i}^{N_{i}} \neq 1 \text{ and }
\chi_{i}^{N_{i}} = 1 \text{ but 0 otherwise}.&  \end{flalign}

Now, if the root $\alpha$ is {\it not} simple  then   $a_{i}^{N_i}$ is not necessarily a skew-primitive, but a skew-primitive ``modulo root vectors of shorter length".

In general, we define the root vector $a_{\alpha}$ for $\alpha \in \Phi_{I}^{+}, I \in \mathcal{X}$, by replacing the $x_{i}$ by $a_{i}$ in the formal expression for $x_{\alpha}$ as a braided commutator in the simple root vectors. Then $a_{\alpha}^{N_{I}}, \alpha \in I$, should be an element $u_{\alpha}$ in the group algebra of the subgroup generated by the $N_{I}$-th powers of the elements in $\Gamma$.

Finally, the Hopf algebra generated by $a_{i}, 1 \leq i \leq \theta$, and $h_{l}, 1 \leq l \leq \sigma$ with the relations \eqref{qls5},\eqref{qls6},\eqref{qls7} (with $a_{i}$ instead of $x_{i}$),
\begin{itemize}
\item the lifted root vector relations $a_{\alpha}^{N_{I}} = u_{\alpha}, \alpha \in \Phi_{I}^{+},I \in \mathcal{X}$,

\medbreak
\item the quantum Serre relations  \eqref{qls9} for $i \neq j$  and $i \sim j$ (with $a_{i}$ instead of $x_{i}$),

\medbreak
\item the lifted quantum Serre relations  \eqref{linkSerre},
\end{itemize}
should have the correct dimension $\text{dim}(\toba(V)) \cdot \text{ord}(\Gamma)$.

We carried out all the steps of this program in the following cases:
\begin{itemize}
\item[(a)] All connected components of the Dynkin diagram are of type $A_{1}$ \cite{AS1}.

\medbreak
\item[(b)] The Dynkin diagram is of type $A_{2}$, and $N>3$ is odd \cite{AS3}.

\medbreak
\item[(c)] The Dynkin diagram is arbitrary, but we assume $g_{{i}}^{N_{i}} = 1$ for all $i$ \cite{AS4}.

\medbreak
\item[(d)] The Dynkin diagram is of type $A_{n}$, any $n \geq 2$, and $N>3$ , see Section \ref{A_n} of this paper.
\end{itemize}
The cases $A_{2}, N=3$ and $B_{2},N \text{ odd and } \neq 5$, were recently done in \cite{BDR}. Here $N$ denotes the common order of $\chi_{i}(g_{i})$ for all $i$ when the Dynkin diagram is connected.

\medbreak
\subsection{Generation in degree one}\label{degone}

\

Let us now discuss step (d) of the Lifting method.

It is not difficult to show that our conjecture \ref{conjecturedegree1} about Nichols
algebras, in the setting of $H = \ku \Gamma$,  is equivalent to

\begin{conjecture}\label{fingenone} \cite{AS3}. Any pointed finite dimensional  Hopf algebra over $\ku$ is generated by group-like and skew-primitive elements.\end{conjecture}

We have seen in Section \ref{DefNichols} that the corresponding conjecture is false  when the Hopf algebra is infinite-dimensional or when the Hopf algebra is finite-dimensional and the characteristic of the field is $>0$.
A strong indication that the conjecture is true is given by:

\begin{theorem}\label{degree1} \cite[Theorem 7.6]{AS4}.
Let $A$ be a finite-dimensional pointed Hopf algebra with coradical
$\ku \Gamma$ and  diagram $R$, that is
 $$\gr A \simeq R \# \ku \Gamma.$$
Assume that $R(1)$ is a Yetter-Drinfeld module of finite Cartan type with  braiding
$(q_{ij})_{1 \leq i,j \leq \theta}$. For all i, let $q_{i} = q_{ii}, N_i = \text {ord}(q_{i})$.
Assume that $\text{ord}(q_{ij})$ is odd and  $N_i$ is not divisible by 3 and $> 7$
for all $1 \leq i,j \leq \theta$.
\begin{enumerate}
\item For any $1 \leq i \leq \theta$ contained in a connected component of type $B_{n}$,
$C_{n}$ or $F_{4}$ resp. $G_{2}$, assume that $N_{i}$ is not divisible by 5 resp. by 5 or 7.
\label{degree2}

\medbreak
\item If  $i$ and $j$ belong to different components,  assume $q_{i}q_{j} = 1$
or $\ord(q_{i}q_{j}) = N_{i}$.\label{degree3}
\end{enumerate}
Then $R$ is generated as an algebra by $R(1)$, that is $A$ is generated
by skew-primitive and group-like elements. \qed
\end{theorem}

Let us discuss the idea of the proof of Theorem \ref{degree1}. At one decisive point, we use our previous results about braidings of Cartan type of rank 2.

Let $S$ be the graded dual of $R$. By the duality principle in Lemma \ref{dualityprinciple}, $S$ is generated in degree one since $P(R) = R(1)$. Our problem is to show that $R$ is generated in degree one, that is $S$ is a Nichols algebra.

Since $S$ is generated in degree one, there is a surjection of graded braided Hopf algebras $S \to \toba (V)$, where $V = S(1)$ has the same braiding as $R(1)$. But we know the defining relations of
$\toba (V)$, since it is of finite Cartan type. So we have
to show that these relations also hold in $S$.

In the case of a quantum Serre relation
$\ad _{c}(x_{i})^{1-a_{ij}}(x_{j}) = 0$, $i\neq j$, we consider the Yetter-Drinfeld submodule
$W$ of $S$ generated by $x_{i}$ and $\ad _{c}(x_{i})^{1-a_{ij}}(x_{j})$ and assume that $\ad _{c}(x_{i})^{1-a_{ij}}(x_{j})$  $\neq 0$. The assumptions (1)
and (2) of the Theorem guarantee that $W$ also is of Cartan type,
but not of finite Cartan type. Thus
$\ad _{c}(x_{i})^{1-a_{ij}}(x_{j}) = 0$ in $S$.

Since the quantum Serre relations hold in $S$, the root vector relations follow automatically from the next Lemma which is a consequence of Theorem \ref{powrootvec-subalg}.
\begin{lema}\label{rootv}\cite[Lemma 7.5]{AS4}
Let $S = \oplus_{n \ge 0} S(n)$ be a finite-dimensional graded Hopf algebra in
$^{\Gamma}_{\Gamma}\mathcal {YD}$ such that $S(0) = \ku 1$. Assume that $V = S(1)$
is of Cartan type with basis $(x_{i})_{1 \leq i,j \leq \theta}$ as described in
the beginning of this Section.
Assume the Serre relations $$(\ad_c x_{i})^{1 - a_{ij}} x_{j} = 0
\text{ for all } 1 \leq i,j \leq \theta, i \neq j \text{ and } i \sim j.$$
Then the root vector relations $$x_{\alpha}^{N_{I}} = 0, \quad\alpha \in \Phi^{+}_{I},
\quad I\in \mathcal X,$$ hold in $S$. \qed
\end{lema}

Another result supporting Conjecture \ref{fingenone} is:

\begin{theorem}\label{aeg-cotr}\cite[6.1]{AEG}
Any finite dimensional cotriangular pointed Hopf algebra is generated by
skew-primitive and group-like elements.
\end{theorem}
\qed

\medbreak
\subsection{Applications}\label{appl}

\

As a special case of the theory explained above we obtain a complete answer to the classification problem in a significant case.

\begin{theorem}\label{expp}\cite[Th. 1.1]{AS4}
Let $p$ be a prime $>17$, $s\geq1$, and $\Gamma = (\mathbb{Z}/(p))^{s}$.
Up to isomorphism there are only finitely many finite-dimensional pointed
Hopf algebras $A$ with $G(A) \simeq \Gamma$. They all have the form
$$A \simeq \mathfrak{u}(\mathcal{D}), \text{ where } \mathcal{D}
\text{ is a linking datum of finite Cartan type for }\Gamma.$$
\end{theorem}
\qed

If we really want to write down all these Hopf algebras we still have to solve the following serious problems:
\begin{itemize}
\item Determine all Yetter-Drinfeld modules $V$ over $\Gamma = (\mathbb{Z}/(p))^{s}$ of finite Cartan type.

\medbreak
\item Determine all the possible linkings for the modules $V$ over $(\mathbb{Z}/(p))^{s}$ in (a).
\end{itemize}
By \cite[Proposition 8.3]{AS2}, $\text{dim}V \leq 2s\dfrac{p-1}{p-2}$, for
all the possible $V$ in (a). This proves the finiteness statement in Theorem \ref{expp}.

\medbreak
Note that we have precise information about the dimension of the Hopf algebras
in \ref{expp}:
$$\dim \mathfrak{u}(\mathcal{D}) = p^{s|\phi^{+}|},$$
where $|\phi^{+}|$ is the number of the positive roots of the root system of
rank $\theta \leq 2s\dfrac{p-1}{p-2}$ of the Cartan matrix of $\mathcal{D}$.

\medbreak
For arbitrary finite abelian groups $\Gamma$, there usually are infinitely many non-isomorphic pointed Hopf algebras of the same finite dimension. The first examples were found in 1997 independently in \cite{AS1}, \cite{BDG}, \cite{G}. Now it is very easy to construct lots of examples by lifting. Using \cite[Lemma 1.2]{AS3} it is possible to decide when two liftings are non-isomorphic.

\medbreak
But we have a bound on the dimension of $A$:
\begin{theorem}\label{bounddim}\cite[Th. 7.9]{AS4}
For any finite (not necessarily abelian) group $\Gamma$ of odd order
there is a natural number $n(\Gamma)$ such that
$$\text{dim}A \leq n(\Gamma)$$
for any finite-dimensional pointed Hopf algebra $A$ with $G(A) = \ku \Gamma$.
\qed\end{theorem}

\begin{obs}\label{pointedp} As a corollary of the Theorem and its proof,
we get the complete  classification of all finite dimensional pointed Hopf
algebras with coradical of prime dimension $p$, $p \neq 2, 5, 7$.
By \cite[Theorem 1.3]{AS2}, the only possibilities for the Cartan matrix of
$\mathcal{D}$ with $\Gamma$ of odd prime order $p$ are
\begin{itemize}
\item[(a)] $A_{1}$ and $A_{1} \times A_{1}$,

\medbreak
\item[(b)]$A_{2}$, if  $p=3$  or $ p \equiv 1 \mod 3$,

\medbreak
\item[(c)] $B_{2}$,  if $p \equiv 1 \mod 4$,

\medbreak
\item[(d)] $G_{2}$,  if $  p \equiv 1 \mod 3$,

\medbreak
\item [(e)] $A_{2} \times A_{1}$ and $A_{2} \times A_{2}$, if $p=3$.
\end{itemize}
The Nichols algebras over $\mathbb{Z}/(p)$ for these Cartan matrices are listed in
\cite[Theorem 1.3]{AS2}. Hence we obtain from Theorem \ref{expp} for
$p \neq 2, 5, 7$ the bosonizations of the Nichols algebras, the liftings in case (a),
that is quantum lines and quantum planes \cite{AS1}, and the liftings of
type $A_{2}$ \cite{AS3} in case (b).

This result was also obtained by Musson \cite{Ms}, using the lifting method and \cite{AS2}.

The case $p=2$ was already done in \cite{N}. In this case the dimension of the pointed Hopf algebras with 2-dimensional coradical is not bounded.
\end{obs}

\medbreak
Let us mention briefly some classification results for Hopf algebras of special
order which can be obtained by the methods we have described. Let $p>2$ be a prime.
Then all pointed Hopf algebras $A$ of dimension $p^{n}, 1 \leq n \leq 5$ are known.
If the dimension is $p$ or $p^{2}$, then $A$ is a group algebra or a Taft Hopf algebra.
The cases of dimension $p^{3}$ and $p^{4}$ were treated in \cite{AS1} and \cite{AS3},
and the classification of dimension $p^{5}$ follows from \cite{AS4} and
\cite{Gn2}. Independently and by other methods, the case $p^{3}$ was also
solved in \cite{CD} and \cite{SvO}.

\medbreak
See \cite{bariloche} for a discussion of what is known on classification
of finite dimensional Hopf algebras.

\medbreak
\subsection{The infinite-dimensional case}\label{infinite-dim}

\

Our methods are also useful in the infinite-dimensional case.
Let us introduce the analogue to FL-type for infinite-dimensional Hopf algebras.

\begin{definition}
Let $(V,c)$ be a braided vector space of Cartan type with Cartan matrix $(a_{ij})$.
We say that $(V,c)$ is of {\it DJ-type} (or Drinfeld-Jimbo type)
if there exist positive integers $d_{1}, \dots, d_{\theta}$ such that
\begin{align}
 &\text{For all } i,j,\  d_{i} a_{ij} = d_{j} a_{ji} \text{ (thus } (a_{ij})
\text{ is symmetrizable).}\\
&\text{There exists  } q \in \ku \text{, which is not a root of unity,  such that }
q_{ij} = q^{d_{i}a_{ij}} \text{ for all } i,j.
\end{align}
\end{definition}

To formulate a classification result for infinite-dimensional Hopf algebras,
we now assume that $\ku$ is the field of complex numbers and we introduce
a notion from \cite{AS5}.

\begin{definition}
The collection $\mathcal{D}$ formed by a
{\it free abelian group} $\Gamma$ of finite rank, a finite Cartan matrix
$(a_{ij})_{1 \leq i,j \leq \theta}$, $g_{1},\dots,g_{\theta} \in \Gamma, \chi_{1},\dots,\chi_{\theta} \in \widehat{\Gamma}$, and a linking datum
$(\lambda_{ij})_{1 \le i < j \le \theta, \, i\not\sim j}$, will be called a
{\it positive datum of finite Cartan type } if
$$\chi_{i}(g_{j}) \chi_{j}(g_{i}) = \chi_{i}(g_{i})^{a_{ij}}, \text{ and }
1 \neq \chi_{i}(g_{i}) > 0, \text{ for all } 1 \leq i,j,\leq \theta. $$
\end{definition}

Notice that the restriction of the braiding of a positive datum of finite Cartan type
to each connected component is twist-equivalent to a braiding of DJ-type.

\medbreak
If $\mathcal{D}$ is a positive datum we define the Hopf algebra $U(\mathcal{D})$
by generators $a_{i}, 1 \leq i \leq \theta$, and $h_{l}^{\pm}, 1 \leq l \leq \sigma$
and the relations $h_{m}^{\pm}h_{l}^{\pm} =
h_{l}^{\pm}h_{l}^{\pm}, h_{l}^{\pm} h_{l}^{\mp} = 1,
\text{ for all } 1 \leq l,m \leq \sigma$, defining the free abelian group of
rank $\sigma$, and
\eqref{qls7}, the quantum Serre relations \eqref{qls9} for $i \neq j$ and $i \sim j$, \eqref{qls10},\eqref{qls11} (with $a_{i}$ instead of $x_{i}$), and the
lifted quantum Serre relations \eqref{linkSerre}.

\medbreak
If $(V,c)$ is a finite-dimensional braided vector space, we will say that
the {\it braiding is positive} if it is diagonal with matrix $(q_{ij})$,
and the scalars $q_{ii}$ are positive and different from 1,
for all $i$.

The next theorem follows from a result of Rosso \cite[Theorem 21]{Ro2} and the
theory described in the previous Sections.
\begin{theorem}\label{positive}\cite{AS5}
Let  $A$ be a pointed Hopf algebra with abelian group $\Gamma=G(A)$ and diagram $R$.
Assume that $R(1)$ has finite dimension and positive braiding.
Then the following are equivalent:

(a). $A$ is a domain of finite Gelfand-Kirillov dimension,
and the adjoint action of $G(A)$ on $A$ (or on $A_1$) is semisimple.

(b). The group $\Gamma$ is free abelian of finite rank, and
$$ A \simeq U(\mathcal{D}), \text{ where }\mathcal{D} \text{ is a positive
datum of finite Cartan type for } \Gamma.$$
\end{theorem}
\qed

\medbreak
It is likely that the positivity assumption on the infinitesimal braiding in the
last theorem is related to the existence of a compact involution.


\section{Pointed Hopf algebras of type $A_{n}$}\label{A_n}

\

In this Chapter, we develop from scratch, {\it i. e.} without using
Lusztig's results, the classification of all finite dimensional pointed Hopf
algebras whose infinitesimal braiding is of type $A_n$. The main results of this Chapter are
new.

\subsection{Nichols algebras of type $A_{n}$}\label{nich-A_n}

\

Let $N$ be an integer, $N >2$, and let $q$ be a root of 1 of order $N$.
For the case $N = 2$, see \cite{AnDa}.

\medbreak
Let $ q_{ij}$, $1\le i, j\le n$, be roots of 1 such that
\begin{equation}\label{defqij}
q_{ii} = q, \quad q_{ij}q_{ji} = \cases q^{-1}, \quad &\text{if } \vert
i-j\vert = 1, \\  1, \quad &\text{if } \vert
i-j\vert \ge 2. \endcases
\end{equation}
for all $1\le i, j\le n$.
For convenience, we denote
$$ B^{i,j}_{p, r} := \prod_{ i\le \ell \le j-1,
\, p\le h \le r-1}
q_{\ell, h},$$
for any $1\le i < j \le n + 1$,  $1\le p < r\le n +1$. Then we have the
following identities, whenever $i < s < j$, $p < t < r$:
\begin{align}\label{2.2.01}
B^{i,s}_{p, r} B^{s,j}_{p, r} &= \prod_{ i\le \ell \le s-1,  p\le h \le r-1}
q_{\ell, h}
\prod_{ s\le \ell \le j-1, p\le h \le r-1} q_{\ell,h} = B^{i,j}_{p, r};
\\ \label{2.2.02}
B_{p,t}^{i, j} B_{t,r}^{i, j} &= B_{p,r}^{i, j};
\end{align}  also,
 \begin{align}\label{2.2.03}
B^{i,j}_{j, j+1}B_{i,j}^{j, j+1} &= \prod_{i\le \ell \le j-1} q_{\ell, j}
\prod_{i\le h \le j-1} q_{j, h} = q^{-1}; \\
\label{2.2.04} B_{i,j}^{i,j} & = q.
 \end{align}

\medbreak
We consider in this Section a vector space $V = V_{n}$ with a basis $x_{1},
\dots,
x_{n}$ and braiding  determined by:
$$ c(x_{i}\otimes x_{j}) = q_{ij} \, x_{j}\otimes x_{i}, \quad
1\le i,j \le n;$$
that is, $V$ is of type $A_{n}$.

\begin{obs}\label{ochopuntouno}
Let  $\Gamma$ be a group, $g_{1}, \dots, g_{n}$ in the center of $\Gamma$,
and $\chi_{1}, \dots, \chi_{n}$ in $\VGamma$
such that
$$q_{ij} = \langle \chi_{j}, g_{i}\rangle,\quad 1 \leq i,j \leq n.$$
Then $V$ can be realized as a Yetter-Drinfeld module over $\Gamma$ by
declaring
\begin{equation}\label{structYD} x_{i}\in V^{\chi_{i}}_{g_{i}}, \qquad 1\le
i \le n.
\end{equation}
For example, we could consider $\Gamma = (\Z/P)^n$, where $P$ is divisible
by the orders of all the $q_{ij}$'s; and take $g_{1}, \dots, g_{n}$ as the
canonical basis of $\Gamma$.
\end{obs}

\medbreak
We shall consider a braided Hopf algebra $R$   provided with an inclusion of braided
vector spaces  $V\rightarrow P(R)$. We identify the elements $x_{1},
\dots, x_{n}$ with their images in $R$. Distinguished examples of such $R$ are
the tensor algebra $T(V)$ and the Nichols algebra $\toba (V)$.
Additional hypotheses on $R$ will be stated when needed.

\medbreak

We introduce the family $(e_{ij})_{1\le i < j \le n + 1}$
of elements of  $R$ as follows:
\begin{flalign} \label{2.2.1}
&e_{i, i+1} := x_{i};& \end{flalign}
\begin{flalign}     \label{2.2.2}
&e_{i, j} :=
[e_{i, j-1} , e_{j-1, j}]_{c},  \quad 1\le i < j \le n + 1, \, j-i \ge 2.&
\end{flalign}

The braiding between elements of this family is given by:
\begin{flalign}\label{2.2.3}
&c(e_{i, j}\otimes e_{p, r}) = B^{i,j}_{p,r}e_{p, r}\otimes e_{i, j},
\quad 1\le i < j \le n + 1, \, 1\le p < r\le n +1.&
\end{flalign}
\noindent In particular,
\begin{align*} e_{i, j} &=  e_{i, j-1} e_{j-1, j}  - B_{j-1, j}^{i, j-1}
e_{j-1, j} e_{i, j-1}.\end{align*}

\begin{obs}\label{ochopuntodos}
When $V$ is realized as a Yetter-Drinfeld module over $\Gamma$ as in Remark
\ref{ochopuntouno},
we have $e_{i, j} \in R^{\chi_{i,j}}_{g_{i,j}}$, where
\begin{equation}\label{generalgchi} \chi_{i,j} = \prod_{i \le \ell \le
j-1}\chi_{\ell},
\qquad g_{i,j} = \prod_{i \le \ell \le j-1} g_{\ell}, \qquad 1\le i < j \le
n + 1.
\end{equation}
\end{obs}

\begin{lema}\label{lema2.2.1}  (a). If $R$ is finite dimensional or  $R
\simeq\toba
(V)$, then   \begin{equation}\label{2.1.7} e_{i, i+1}^{N} =0, \quad \text{if} \quad 1\le i \le n.
\end{equation}
(b). Assume that $R \simeq\toba (V)$. Then
\begin{equation}\label{2.2.40} [e_{i,i+1}, e_{p,p+1}]_{c} = 0, \quad \text{that is} \quad e_{i,i+1} e_{p,p+1} = q_{ip} e_{p,p+1} e_{i,i+1},\quad \text{if} \quad 1\le i < p  \le n, p-i \ge 2.
\end{equation}
(c). Assume that $R \simeq\toba (V)$. Then
\begin{align}
[e_{i,i+1},
[e_{i,i+1}, e_{i+1,i+2}]_{c}]_{c} &=0,  \quad \text{if} \quad 1\le i < n;\label{2.2.400}\\
 \label{2.2.401} [e_{i+1,i+2},
[e_{i+1,i+2}, e_{i,i+1}]_{c}]_{c} &=0, \quad \text{if} \quad 1\le i < n;
\end{align}
that is
\begin{align}
\label{2.2.402}
 e_{i,i+1} e_{i,i+2} &= B^{i,i+1}_{i,i+2} e_{i,i+2} e_{i,i+1}, \\
\label{2.2.403}
e_{i,i+2} e_{i+1,i+2} &= B^{i,i+2}_{i+1,i+2} e_{i+1,i+2} e_{i,i+2}.
\end{align}
\end{lema}

\pf (a). This follows from Lemma \ref{primitivos} (a),
use $c(e_{i, i+1}^{N}\otimes e_{i, i+1}^{N}) = e_{i,
i+1}^{N}\otimes e_{i, i+1}^{N}$ in the finite dimensional case.

(b) and (c). By Lemma \ref{primitivos} (b), the elements
   $$[e_{i,i+1}, e_{p,p+1}]_{c}, \quad  [e_{i,i+1},
[e_{i,i+1}, e_{i+1,i+2}]_{c}]_{c}\quad  \text{and } [e_{i+1,i+2},
[e_{i+1,i+2}, e_{i,i+1}]_{c}]_{c}$$ are primitive. Since they are
homogeneous of degree
2, respectively of degree 3, they should be 0. To derive \eqref{2.2.403}
from \eqref{2.2.401}, use \eqref{brcomm1}.
 \epf

\begin{lema}\label{lema2.2.200}  Assume that \eqref{2.2.40} holds in $R$.
Then
\begin{flalign}
\label{2.2.4} &\qquad [e_{i,j}, e_{p,r}]_{c} = 0, \qquad \;\;\text{ if} \quad
1\le i < j < p < r \le n + 1.&\\
 \label{2.2.4bis}&\qquad [e_{p, r}, e_{i,j}]_{c} = 0,\qquad  \;\;\text{ if} \quad 1\le i < j < p < r \le n + 1.&\\
 \label{2.2.6}
 &\qquad[e_{i,p}, e_{p,j}]_{c} = e_{i,j}, \qquad  \text{if} \quad 1\le i < p < j \le n + 1.&
\end{flalign}
\end{lema}

\pf \eqref{2.2.4}. For $j = i+1$ and $r = p+1$, this is  \eqref{2.2.40};
the general case follows recursively using  \eqref{brcomm2}. \eqref{2.2.4bis} follows
from \eqref{2.2.4}, since
$B^{i,j}_{p,r}B_{i,j}^{p,r} = 1$ in this case.

\eqref{2.2.6}. By induction on $j - p$; if $p = j - 1$ then \eqref{2.2.6} is just
\eqref{2.2.2}.
For $p < j$, we have
\begin{equation*}
e_{i, j+1} = [e_{i,j}, e_{j, j+1}]_{c} =
[[e_{i,p}, e_{p,j}]_{c}, e_{j, j+1}]_{c}
= [e_{i,p}, [e_{p,j}, e_{j, j+1}]_{c}]_{c} = [e_{i,p}, e_{p,j + 1}]_{c}
\end{equation*}
by \eqref{brcomm2}, since $ [e_{i,p}, e_{j, j+1}]_{c} = 0$ by \eqref{2.2.4}.

\epf

\begin{lema}\label{lema2.2.2} Assume that \eqref{2.2.40} holds in $R$. Then for any $1\le i < j \le n + 1,
$\begin{equation}\label{2.2.5}
\Delta(e_{i,j}) = e_{i,j}\otimes 1 + 1\otimes
e_{i,j} +
(1 - q^{-1}) \sum_{i
< p < j} e_{i,p} \otimes e_{p,j}.
\end{equation}
\end{lema}

\pf We proceed by induction on $j-i$.  If $j-i = 1$, the formula
just tells that the $x_{i}$'s are primitive. For the inductive step, we
compute
\begin{align*}
\Delta(e_{i,j}e_{j, j+1}) &=
 \left(e_{i,j}\otimes 1 + 1\otimes e_{i,j} + (1 -
q^{-1}) \sum_{i < p < j}
e_{i,p} \otimes e_{p,j}\right)
(e_{j, j+1}\otimes 1 + 1\otimes e_{j,j+1})
\\ &= e_{i,j}e_{j, j+1}\otimes 1  + B^{i,j}_{j, j+1}  e_{j, j+1}\otimes
e_{i,j}
+ (1- q^{-1}) \sum_{i < p < j} B^{p,j}_{j,j+1} e_{i,p} e_{j, j+1} \otimes
e_{p,j}
\\ &\quad + e_{i,j}\otimes e_{j, j+1}  + 1\otimes e_{i,j} e_{j,j+1}
 + (1- q^{-1}) \sum_{i < p < j} e_{i,p}  \otimes e_{p,j} e_{j, j+1};
\end{align*}
and also
\begin{align*}
\Delta(e_{j, j+1} e_{i,j}) &=  (e_{j, j+1}\otimes 1 + 1\otimes e_{j,
j+1})  \left(e_{i,j}\otimes 1 + 1\otimes e_{i,j} + (1- q^{-1}) \sum_{i < p <
j}
e_{i,p} \otimes e_{p,j}\right)\\
&= e_{j, j+1} e_{i,j} \otimes 1 + e_{j, j+1}\otimes e_{i,j} + (1- q^{-1})
\sum_{i < p < j} e_{j, j+1} e_{i,p} \otimes e_{p,j}
\\ &\quad + B^{j, j+1}_{i,j} e_{i,j} \otimes e_{j, j+1} +1 \otimes e_{j,
j+1}
e_{i,j}  + (1- q^{-1})  \sum_{i < p < j} (B_{j, j+1}^{i,p})^{-1} e_{i,p}
\otimes e_{j, j+1}e_{p,j}. \end{align*}
Hence
\begin{align*}
\Delta(e_{i, j+1}) &=
 e_{i, j+1} \otimes 1 + 1 \otimes e_{i, j+1} + (1 - B^{ij}_{j, j+1}
B_{ij}^{j,
j+1}) e_{i,j} \otimes e_{j, j+1} \\
&\quad  + (1- q^{-1})
\sum_{i < p < j} \left(B^{p,j}_{j, j+1}B^{i,p}_{j,j+1} - B^{i,j}_{j,j
+1}\right)e_{j, j+1} e_{i,p} \otimes e_{p,j} \\ &\quad
+ (1- q^{-1})  \sum_{i < p < j}  e_{i,p} \otimes
\left(e_{p,j}e_{j, j+1} - B_{j, j+1}^{i,j} (B_{j, j+1}^{i,p})^{-1} e_{j,
j+1}e_{p,j}
\right)
\\ &=
 e_{i, j+1} \otimes 1 + 1 \otimes e_{i, j+1}
+ (1- q^{-1})
\sum_{i < p < j + 1}  e_{i,p} \otimes e_{p,j+1};
\end{align*}
by \eqref{2.2.03}, \eqref{2.2.4} and the hypothesis. \epf

\begin{obs}\label{ochopuntotres}
Let  $\Gamma$ be a group with $g_{1}, \dots, g_{n}$ in the center of
$\Gamma$,
$\chi_{1}, \dots, \chi_{n}$ in $\VGamma$, as in \ref{ochopuntouno}. Let
$R$ be a braided Hopf algebra in $\ydg$ such that \eqref{2.2.40} holds in
$R$.
It follows from \eqref{2.2.5} and the reconstruction formulas
for the bosonization  \eqref{smash1} that
\begin{equation}\label{boson-delta}
\Delta_{R\# \ku \Gamma}(e_{i, j}) = e_{i, j} \otimes 1 +
g_{i, j}  \otimes
e_{i, j} + (1 - q^{-1}) \sum_{i < p < j} e_{i,p} g_{p, j} \otimes e_{p, j}.
\end{equation}
\end{obs}


\begin{lema}\label{lema2.2.201}
 Assume that \eqref{2.2.40}, \eqref{2.2.400}, \eqref{2.2.401} hold in $R$. Then
\begin{flalign}\label{2.2.41} &\qquad [e_{i,j}, e_{p,r}]_{c} = 0,
\qquad\text{ if } \quad
1\le i   <  p < r < j \le n + 1;& \\
\label{2.2.42} &\qquad [e_{i,j}, e_{i,r}]_{c} = 0, \qquad\text{ if } \quad
1\le i  < j  < r\le n + 1;& \\
\label{2.2.43} &\qquad [e_{i,j}, e_{p, j}]_{c} = 0, \qquad\text{ if } \quad
1\le i < p < j \le n + 1.&
\end{flalign}
\end{lema}
\pf (a).  We prove  \eqref{2.2.41} by induction on $j -i$. If $j -i = 3$,
then
$$
[e_{i,i + 3}, e_{i+1,i+2}]_{c} = [[e_{i,i + 2}, e_{i+2,i+3}]_{c},
e_{i+1,i+2}]_{c}
= [[[e_{i,i + 1}, e_{i+1,i+2}]_{c}, e_{i+2,i+3}]_{c}, e_{i+1,i+2}]_{c} = 0,
$$
by Lemma \ref{brcomm-an}. If $j -i > 3$ we argue by induction on $r- p$.
If $r- p = 1$, then there exists an index $h$ such that either $ i < h  <  p
< r = p+1 < j$
or $ i   <  p < r = p+1 < h < j$. In the first case, by \eqref{2.2.6}, we
have
$$
[e_{i,j}, e_{p,r}]_{c} = [[e_{i,h}, e_{h,j}]_{c}, e_{p,p+1}]_{c} = 0;$$
the last equality follows from Lemma \ref{brcomm} (c), because
of  \eqref{2.2.4} and the induction hypothesis. In the second case, we have
$$
[e_{i,j}, e_{p,r}]_{c} = [[e_{i,h}, e_{h,j}]_{c}, e_{p,p+1}]_{c} = 0;$$
the last equality follows from Lemma \ref{brcomm} (c), because
of the induction hypothesis and \eqref{2.2.4bis}. Finally, if $r- p > 1$
then
$$
[e_{i,j}, e_{p,r}]_{c} = [e_{i,j}, [e_{p,r-1}, e_{r-1,r}]_{c}]_{c} = 0$$
by Lemma \ref{brcomm} (b) and the induction hypothesis.

\medbreak
(b). We prove  \eqref{2.2.42} by induction on $r -i$. If $r -i = 2$, then
the claimed equality
is just \eqref{2.2.400}. If $r -i > 2$ we argue by induction on $j-i$.
If $j-i = 1$ we have $$ [e_{i,i +1}, e_{i, r}]_{c} =
[e_{i,i +1}, [e_{i,r-1}, e_{r-1, r}]_{c}]_{c} = [[e_{i,i +1},
e_{i,r-1}]_{c}, e_{r,r-1}]_{c} = 0
$$
by \eqref{brcomm2}, since $ [e_{i,i +1}, e_{r-1, r}]_{c} = 0$ by
\eqref{2.2.4}.
If $j -i > 2$, we have
$$ [e_{i,j}, e_{i, r}]_{c} =
[[e_{i,j-1}, e_{j-1, j}]_{c}, e_{i, r}]_{c}  = 0
$$
by Lemma \ref{brcomm} (c), because
of the induction hypothesis and \eqref{2.2.41}.

\medbreak

The proof of \eqref{2.2.43} is analogous to the proof of \eqref{2.2.42},
using \eqref{2.2.401} instead of \eqref{2.2.400}.
\epf

\begin{lema}\label{lema2.2.201d}
 Assume that \eqref{2.2.40}, \eqref{2.2.400}, \eqref{2.2.401} hold in $R$.
Then
\begin{equation}\label{2.2.7}
[e_{i,j}, e_{p,r}]_{c} = B^{pj}_{jr}(q-1) e_{ir}e_{pj}, \qquad\text{ if }
\quad
1\le i < p < j < r\le n + 1. \end{equation}
\end{lema}

\pf We compute:
\begin{align*}
[e_{i,j}, e_{p,r}]_{c} &=  [[e_{i,p}, e_{p,j}]_{c}, e_{p,r}]_{c}\\
&= [e_{i,p}, [e_{p,j}, e_{p,r}]_{c}]_{c} + B_{p,r}^{p,j} [e_{i,p},
e_{p,r}]_{c}e_{p,j}  -
B_{p,j}^{i,p} e_{p,j} [e_{i,p}, e_{p,r}]_{c}
 \\&=  B_{p,r}^{p,j} e_{i,r} e_{p,j}  - B_{p,j}^{i,p}  e_{p,j} e_{i,r}
\\&= \left(B_{p,r}^{p,j} - B_{p,j}^{i,p} (B_{p,j}^{i,r})^{-1} \right)
e_{i,r} e_{p,j}
\\&= \left(B_{p,r}^{p,j} - ( B_{p,j}^{p,r} )^{-1}  \right) e_{i,r} e_{p,j}
\\&= B^{pj}_{jr}(q-1) e_{ir}e_{pj}.
 \end{align*}
Here, the first equality is by \eqref{2.2.6}; the second, by Lemma
\ref{brcomm} (a);
the third,  by \eqref{2.2.42} and  by \eqref{2.2.6}; the fourth, by
\eqref{2.2.41}.
\epf

\begin{lema}\label{comul-powrootvec} Assume that \eqref{2.2.40},
\eqref{2.2.400} and \eqref{2.2.401} hold in $R$.
For any $1\le i < j \le n + 1$ we have
\begin{equation}\label{2.2.500} \Delta(e_{i,j}^N) =
e_{i,j}^N\otimes 1 + 1\otimes e_{i,j}^N + (1 - q^{-1})^N \sum_{i
< p < j} \left(B_{i,p}^{p,j}\right)^{N(N-1)/2} e_{i,p}^N \otimes e_{p,j}^N.
\end{equation}\end{lema}
\pf By \eqref{2.2.5}, and using several times the quantum binomial formula
\eqref{qbfpwr},
we have
\begin{align*}
\Delta(e_{i,j}^N) &= \left(e_{i,j}\otimes 1  + (1 - q^{-1}) \sum_{i
< p < j} e_{i,p} \otimes e_{p,j}\right)^N + 1\otimes e_{i,j}^N \\
&= e_{i,j}^N\otimes 1  + (1 - q^{-1})^N  \left(\sum_{i
< p < j} e_{i,p} \otimes e_{p,j}\right)^N + 1\otimes e_{i,j}^N
 \\&=
e_{i,j}^N\otimes 1  + (1 - q^{-1})^N  \sum_{i
< p < j} \left(e_{i,p} \otimes e_{p,j} \right)^N + 1\otimes e_{i,j}^N\\&=
e_{i,j}^N\otimes 1  + (1 - q^{-1})^N  \sum_{i
< p < j} \left(B_{i,p}^{p,j}\right)^{N(N-1)/2} e_{i,p}^N \otimes e_{p,j}^N +
1\otimes e_{i,j}^N.
\end{align*}
Here, in the first equality we use that $(1\otimes e_{i,j})(e_{i,j}\otimes
1) =
q(e_{i,j}\otimes 1)(1\otimes e_{i,j})$  and $(1\otimes e_{i,j})(e_{i,p}
\otimes e_{p,j}) =
q(e_{i,p} \otimes e_{p,j})(1\otimes e_{i,j})$, this last by \eqref{2.2.43};
in the second, we use $(e_{i,p} \otimes e_{p,j}) (e_{i,j} \otimes 1) = q
(e_{i,j} \otimes 1)
(e_{i,p} \otimes e_{p,j})$, which follows from \eqref{2.2.42}; the third,
that
$(e_{i,p} \otimes e_{p,j}) (e_{i,s} \otimes e_{s,j}) = q^2 (e_{i,s} \otimes
e_{s,j})
(e_{i,p} \otimes e_{p,j})$ for $p<s$, which follows from \eqref{2.2.42} and
\eqref{2.2.43};
the fourth, from $(e_{i,p} \otimes e_{p,j})^h =
\left(B_{i,p}^{p,j}\right)^{h(h-1)/2}
e_{i,p}^h \otimes e_{p,j}^h$.
\epf

\begin{obs}\label{ochopuntocuatro}
Let  $\Gamma$ be a group with $g_{1}, \dots, g_{n}$ in the center of
$\Gamma$,
$\chi_{1}, \dots, \chi_{n}$ in $\VGamma$, as in \ref{ochopuntouno}. Let
$R$ be a braided Hopf algebra in $\ydg$ such that \eqref{2.2.40},
\eqref{2.2.400} and \eqref{2.2.401} hold in $R$.
By \eqref{2.2.500} and the reconstruction formulas  \eqref{smash1}, we have
\begin{equation}\label{boson-deltapow}
\Delta_{R\# \ku \Gamma}(e_{i, j}^N) = e_{i, j}^N \otimes 1 +
g_{i, j}^N  \otimes
e_{i, j}^N + (1 - q^{-1})^N \sum_{i < p < j}
\left(B_{i,p}^{p,j}\right)^{N(N-1)/2}
e_{i,p}^N g_{p, j}^N \otimes e_{p, j}^N.
\end{equation}
\end{obs}

\begin{lema}\label{powrootvec0} Assume that  $R = \toba(V)$.
Then
\begin{equation}\label{2.2.501} e_{i,j}^N = 0, \qquad
1\le i < j \le n + 1. \end{equation}\end{lema}
\pf This follows from Lemma \ref{comul-powrootvec} by induction
on $j-i$, the case  $j-i = 1$ being Lemma \ref{lema2.2.1} (c). \epf

\begin{lema}\label{Lemma2.2.7} Assume that \eqref{2.2.40}, \eqref{2.2.400},
\eqref{2.2.401},
\eqref{2.2.501}, hold in $R$. Assume, furthermore, that $R$ is generated as
an algebra by the elements $x_{1}, \dots, x_{n}$.
Then the algebra $R$ is spanned as a vector space by the elements
\begin{equation}\label{2.2.10}e_{1,2}^{\epsilon_{1,2}}
e_{1,3}^{\epsilon_{1,3}}
\dots e_{1,n + 1}^{\epsilon_{1,n + 1}} e_{2,3}^{\epsilon_{2,3}}
\dots e_{2,n + 1}^{\epsilon_{2,n +1}} \dots
e_{n,n +1}^{\epsilon_{n,n + 1}}, \qquad \text{with }
\epsilon_{i,j}\in \{0,1, \dots, N -1\}. \end{equation}\end{lema}

\pf  We order the family $(e_{ij})$  by   $$e_{1,2} \prec e_{1,3}
\prec  \dots e_{1,n + 1} \prec  e_{2,3} \prec  \dots e_{2,n
+ 1} \prec \dots e_{n,n +1};$$ this induces an ordering in
the monomials \eqref{2.2.10}. If $M$ is an ordered monomial, we
set $\sigma(M) := e_{r,s}$ if $e_{r, s}$ is the first element
appearing in $M$. Let $B$ be the subspace generated by the
monomials in \eqref{2.2.10}.  We show by induction on the length
that, for any ordered monomial $M$ and for any $i$, $e_{i, i+1} M
\in B$ and it is 0 or a combination of monomials $N$ with
$\sigma(N) \succeq \min \{e_{i, i+1}, \sigma(M)\}$, length of $N
\le$ length of $M + 1$. The statement is evident if the length of
$M$ is 0; so that assume that the length is positive. Write $M
=
e_{p,q} M'$ where $e_{p,q} \preceq M'$. We have several cases:

\medbreak

If $i < p$ or $i = p$ and $i+1 < q$, $e_{i, i+1} \prec e_{p,q}$ and we
are done.

If $i = p$ and $i+1 = q$, then the claim is clear.

 If $p < q < i$ then $e_{i, i+1} e_{p,q} = B^{i, i+1}_{p,q} \,
e_{p,q}e_{i, i+1}$ by
\eqref{2.2.4bis}; hence $e_{i, i+1} M = e_{i, i+1} e_{p,q}M' =
B^{i, i+1}_{p,q} \,e_{p,q}e_{i, i+1} M'$;
by the inductive hypothesis and the fact that $e_{p,q} \preceq \min \{e_{i,
i+1}, \sigma(M')\}$, the claim follows.

 If $p < i = q$ then $e_{i, i+1} e_{p,i} =
(B^{p,i}_{i,i+1})^{-1}(e_{p,i} e_{i,i+1}
- e_{p, i+1})$ by \eqref{2.2.6}; again, the inductive hypothesis and
$e_{p,q} \preceq
\min \{e_{i, i+1}, \sigma(M')\}$ imply that $e_{p,i} e_{i,i+1} M'$ has the form we want. To see that $e_{p,i+1} M"$ satisfies the claim when $e_{p,i} = \sigma(M')$, we use $e_{p,i+1} e_{p,i} = (B_{p,i+1}^{p,i})_{-1} e_{p,i} e_{p,i+1}$ by \eqref{2.2.42}.

 If $p < i < q$ then $e_{i, i+1} e_{p,q} =
(B^{p,q}_{i,i+1})^{-1}e_{p,q}e_{i, i+1}$ by
\eqref{2.2.41} or \eqref{2.2.43}; we then argue as in the two preceding
cases.

\medbreak

Therefore, $B = R$ since it is a left ideal containing 1.
\epf

We shall say that the elements $e_{1,2}, e_{1,3}, \dots,
e_{1,n +
1} e_{2,3} \dots e_{2,n + 1}\dots e_{n,n +1} $, in this order,
form a {\it PBW-basis} for $R$ if the monomials \eqref{2.2.10}
form a basis of $R$. Then we have:

\begin{theorem}\label{Theorem2.2.8}The elements $e_{1,2}, e_{1,3}, \dots,
e_{1,n +
1} e_{2,3} \dots e_{2,n + 1}\dots e_{n,n +1} $, in this order,
form a PBW basis for $\mathfrak B(V_{n})$. In particular,
$$\dim \mathfrak B(V_{n}) = N^{\frac{n(n + 1)}2}.$$
\end{theorem}
\pf We proceed by induction on $n$. The case $n = 1$ is clear, see
\cite[Section 3]{AS1}
for details.
We assume the statement for $n - 1$. We consider $V_{n}$ as a
Yetter-Drinfeld module over
$\Gamma =(\Z/P)^n$, as explained in Remark \ref{ochopuntouno}.
Let $Z_{n}= \mathfrak B(V_{n}) \# \ku\Gamma$. Let $i_{n}: V_{n - 1}
\to
V_{n}$ be given by $x_{i} \mapsto x_{i}$ and $p_{n}: V_{n} \to
V_{n - 1}$ by $x_{i} \mapsto x_{i}$, $1\le i\le n - 1$ and
$x_{n}\mapsto 0$. The splitting of Yetter-Drinfeld modules $\id_{V_{n
- 1}} = p_{n} i_{n}$ gives rise to a splitting of Hopf algebras
$\id_{Z_{n - 1}} = \pi_{n} \iota_{n}$, where $\iota_{n}: Z_{n
-
1} \to Z_{n}$  and $\pi_{n}: Z_{n} \to Z_{n - 1}$ are
respectively
induced by $i_{n}$, $p_{n}$. Let $$R_{n} =  Z_{n}^{\text{co }
\pi_{n}} = \{z\in Z_{n}: (\id \otimes \pi_{n})\Delta(z) = z\otimes
1\}. $$

Then $R_{n}$ is a braided Hopf algebra in the category
${}_{Z_{n - 1}}^{Z_{n - 1}}\mathcal{YD}$;
we shall denote by $c_{R_{n}}$ the corresponding braiding of $R_{n}$.
We have   $ Z_{n} \simeq R_{n}\#Z_{n - 1}$ and in particular
$\dim  Z_{n} = \dim  R_{n} \dim Z_{n - 1}$.

\medbreak

For simplicity, we denote $h_{i} = e_{i, n + 1}$, $1\le i \le n$.
We have
$h_{i} h_{j} = B^{i, n+1}_{j, n+1}h_{j}h_{i}$, for $i < j$, by
\eqref{2.2.43}.
We claim that
 $h_{1}, \dots, h_{n}$ are linearly
independent primitive elements of the braided Hopf algebra $R_{n}$.

\medbreak
Indeed, it follows from \eqref{2.2.2} that $\pi_{n}(h_{i}) = 0$; by
\eqref{2.2.5}, we conclude that $h_{i}\in R_{n}$.
We prove by induction on $j = n + 1 -i$ that $h_{i}$ is a
primitive element of $R_{n}$, the case $j = 1$ being clear. Assume the
statement for $j$.  Now
\begin{align*}
x_{i-1}\rightharpoonup h_{i} &=  x_{i-1}h_{i} + g_{i-1}h_{i} \Ss (x_{i-1})
=
x_{i-1}h_{i} -g_{i-1}h_{i}g_{i-1}^{-1}x_{i-1} \\
&= x_{i-1}h_{i} - B^{i-1, i}_{i, n+1} h_{i}x_{i-1} = [x_{i-1}, h_{i}]_{c} =
h_{i-1}.\end{align*}
So
\begin{align*}
\Delta_{R_{n}}(h_{i-1}) &=  \Delta_{R_{n}}(x_{i-1}\rightharpoonup
h_{i}) = x_{i-1}\rightharpoonup \Delta (h_{i}) \\
&=  g_{i-1}
\rightharpoonup 1 \otimes  x_{i-1}\rightharpoonup h_{i} + x_{i-1}
\rightharpoonup h_{i}\otimes  1 = h_{i-1}\otimes 1
+ 1\otimes h_{i-1}. \end{align*}
We prove also by induction on $j = n + 1 -i$ that $h_{i} \neq 0$ using
\eqref{2.2.5} and the
induction hypothesis on $Z_{n - 1}$. Since $h_{i}$ is homogeneous of
degree $j$ (with respect to the grading of $Z_{n}$),
we conclude that $h_{1}, \dots, h_{n}$ are linearly
independent.

\medbreak

We next claim
that $c_{R_{n}}(h_{i}\otimes h_{j}) = B^{i, n+1}_{j, n+ 1} h_{j}\otimes
h_{i}$, for any $i>j$.

By \eqref{boson-deltapow}, the coaction of $Z_{n - 1}$ on $R_{n}$
satisfies
$$
\delta(h_{i}) = g_{i, n + 1} \otimes e_{i, n + 1}
+ (1 - q^{-1}) \sum_{i < p < n + 1} e_{i,p} g_{p, n + 1}  \otimes
e_{p, n + 1}.
$$
If $j < i$, we compute the action on $R_{n}$:
\begin{align*}
e_{i, p} \rightharpoonup h_{j} &=  e_{i, p} h_{j} +
g_{i, p}  h_{j} \Ss (e_{i, n + 1}) + (1 - q^{-1}) \sum_{i < t < p} e_{i, t}
g_{t, p}
 h_{j} \Ss(e_{t, p}) \\
&= (B_{i,p}^{j, n+1})^{-1} h_{j} e_{i, p}  + B^{i,p}_{j, n+1} h_{j}
g_{i, p}   \Ss (e_{i, n + 1}) \\ & \qquad+  h_{j} (1 - q^{-1}) \sum_{i < t <
p}
B^{t,p}_{j, n+1} (B_{i,t}^{j, n+1})^{-1} e_{i, t} g_{t, p}   \Ss(e_{t, p})
\\ &= B^{i,p}_{j, n+1} h_{j} {e_{i, p}}_{(1)} \Ss({e_{i, p}}_{(2)})
\\ &=  0,
\end{align*}
by \eqref{2.2.43}. Thus
$$
c_{R_{n}}(h_{i}\otimes h_{j}) = g_{i, n + 1} \rightharpoonup
h_{j}\otimes h_{i} = B^{i, n+1}_{j, n+ 1} h_{j}\otimes h_{i}.$$

\medbreak

We next claim
that the dimension of the subalgebra of $R_{n}$ spanned by $h_{1}, \dots,
h_{n}$ is $\ge N^{n}$.

We already know that $$\Delta
(h_{j}^{m_{j}}) = \sum_{0 \le i_{j} \le m_{j}}
{\binom{m_{j}}{i_{j}}}_{q}h_{j}^{i_{j}}\otimes  h_{j}^{m_{j}-i_{j}}, \qquad
m_j \le N.$$

Let us denote $\mathbf m = (m_{1}, \dots, m_{j}, \dots , m_{n})$,
$\mathbf 1 = (1, \dots,1, \dots , 1)$, $\mathbf N = (N, \dots , N)$. We
consider the partial order
$\mathbf i \le \mathbf m$, if $i_{j} \le m_{j}$, $j = 1, \dots, n$. We set
$h^{\mathbf m} := h_{n}^{m_{n}} \dots h_{j}^{m_{j}} \dots h_{1}^{m_{1}}$.
>From the preceding claim, we deduce   that  $$\Delta
(h^{\mathbf m}) = h^{\mathbf m}\otimes 1 + 1\otimes h^{\mathbf m} +
\sum_{ 0 \le \mathbf i \le \mathbf m,\quad 0 \ne \mathbf i \ne \mathbf m}
c_{\mathbf m, \mathbf i}h^{\mathbf i}\otimes h^{\mathbf m-\mathbf i}, \quad
\mathbf m \le\mathbf N - \mathbf 1;$$ where
$c_{\mathbf m, \mathbf i} \ne 0$ for all $\mathbf i$.
We then argue recursively as in the proof of \cite[Lemma 3.3]{AS1} to
conclude
that the elements $h^{\mathbf m}$,
$\mathbf m \le\mathbf N - \mathbf 1$, are linearly independent; hence the
dimension of the subalgebra of $R_{n}$ spanned by
$h_{1}, \dots, h_{n}$ is $\ge N^{n}$, as claimed.

\medbreak

 We can now finish the proof of the Theorem.
Since $\dim Z_{n} \le N^{\frac{n(n + 1)}2}$ by Lemma \ref{Lemma2.2.7} and
$\dim
Z_{n - 1} =  N^{\frac{n(n - 1)}2}$ by the induction hypothesis, we
have $\dim R_{n}\le N^{n}$. By what we have just seen, this
dimension is exactly $N^{n}$.  Therefore, $\dim  Z_{n} = N^{\frac{n(n +
1)}2}$; in presence of Lemma \ref{Lemma2.2.7}, this implies the Theorem.
\epf

\begin{theorem}\label{nichalg-typeA_n} The Nichols algebra $\toba(V)$
can be presented by   generators $e_{i, i+1}$, $1\le i \le n$,
and relations \eqref{2.2.40},
\eqref{2.2.400}, \eqref{2.2.401} and \eqref{2.2.501}.
\end{theorem}
\pf Let $\mathfrak B'$ be the algebra presented by generators
$e_{i, i+1}$, $1\le i \le n$,
and relations \eqref{2.2.40},
\eqref{2.2.400}, \eqref{2.2.401} and \eqref{2.2.501}.
We claim that $\mathfrak B'$ is
 is a braided Hopf algebra with the $e_{i, i+1}$'s
primitive. Indeed, the claim follows
without difficulty; use Lemma \ref{comul-powrootvec} for relations
\eqref{2.2.501}.

By  Lemma \ref{Lemma2.2.7}, we see that the monomials \eqref{2.2.10}
span $\mathfrak B'$ as a vector space, and in particular that $\dim \mathfrak
B' \le N^{\frac{n(n + 1)}2}$. By Lemmas \ref{lema2.2.1} and
\ref{powrootvec0},
there is a surjective algebra map $\psi: \mathfrak B' \to \mathfrak
B(V)$. By Theorem \ref{Theorem2.2.8}, $\psi$ is an isomorphism. \epf

\medbreak
\subsection{Lifting of Nichols algebras of type $A_{n}$}\label{lift-A_n}

\

We fix in this Section a finite {\it abelian} group $\Gamma$ such that our
braided vector
space $V$ can be realized in $\ydg$, as in Remark \ref{ochopuntouno}.
That is, we have
$g_{1}, \dots, g_{n}$ in $\Gamma$,  $\chi_{1}, \dots,
\chi_{n}$ in $\VGamma$, such that $q_{ij} = \langle \chi_{j}, g_{i}\rangle$
for all i,j, and $V$ can be realized as a Yetter-Drinfeld module over $\Gamma$ by
\eqref{structYD}.

\medbreak
We also fix a finite dimensional pointed Hopf algebra $A$ such that
$G(A)$ is isomorphic to $\Gamma$, and the infinitesimal braiding  of $A$
is isomorphic to $V$ as a Yetter-Drinfeld module over $\Gamma$.
That is, $\gr A \simeq R \# \ku \Gamma$, and the subalgebra $R'$ of $R$
generated by $R(1)$ is isomorphic to $\toba (V)$.
We choose elements $a_{i} \in (A_{1})^{\chi_{i}}_{g_{i}}$
such that $\pi(a_{i}) = x_{i}$, $1\le i \le n$.

\medbreak
We shall consider, more generally, Hopf algebras $H$ provided with

\begin{itemize} \item a group isomorphism $\Gamma \to G(H)$;

\medbreak
\item  elements $a_{1}, \dots, a_{n}$ in $\SkPr(H)^{\chi_{i}}_{g_{i}, 1}$.
\end{itemize}
Further hypotheses on $H$ will be stated when needed. The examples of such
$H$ we are thinking of are the Hopf algebra $A$, and any bosonization $R \#
\ku \Gamma$,
where $R$ is any braided Hopf algebra in $\ydg$ provided with a monomorphism
of Yetter-Drinfeld
modules $V \to P(R)$; so that $a_{i} := x_{i} \# 1$, $1\le i \le n$.
This includes notably the Hopf algebras $T(V) \# \ku \Gamma$,
$\wtoba (V) \# \ku \Gamma$, $\toba (V) \# \ku \Gamma$.

Here $\wtoba(V)$ is the braided Hopf algebra in $\ydg$ generated by
$x_1, \dots, x_{\theta}$ with relations \eqref{2.2.40},
\eqref{2.2.400} and  \eqref{2.2.401}.

\medbreak

We introduce inductively the following elements of $H$:
\begin{flalign} \label{2.2.1lift}
&E_{i, i+1} := a_{i};& \end{flalign}
\begin{flalign}     \label{2.2.2lift}
&E_{i, j} :=
\ad (E_{i, j-1}) (E_{j-1, j}),  \quad 1\le i < j \le n + 1, \, j-i \ge 2.&
\end{flalign}
Assume that $H = R\# \ku \Gamma$ as above. Then, by the relations
between braided commutators and the adjoint \eqref{bradj}, the
relations \eqref{2.2.40}, \eqref{2.2.400} and  \eqref{2.2.401} translate
respectively to
\begin{align}\label{Anserre1} \ad E_{i,i+1} (E_{p,p+1}) &= 0;
\qquad 1\le i < p  \le n,\qquad p-i \ge 2;
\\
\label{Anserre2} (\ad E_{i,i+1})^2 (E_{i+1,i+2}) &=0, \qquad 1\le i < n;\\
 \label{Anserre3} (\ad E_{i+1,i+2})^2 (E_{i,i+1}) &=0, \qquad 1\le i < n.
\end{align}
\begin{obs}\label{ochopuntocinco}
Relations \eqref{Anserre1}, \eqref{Anserre2} and  \eqref{Anserre3}
can be considered, more generally, in any $H$ as above.
If these relations hold in $H$, then we have a Hopf algebra map
$\pi_H: \wtoba (V) \# \ku \Gamma \to H$. On the other hand, we know by
Remark
\ref{ochopuntocuatro} that the comultiplication of the elements
$E_{ij}^N$ is given by
\eqref{boson-deltapow}. Hence, the same formula is valid in $H$, provided
that relations
\eqref{Anserre1}, \eqref{Anserre2} and  \eqref{Anserre3} hold in it.
In particular, the subalgebra of H generated by the elements $E_{ij}^N$,
$g_{i, j}^N$,
$1\le i < j\le n +1$, is a Hopf subalgebra of $H$.
\end{obs}

\begin{lema}\label{valenenA} Relations
\eqref{Anserre1}, \eqref{Anserre2} and  \eqref{Anserre3} hold in $A$ if $N>3$.

\end{lema}

\pf This is a particular case of Theorem \ref{QSR}; we include the proof for
completeness.
We know, by Lemma \ref{twistnichols}, that
\begin{align*}  \ad E_{i,i+1} (E_{p,p+1}) &\in {\mathcal
P}_{g_{i}g_{p},1}(A)^{\chi_i\chi_p},
\qquad 1\le i < p  \le n,\qquad p-i \ge 2,
\\  (\ad E_{i,i+1})^2 (E_{p,p+1})
&\in {\mathcal P}_{g_{i}^{2}g_{p},1}(A)^{\chi_{i}^{2}\chi_{p}},
\qquad 1\le i, p  \le n,\qquad \vert p-i\vert = 1.
\end{align*}

Assume  that $ \ad E_{i,i+1} (E_{p,p+1})\neq 0$, and $\chi_i\chi_p \neq \varepsilon$, where  $1\le i < p  \le n$, $p-i
\ge 2$.
By Lemma \ref{start}, there exists $\ell$, $1\le \ell \le n$,
such that $g_{i} g_{p} = g_{\ell}$, $\chi_{i} \chi_{p} = \chi_{\ell}$. But
then
\begin{equation*}
q =
\chi_{\ell}(g_{\ell}) = \chi_{i}(g_{i}) \chi_{i}(g_{p}) \chi_{p}(g_{i})
\chi_{p}(g_{p}) =q^2.
\end{equation*}
Hence $q= 1$, a contradiction.

Assume next that $ \ad E_{i,i+1}^2 (E_{p,p+1})\neq 0$,   $\vert p-i\vert =
1$.\, and $\chi_{i}^{2} \chi_{p} \neq \varepsilon$. By Lemma \ref{start} , there exists $\ell$, $1\le \ell \le n$,
such that $g_{i}^2 g_{p} = g_{\ell}$, $\chi_{i}^2 \chi_{p} = \chi_{\ell}$.
But then
\begin{equation*}
q =
\chi_{\ell}(g_{\ell}) = \chi_{i}(g_{i})^4 \chi_{i}(g_{p})^2
\chi_{p}(g_{i})^2 \chi_{p}(g_{p}) = q^3.
\end{equation*}
Hence $q= \pm 1$, a contradiction (we assumed $N >2$).

It remains to exclude the cases $\chi_i\chi_p =\varepsilon, |p-i| \geq 2$, and
 $\chi_{i}^{2} \chi_{p} = \varepsilon, |p-i|=1$. The first case leads to the contradiction $N=3$. In the second case it follows from the connectivity of $A_{n}$ that $N$ would divide 2 which is also impossible.
\epf

\begin{lema}\label{vivenkgamma} If $H = A$, then
$E_{i,j}^N \in \ku \Gamma^N$, for any $1\le i < j \le n + 1$.
\end{lema}
\pf We first show that $E_{i,j}^N \in \ku \Gamma$, $1\le i < j \le n + 1$.
(For our further purposes, this is  what we really need).

Let $i<j$. We claim that there exists no $\ell$, $1\le \ell \le n$,
such that $g_{i, j}^N = g_{\ell}$, $\chi_{i, j}^N = \chi_{\ell}$. Indeed,
otherwise
we would have
\begin{equation*}
q = \chi_{\ell}(g_{\ell}) = \chi_{i, j}(g_{i, j})^{N^2} =q^{N^2} = 1.
\end{equation*}

By Lemma \ref{valenenA} and Remark \ref{ochopuntocinco}, we have
\begin{equation}\label{coproduct}
\Delta(E_{i, j}^N) = E_{i, j}^N \otimes 1 +
g_{i, j}^N  \otimes
E_{i, j}^N + (1 - q^{-1})^N \sum_{i < p < j}
\left(B_{i,p}^{p,j}\right)^{N(N-1)/2}
E_{i,p}^N g_{p, j}^N \otimes E_{p, j}^N.
\end{equation}

We proceed by induction on $j-i$.
If $j-i=1$, then, by Lemma \ref{start} ,
 either $E_{i,i+1}^N \in \ku \Gamma$ or $E_{i,i+1}^N \in {\mathcal
P}_{g_{i}^N,1}(A)^{\chi_{i}^N}$ and $\chi_{i}^{N} \neq \varepsilon$, hence $g_{i}^{N} = g_{l}, \chi_{i}^{N}= \chi_{l}$ for some $l$; but this last possibility contradicts the
claim above.
Assume then that $j-i > 1$. By the induction hypothesis,
$\Delta(E_{i, j}^N) = E_{i, j}^N \otimes 1 +
g_{i, j}^N  \otimes
E_{i, j}^N + u$, for some $u \in \ku \Gamma \otimes \ku \Gamma$.
In particular, we see that $E_{i, j}^N   \in (A_{1})^{\chi_{i}^N}$.
Then, by Lemma \ref{start},
 either $\chi_{i}^N = \varepsilon$ and hence $E_{i,i+1}^N \in \ku \Gamma$,
or else $\chi_{i}^N \neq \varepsilon$, which implies $u= 0$ and
$E_{i,i+1}^N \in {\mathcal P}_{g_{i}^N,1}(A)^{\chi_{i}^N}$. Again,
this last possibility contradicts the claim above.

Finally, let $C$ be the subalgebra of $A$ generated by the elements
$E_{ij}^N$, $g_{i, j}^N$,
$1\le i < j\le n +1$, which is a Hopf subalgebra of $H$.
Since $E_{i,j}^N \in \ku \Gamma \cap C$, we conclude that $E_{i,j}^N \in C_0
= \ku \Gamma^N$.
\epf

To solve the lifting problem, we see from Lemma \ref{vivenkgamma} that 
we first have to answer a combinatorial question in the group algebra 
of an abelian group. To simplify the notation we define

$$h_{ij} = g_{i, j}^N, \qquad C^{j}_{i,p} = (1 - q^{-1})^N
\left(B_{i,p}^{p,j}\right)^{N(N-1)/2}.
$$

We are looking for families $(u_{ij})_{1\le i < j\le n +1}$ of elements
in $\ku \Gamma$ such that
\begin{equation}\label{combinquestion}
\Delta(u_{ij}) = u_{ij} \otimes 1 +
h_{i, j}  \otimes
u_{ij} +  \sum_{i < p < j}  C^{j}_{i,p}
u_{i, p} \, h_{p, j} \otimes u_{p, j},\text{  for all } 1 \leq i <  j \leq n+1.
\end{equation}

The coefficients $C^{j}_{i,p}$ satisfy the rule
\begin{equation}\label{ruleC}
 C^{j}_{is}C^{j}_{st} = C^{t}_{is}C^{j}_{it},
 \text{  for all } 1 \leq i < s < t < j \leq n+1.
\end{equation}

This follows from \eqref{2.2.01} and \eqref{2.2.02} since
$$ B^{sj}_{is} B^{tj}_{st} = B^{st}_{is} B^{tj}_{is} B^{tj}_{st} 
= B^{st}_{is} B^{tj}_{it}.$$

\begin{theorem}\label{lifting-combinator} Let $\Gamma$ be a finite abelian
group and
 $h_{ij}\in \Gamma$, $1\le i < j\le n +1$, a family of elements such
that
\begin{equation}\label{combin1}
h_{ij} = h_{i, p} \, h_{p, j}, \qquad    \text{if } \quad i < p < j.
\end{equation}
Let $C^{j}_{i,p} \in \ku^{\times}$, $1\le i < p < j\le n +1$, be a family of
elements satisfying \eqref{ruleC}.
Then the solutions $(u_{ij})_{1 \leq i < j \leq n+1}$ of \eqref{combinquestion}, $u_{ij} \in \ku \Gamma \text{ for all } i < j,$ have the form $(u_{ij}(\gamma))_{1 \leq i < j \leq n+1}$ where $\gamma = (\gamma_{ij})_{1 \leq i < j \leq n+1}$ is an arbitrary family of scalars $\gamma_{ij} \in \ku$ such that

\medbreak
\begin{equation}\label{hij}
\text{ for all } 1 \leq i < j \leq n+1, \ \gamma_{ij} = 0 \text{ if } h_{ij} = 1,
\end{equation}
and where the elements $u_{ij}(\gamma)$ are defined by induction on $j-i$ by

\medbreak
\begin{equation}\label{combin2}
u_{ij}(\gamma) = \gamma_{ij}(1 - h_{ij}) + \sum_{i < p < j} C_{ip}^{j} \gamma_{ip} u_{pj}(\gamma) \text{ for all } 1 \leq i < j \leq n+1.
\end{equation}
\end{theorem}

\pf We proceed by induction on $k$. We claim that the solutions $u_{ij} \in \ku \Gamma, 1 \leq i < j \leq n+1, j - i \leq k,$ of \eqref{combinquestion} for all $i < j$ with $j - i \leq k$ are given by arbitrary families of scalars $\gamma_{ij}, 1 \leq i < j \leq n+1, j - i \leq k$ such that
$$u_{ij} = \gamma_{ij}(1 - h_{ij}) + \sum_{i < p < j} C_{ip}^{j} \gamma_{ip} u_{pj} \text{ for all } 1 \leq i < j \leq n+1 \text{ with } j - i \leq k.$$

\medbreak

Suppose $k = 1$. For any $1 \leq i < n, j = i+1,$ $u_{i,i+1}$ is a solution of \eqref{combinquestion} if and only if $u_{i,i+1}$ is $(h_{i,i+1},1)$-primitive  in $\ku \Gamma$, that is $u_{i,i+1} = \gamma_{i,i+1}(1 - h_{i,i+1})$ for some $\gamma_{i,i+1} \in \ku.$ We may assume that $\gamma_{i,i+1} = 0$, if $h_{i,i+1} = 1.$

\medbreak

For the induction step, let $k > 2$. We assume that $\gamma_{ab} \in \ku, 1 \leq a < b \leq n+1, b-a \leq k-1,$ is a family of scalars with $\gamma_{ab} = 0$, if $h_{ab} = 1$, and that the family $u_{ab} \in \ku \Gamma, 1 \leq a < b \leq n+1, b-a \leq k-1,$ defined inductively by the $\gamma_{ab}$ by \eqref{combin2} is a solution of \eqref{combinquestion}. Let $1 \leq i < j \leq n+1,$ and $j-i = k.$ We have to show that
\begin{equation}\label{equationij}
\Delta(u_{ij}) = u_{ij} \otimes 1 + h_{ij} \otimes u_{ij} + \sum_{i < p < j} C_{ip}^{j} u_{ip} h_{pj} \otimes u_{pj}
\end{equation}
is equivalent to
\begin{equation}\label{solutionij}
u_{ij} = \gamma_{ij}(1 - h_{ij}) + \sum_{i < p < j} C_{ip}^{j} \gamma_{ip} u_{pj} \text{ for some } \gamma_{ij} \in \ku.
\end{equation}
We then may define $\gamma_{ij} = 0$ if $h_{ij} = 1.$

We denote
$$z_{ij}:= u_{ij} - \sum_{i < p < j} C_{ip}^{j} \gamma_{ip} u_{pj}.$$
Then \eqref{solutionij} is equivalent to
\begin{equation}\label{primitiveij}
\Delta(z_{ij}) = z_{ij} \otimes 1 + h_{ij} \otimes z_{ij}.
\end{equation}

For all $i < p < j$ we have $\Delta(u_{pj}) = u_{pj} \otimes 1 + h_{pj} \otimes u_{pj} + \sum_{p < s < j} C_{ps}^{j} u_{ps} h_{sj} \otimes u_{sj}$, since $j-p < k$. Using this formula for $\Delta(u_{pj})$ we compute
\begin{align}
&\Delta(z_{ij}) - z_{ij} \otimes 1 - h_{ij} \otimes z_{ij} \notag \\
&= \Delta(u_{ij}) - \sum_{i < p < j} C_{ip}^{j} \gamma_{ip} \Delta(u_{pj}) - z_{ij} \otimes 1 - h_{ij} \otimes z_{ij} \notag \\
&= \Delta(u_{ij}) - \sum_{i < p < j} C_{ip}^{j} \gamma_{ip}(u_{pj} \otimes 1 + h_{pj} \otimes u_{pj} + \sum_{p < s < j} C_{ps}^{j} u_{ps} h_{sj} \otimes u_{sj}) \notag\\
&- (u_{ij} - \sum_{i < p < j} C_{ip}^{j} \gamma_{ip} u_{pj}) \otimes 1 - h_{ij} \otimes  (u_{ij} - \sum_{i < p < j} C_{ip}^{j} \gamma_{ip} u_{pj}) \notag \\
 &= \Delta(u_{ij}) - u_{ij} \otimes 1 - h_{ij} \otimes u_{ij} \notag \\
&+ \sum_{i < p < j} C_{ip}^{j} \gamma_{ip} h_{ij} \otimes u_{pj} - \sum_{i < p < j} C_{ip}^{j} \gamma_{ip} h_{pj} \otimes u_{pj} - \sum_{i < p < s < j} C_{ip}^{j} C_{ps}^{j} \gamma_{ip} u_{ps} h_{sj} \otimes u_{sj}.  \notag
\end{align}

Therefore, \eqref{equationij} and  \eqref{solutionij} are equivalent if and only if the identity
\begin{align}
&\sum_{i < p < j} C_{ip}^{j} \gamma_{ip} h_{ij} \otimes u_{pj} - \sum_{i < p < j} C_{ip}^{j} \gamma_{ip} h_{pj} \otimes u_{pj} - \sum_{i < p < s < j} C_{ip}^{j} C_{ps}^{j} \gamma_{ip} u_{ps} h_{sj} \otimes u_{sj} \label{identityij} \\
&= -\sum_{i < p < j} C_{ip}^{j} u_{ip} h_{pj} \otimes u_{pj} \notag
\end{align}
holds.

To prove \eqref{identityij} we use \eqref{combin2} for all $i < p$, where $i < p < j$, that is $u_{ip} = \gamma_{ip} (1 - h_{ip}) + \sum_{i < s < p} C_{is}^{p} \gamma_{is} u_{sp}$. Then
\begin{align}
&\sum_{i < p < j} C_{ip}^{j} h_{pj} u_{ip} \otimes u_{pj} + \sum_{i < p < j} C_{ip}^{j} \gamma_{ip} h_{ij} \otimes u_{pj} - \sum_{i < p < j} C_{ip}^{j} \gamma_{ip} h_{pj} \otimes u_{pj} \notag \\
 &=\sum_{i < p < j} (C_{ip}^{j} h_{pj} (\gamma_{ip}(1-h_{ip}) + \sum_{i < s < p} C_{is}^{p} \gamma_{is} u_{sp}) + C_{ip}^{j} \gamma_{ip} (h_{ij} - h_{pj})) \otimes u_{pj} \notag \\
&= \sum_{i < p < j} C_{ip}^{j} h_{pj} \sum_{i < s < p} C_{is}^{p} \gamma_{is} u_{sp} \otimes u_{pj}, \text{ since } h_{pj}(1-h_{ip}) = h_{pj} - h_{ij} \text{ by } \eqref{combin1}, \notag \\
&= \sum_{i < s < p < j} C_{is}^{j} C_{sp}^{j} \gamma_{is} u_{sp} h_{pj} \otimes u_{pj}, \text{ since } C_{ip}^{j} C_{is}^{p} = C_{is}^{j} C_{sp}^{j} \text{ by } \eqref{ruleC}. \notag
\end{align}
This proves \eqref{identityij} by interchanging $s$ and $p$.
\epf

\begin{rmk}\label{explicit}
1) Let $\gamma = (\gamma_{ij})_{1 \le1q i < j \leq n+1}$ be an arbitrary 
family of scalars. Then it is easy to see that the family $u_{ij}(\gamma) 
\in \ku \Gamma, 1 \leq i < j \leq n+1,$ can be defined explicitly as follows:
$$u_{ij}(\gamma) = \sum_{i \leq p < j} \phi_{ip}^{j}(\gamma) (1 - h_{pj}) 
\text{ for all } i < j,$$
where
$$ \phi_{ip}^{j}(\gamma) = \sum_{i=i_{1} < \dots < i_{k} 
= p, k\geq 1} C_{i_{1},i_{2}}^{j} \dots C_{i_{k-1},i_{k}}^{j} 
\gamma_{i_{1},i_{2}} \dots \gamma_{i_{k-1},i_{k}} \gamma_{pj} 
\text{ for all } i \leq p < j$$
is a polynomial of degree $p$ in the free variables 
$(\gamma_{ij})_{1 \leq i < j \leq n+1}$.

\medbreak
2) Let $\gamma = (\gamma_{ij})_{1 \leq i < j \leq n+1}$ and 
$\widetilde{\gamma} =( \widetilde{\gamma}_{ij})_{1 \leq i < j \leq n+1}$ 
be families of scalars in $\ku$ satisfying \eqref{hij}. Assume that 
for all $i < j$, $u_{ij}(\gamma) = u_{ij}(\widetilde{\gamma}).$ 
Then  $\gamma = \widetilde{\gamma}$. This follows easily by induction on $j-i$ 
from \eqref{combin2}.
\end{rmk}

\begin{lema}\label{centraluij}
Assume the situation of Theorem \ref{lifting-combinator}. Let 
$\gamma = (\gamma_{ij})_{1 \leq i < j \leq n+1}$ be a 
family of scalars in $\ku$ satisfying \eqref{hij} and define 
$u_{ij} = u_{ij}(\gamma)$ for all ${1 \leq i < j \leq n+1}$ by \eqref{combin2}.

\begin{itemize}
\item[1)] The following are equivalent:
\begin{itemize}
\item[(a)] For all $i < j, u_{ij} = 0$ if $\chi_{ij}^{N} \neq \varepsilon.$

\medbreak
\item[(b)] For all $i < j, \gamma_{ij} = 0$ if $\chi_{ij}^{N} \neq \varepsilon.$
\end{itemize}
\item[2)] Assume that $h_{ij} = g_{ij}^{N}$ for all $i < j$. Then the following are equivalent:
\begin{itemize}
\item[(a)] For all $i < j, u_{ij} = 0$ if $\chi_{ij}^{N}(g_{l}) \neq 1$ for some $1 \leq l \leq n.$

\medbreak
\item[(b)] For all $i < j, \gamma_{ij} = 0$ if $\chi_{ij}^{N}(g_{l}) \neq 1$ for some $1 \leq l \leq n.$

\medbreak
\item[(c)] The elements $u_{ij}, 1 \leq i < j \leq n+1,$ are central in $\widehat{\mathfrak{B}}(V) \# \ku \Gamma.$
\end{itemize}
\end{itemize}
\end{lema}

\pf
1) follows by induction on $j-i$.

Suppose $j=i+1$. Then $u_{i,i+1} = \gamma_{i,i+1}(1-h_{i,i+1}).$ If $ h_{i,i+1}= 1$, then both $u_{i,i+1}$ and $ \gamma_{i,i+1}$ are $0$. If $ h_{i,i+1} \neq 1$, then  $u_{i,i+1}=0$ if and only if $ \gamma_{i,i+1}=0$.

The induction step follows in the same way from \eqref{combin2}, since for all $i < p < j$, if $\chi_{ij}^{N} \neq \varepsilon$, then $\chi_{ip}^{N} \neq \varepsilon$ or $\chi_{pj}^{N} \neq \varepsilon$, hence by induction $\gamma_{ip} =0$ or $u_{pj}=0$, and $u_{ij} = \gamma_{ij}(1-h_{ij}).$

2) Suppose that for all $i < p < j$, $u_{pj}$ is central in $\widehat{\mathfrak{B}}(V) \# \ku \Gamma,$ and let $1 \leq l \leq n$. Then
$$h_{ij} x_{l} = x_{l} \chi_{l}(h_{ij})h_{ij},$$
and we obtain from \eqref{combin2}
$$u_{ij} x_{l} = x_{l} \gamma_{ij}(1 - \chi_{l}(h_{ij})h_{ij}) + x_{l} \sum_{i < p < j} C_{ip}^{j} \gamma_{ip} u_{pj}.$$
Hence $u_{ij}$ is central in $\widehat{\mathfrak{B}}(V) \# \ku \Gamma$ if and only if
$\gamma_{ij} = \gamma_{ij}\chi_{l}(h_{ij})$ for all $1 \leq l \leq n$. Since the braiding is of type $A_{n}$ and the order of $q = \chi_{l}(g_{l})$ is $N$,
$$\chi_{l}(h_{ij})= \chi_{l}(g_{ij}^{N}) = \chi_{ij}^{-N}(g_{l}),$$
and the equivalence of (b) and (c) follows by induction on $j-i$.
The equivalence of (a) and (b) is shown as in 1).
\epf

\medbreak

\subsection{Classification of pointed Hopf algebras of type
$A_{n}$}\label{classif-A_n}

\

Using the previous results we will now determine exactly all finite-dimensional
pointed Hopf algebras of type $A_{n}$ (up to some exceptional cases).
We will find a big new class of deformations of $\mathfrak{u}_{q}^{\ge 0}(sl_{n + 1})$.

\medbreak

As before, we fix a natural number $n$, a finite abelian group $\Gamma$, an integer $N>2$, a root of unity $q$ of order $N$, $g_{1},\dots,g_{n} \in \Gamma, \chi_{1}, \dots, \chi_{n} \in \VGamma$ such that $q_{ij} = \chi_{j}(g_{i})$ for all $i,j$ satisfy \eqref{defqij}, and $V \in \ydg$ with basis $x_{i} \in V_{g_{i}}^{\chi_{i}}, 1 \leq i \leq n$.

We recall that $\widehat{\mathfrak{B}}(V)$ is the braided Hopf algebra in $\ydg$ generated by $x_{1}, \dots, x_{n}$ with the quantum Serre relations
 \eqref{2.2.40},
\eqref{2.2.400} and  \eqref{2.2.401}.

In  $\widehat{\mathfrak{B}}(V)$ we consider the iterated braided commutators $e_{i,j}, 1 \le i < j \le n+1$ defined inductively by \eqref{2.2.2} beginning with $e_{i,i+1} = x_{i}$ for all $i$.

Let $\mathbb{A}$ be the set of all families $(a_{i,j})_{1 \leq i < j \leq n+1}$ of integers $a_{i,j} \geq 0$ for all $1 \leq i < j \leq n+1$. For any $a \in \mathbb{A}$ we define
$$e^{a} := (e_{1,2})^{a_{1,2}} (e_{1,3})^{a_{1,3}} \dots (e_{n,n+1})^{a_{n,n+1}},$$
where the order in the product is the lexicographic order of the index pairs.
We begin with the PBW-theorem for $\widehat{\mathfrak{B}}(V)$.

\begin{theorem}\label{PBWhat}
The elements $e^{a}, a \in \mathbb{A}$, form a basis of the $\ku$-vector space $\widehat{\mathfrak{B}}(V)$.
\end{theorem}
\pf
The proof is similar to the proof of Theorem \ref{Theorem2.2.8}.
For general finite Cartan type the theorem can be derived from the PBW-basis of $U_{q}(\mathfrak{g})$ (see \cite{L3}) by changing the group and twisting as described in \cite[Section 4.2]{AS4}.
\epf

The following commutation rule for the elements $e_{ij}^{N}$ is crucial.

\begin{lema}\label{commutationN}
For all $1 \leq i < j \leq n+1, 1 \leq s < t \leq n+1,$
$$[e_{i,j},e_{s,t}^{N}]_{c} = 0, \text{ that is }\  e_{i,j} e_{s,t}^{N} = \chi_{s,t}^{N}(g_{i,j}) e_{s,t}^{N} e_{i,j}.$$
\end{lema}
\pf
Since $e_{i,j}$ is a linear combination of elements of the form $x_{i_{1}} \dots x_{i_{k}}$ with $k = j-i$ and $g_{i_{1}} \dots g_{i_{k}} = g_{ij}$, it is enough to consider the case when $j = i+1$.

To show $[e_{i,i+1},e_{s,t}^{N}]_{c} = 0$, we will distinguish several cases.

First assume that $(i,i+1) < (s,t)$. If $i+1 < s$ resp. $i = s$ and $i+1 < t$, then $[e_{i,i+1},e_{s,t}]_{c} = 0$ by \eqref{2.2.4} resp. \eqref{2.2.42}, and the claim follows.

If $i+1 = s$, we denote $x =e_{i+1,t}, y= e_{i,i+1}, z=e_{it}$ and $\alpha=\chi_{i+1,t}(g_{i}), \beta= \chi_{i+1,t}(g_{i,t}).$ Then
$$yx=\alpha xy + z, \text{ by } \eqref{2.2.6}, \text{ and } zx = \beta xz, \text{ by } \eqref{2.2.43}.
$$
Moreover, $\alpha = \chi_{i+1,t}(g_{i}) \neq \beta = \chi_{i+1,t}(g_{i,t}) =\chi_{i+1,t}(g_{i}) \chi_{i+1,t}(g_{i+1,t})$, and $\alpha^{N} = \beta^{N}$, since $\chi_{i+1,t}(g_{i+1,t}) = q$ by \eqref{2.2.04}. Therefore it follows from  \cite[Lemma 3.4]{AS4} that $y x^{N} = \alpha^{N} x^{N} y$, which was to be shown.

The claim is clear if $i = s$, and $i+1 = t$, since $\chi_{i}^{N}(g_{i}) = 1$.

It remains to consider the case when $(i,i+1) > (s,t)$. If $s<i$ and $t=i+1$, resp. $s<i$ and $i+1 <t$, resp. $s=i$ and $t < i+1$, then $e_{s,t} e_{i,i+1} = \chi_{i}(g_{s,t}) e_{i,i+1} e_{s,t}$ by \eqref{2.2.43}, resp.\eqref{2.2.41},resp.\eqref{2.2.42}. Hence in all cases, $e_{i,i+1} e_{s,t}^{N} = \chi_{i}^{-N}(g_{s,t}) e_{s,t}^{N} e_{i,i+1}$. This proves the claim $[e_{i,i+1},e_{s,t}^{N}]_{c} = 0$, since $\chi_{i}^{-N}(g_{s,t}) = \chi_{s,t}^{N}(g_{i}).$
\epf

\medbreak
We want to compute the dimension  of certain quotient algebras of
$\wtoba(V) \# \ku \Gamma$. Since this part of the theory works for
any finite Cartan type, we now consider more
generally a left $\ku \Gamma$-module algebra $R$ over any abelian group $\Gamma$
and assume that there are integers $P$ and $N_{i} >1$, elements $y_{i} \in R$,
$h_{i} \in \Gamma$,
$\eta_{i} \in \widehat{\Gamma}$, $1 \leq i \leq P$, such that
\begin{flalign}
&g \cdot y_{i} = \eta_{i}(g) y_{i}, \text{for all } g \in \Gamma, 1 \leq i \leq P.& \label{actioneta} \\
& y_{i} y_{j}^{N_{j}} = \eta_{j}^{N_{j}}(h_{i}) y_{j}^{N_{j}} y_{i}
\text{ for all } 1 \leq i,j \leq P.& \label{skewroot}\\
&\text{The elements }y_{1}^{a_{1}} \dots y_{P}^{a_{P}}, a_{1}, \dots,a_{P}
\geq 0, \text{ form a } \ku-\text{basis of } R. \label{basisy}&
\end{flalign}
Let $\mathbb{L}$ be the set of all $l = (l_{i})_{1 \leq i \leq P} \in \mathbb{N}^{P}$
such that $0 \leq l_{i} < N_{i}$ for all $1 \leq i \leq P$.
For $a = (a_{i})_{1 \leq i \leq P} \in \mathbb{N}^{P}$, we define
$$y^{a} = y_{1}^{a_{1}} \dots y_{P}^{a_{P}}, \text{ and } aN
= (a_{i}N_{i})_{1 \leq i \leq P}.$$
Then by \eqref{skewroot}, \eqref{basisy}, the elements
$$y^{l}y^{aN}, \qquad l \in \mathbb{L}, \  a \in \mathbb{N}^{P},$$
form a $\ku$-basis of $R$.

In the application to $\wtoba(V) \# \ku \Gamma$, $P$ is the number of
positive roots, and the $y_{i}$ play the role of the root vectors $e_{i,j}$.

To simplify the notation in the smash product algebra $R \# \ku \Gamma$, we identify $r \in R$ with $r \# 1$ and $v \in \ku \Gamma$ with $1 \# v$.
For $1 \leq i \leq P$, let $\widetilde{\eta}_{i} : \ku \Gamma \to \ku \Gamma$ be the algebra map defined by $\widetilde{\eta}_{i}(g) = \eta_{i}(g) g$ for all $g \in \Gamma$.  Then
$$v y_{i} = y_{i} \widetilde{\eta}_{i}(v)  \text{ for all } v \in \ku \Gamma.$$
We fix a family $u_{i}, 1 \leq i \leq P$, of elements in $\ku \Gamma$, and denote
$$ u^{a} := \prod_{1 \leq i \leq P} u_{i}^{a_{i}}, \text{ if } a = (a_{i})_{1 \leq i \leq P} \in \mathbb{N}^{P}.$$
Let $M$ be a free right $\ku \Gamma$-module with basis $m(l), l \in \mathbb{L}$.
We then define a right $\ku \Gamma$-linear map
$$\varphi : R \# \ku \Gamma \to M \text{ by } \varphi( y^{l} y^{aN}) := m(l) u^{a} \text{ for all } l \in \mathbb{L}, a \in \mathbb{N}^{P}.$$

\begin{lema}\label{rightideal}
Assume that
\begin{itemize}
\item $u_{i}$ is central in  $R \# \ku \Gamma$,
for all $1 \leq i \leq P$,  and

\medbreak
\item $u_{i} = 0$ if $\eta_{i}^{N_{i}}(h_{j}) \neq 1$ for some $1 \leq j \leq P$.
\end{itemize}

Then the kernel of $\varphi$ is a right ideal of $ R \# \ku \Gamma$ containing $y_{i}^{N_{i}} - u_{i}$ for all $1 \leq i \leq P$.
\end{lema}
\pf
By definition, $ \varphi(y_{i}^{N_{i}}) = m(0) u_{i} = \varphi(u_{i}).$

To show that the kernel of $\varphi$ is a right ideal, let
$$z = \sum_{l \in \mathbb{L}, a \in \mathbb{N}^{P}} y^{l} y^{aN} v_{l,a},\text{ where } v_{l,a} \in \ku \Gamma, \text {for all } l \in \mathbb{L}, a \in \mathbb{N}^{P},$$
be an element with $\varphi(z) = 0$. Then $\varphi(z) = \sum_{l,a} m(l) u^{a} v_{l,a} = 0$, hence
$$  \sum_{a \in \mathbb{N}^{P}} u^{a} v_{l,a} = 0,\text{for all } l \in \mathbb{L}.$$
Fix $1 \leq i \leq P$. We have to show that $\varphi(z y_{i}) = 0$.

For any $l \in \mathbb{L}$, we have the basis representation
$$y^{l} y_{i} = \sum _{t \in \mathbb{L}, b \in \mathbb{N}^{P}} \alpha_{t,b}^{l} y^{t} y^{bN}, \text { where } \alpha_{t,b}^{l} \in \ku \text{ for all }t \in \mathbb{L}, b \in \mathbb{N}^{P}.$$
Since $u^{a}$ is central in $R \# \ku \Gamma$,
\begin{equation}\label{ua}
u^{a} = \widetilde{\eta}_{i}(u^{a}) \text{ for all } a \in \mathbb{N}^{P}.
\end{equation}
For any $a = (a_{i})_{1 \leq i  \leq P} \in \mathbb{N}^{P}$ and
any family $({g}_{i})_{1 \leq i \leq  P}$ of elements in $\Gamma$ we define $\eta^{aN}(({g}_{i})) = \prod_{i} {\eta_{i}}^{a_{i}N_{i}}({g}_{i})$. Then by \eqref{skewroot}, for all $a,b \in \mathbb{N}^{P}$,
\begin{equation}\label{etaaN}
y^{aN} y_{i}= y_{i} y^{aN} \eta^{aN}(g^{a}), \text{ and } y^{bN} y^{aN} = y^{(a+b)N} \eta^{aN}(g^{b}),
\end{equation}
for some families of elements $g^{a}, g^{b}$ in $\Gamma$.

By a reformulation of our assumption,
\begin{equation}\label{reform}
 u^{a} \eta^{aN}((g_{i})) = u^{a} \text{ for any } a \in \mathbb{N}^{P} \text{ and family } (g_{i}) \text{ in } \Gamma.
\end{equation}
Using \eqref{etaaN} we now can compute
\begin{align}
z y_{i} &= \sum_{l,a} y^{l} y^{aN} v_{l,a} y_{i} = \sum_{l,a}  y^{l} y^{aN} y_{i} \widetilde{\eta}_{i}(v_{l,a}) \notag \\
 &= \sum_{l,a} y^{l} y_{i} y^{aN} \eta^{aN}(g^{a}) \widetilde{\eta}_{i}(v_{l,a}) \notag \\
 & = \sum_{l,a}\sum_{t, b} \alpha_{t,b}^{l} y^{t} y^{bN} e^{aN}  \eta^{aN}(g^{a}) \widetilde{\eta}_{i}(v_{l,a}) \notag \\
&= \sum_{l,a}\sum_{t, b} \alpha_{t,b}^{l} y^{t} y^{(a+b)N} \eta^{aN}(g^{a}) \eta^{aN}(g^{b}) \widetilde{\eta}_{i}(v_{l,a}). \notag
\end{align}
Therefore
\begin{align}
\varphi(z y_{i}) &= \sum_{t} m(t) \sum_{l,a,b} \alpha_{t,b}^{l} u^{a+b} \eta^{aN}(g^{a}) \eta^{aN}(g^{b}) \widetilde{\eta}_{i}(v_{l,a}) \notag \\
&= \sum_{t} m(t) \sum_{l,a,b} \alpha_{t,b}^{l} u^{a+b} \widetilde{\eta}_{i}(v_{l,a}), \text{ by } \eqref{reform}, \notag \\
&= \sum_{t} m(t) \sum_{b,l} \alpha_{t,b}^{l} u^{b} \sum_{a} u^{a} \widetilde{\eta}_{i}(v_{l,a}) \notag \\
&= \sum_{t} m(t) \sum_{b,l}  \alpha_{t,b}^{l} u^{b} \widetilde{\eta}_{i}(\sum_{a} u^{a} v_{l,a}), \text{ by } \eqref{ua}, \notag \\
&=0, \text{ since } \sum_{a} u^{a} v_{l,a} = 0. \notag
\end{align}
\epf

\begin{theorem}\label{basis}
Let $u_{i}, 1 \leq i \leq P,$ be a family of elements in $\ku \Gamma$, and $I$ the ideal in $R \# \ku \Gamma$ generated by all $y_{i}^{N_{i}} - u_{i}, 1 \leq i  \leq P.$ Let $A=( R \# \ku \Gamma)/I$ be the quotient algebra. Then the following are equivalent:

1) The residue classes of $y^{l}g, l \in \mathbb{L}, g \in \Gamma,$ form a $\ku$-basis of $A$.

\medbreak
2) $u_{i}$ is central in $R \# \ku \Gamma$ for all $1 \leq i \leq P$,
and $u_{i} = 0$ if $\eta_{i}^{N_{i}} \neq \varepsilon.$
\end{theorem}

\pf
$1) \Rightarrow 2): $ For all $i$ and $g \in \Gamma$, $g y_{i}^{N_{i}} = \eta_{i}^{N_{i}}(g) y_{i}^{N_{i}} g$, hence $u_{i} g = g u_{i} \equiv \eta_{i}^{N_{i}}(g) u_{i} g \mod I$. Since by assumption, $\ku \Gamma$ is a subspace of $A$, we conclude that $u_{i} = \eta_{i}^{N_{i}}(g) u_{i}$, and $u_{i} = 0$ if $\eta_{i}^{N_{i}} \neq \varepsilon$.

Similarly, for all $1 \leq i,j \leq n$, $ y_{i} y_{j}^{N_{j}} = \eta_{j}^{N_{j}}(h_{i}) y_{j}^{N_{j}} y_{i}$ by \eqref{skewroot}, hence
$y_{i} u_{j} \equiv \eta_{j}^{N_{j}}(h_{i}) u_{j} y_{i} \mod I$. Since we already know that $u_{i} = 0$ if $\eta_{j}^{N_{j}} \neq \varepsilon$, we see that $y_{i} u_{j} \equiv u_{j} y_{i} \mod I$. On the other hand $u_{j} y_{i} = y_{i} \widetilde{\eta}_{i}(u_{j})$. Then our assumption in 1) implies that $\widetilde{\eta}_{i}(u_{j}) = u_{j}$. In other words, $u_{j}$ is central in $R \# \ku \Gamma$.

$2) \Rightarrow 1)$: Let $J$ be the right ideal of $R \# \ku \Gamma$ generated by all $y_{i}^{N_{i}} - u_{i}, 1 \leq i  \leq P.$ For any $ 1 \leq i  \leq P$ and $g \in \Gamma,$
$$g(y_{i}^{N_{i}} - u_{i}) = y_{i}^{N_{i}} g \eta_{i}^{N_{i}}(g) - g u_{i}=(y_{i}^{N_{i}} - u_{i})  \eta_{i}^{N_{i}}(g) g,$$
since $gu_{i} = u_{i}  \eta_{i}^{N_{i}}(g) g$ by 2).

And for all $1 \leq i,j \leq P$,
$$y_{i}(y_{j}^{N_{j}} - u_{j}) = \eta_{j}^{N_{j}}(h_{i}) y_{j}^{N_{j}} y_{i} - y_{i}u_{j} = (y_{j}^{N_{j}} - u_{j}) \eta_{j}^{N_{j}}(h_{i}) y_{i},$$
since by 2) $y_{i} u_{j} = u_{j} y_{i} = u_{j} \eta_{j}^{N_{j}}(h_{i}) y_{i}.$

This proves $J = I$.

It is clear that the images of all $y^{l}g, l \in \mathbb{L}, g \in \Gamma,$ generate
the vector space $A$. To show linear independence, suppose
$$\sum_{l\in \mathbb{L}, g \in \Gamma} \alpha_{l,g} y^{l} g \in I,\text{ with }
\alpha_{l,g} \in \ku \text{ for all } l \in \mathbb{L}, g \in \Gamma.$$
Since $I = J$, we obtain from Lemma \ref{rightideal} that $\varphi(I) = 0$. Therefore,
$$0 = \varphi(\sum_{l\in \mathbb{L}, g \in \Gamma} \alpha_{l,g} y^{l} g) = \sum_{l\in \mathbb{L}, g \in \Gamma} \alpha_{l,g} m(l) g,$$
 hence $\alpha_{l,g} = 0$ for all $l,g$.
\epf

We come back to $A_{n}$.  Our main result in this Chapter is

\begin{theorem}\label{classification}
(i). Let  $\gamma=(\gamma_{i,j})_{1 \leq i < j \leq n+1}$ be any family of scalars in $\ku$ such that for any $i < j$, $\gamma_{i,j} = 0 \text{ if } g_{i,j}^{N} = 1 \text{ or } \chi_{i,j}^{N} \neq \varepsilon.$ Define $u_{i,j} = u_{i,j}(\gamma) \in \ku \Gamma, 1 \leq i < j \leq n+1,$ by \eqref{combin2}. Then
$$A_{\gamma} := (\widehat{\mathfrak{B}}(V) \# \ku \Gamma)/(e_{i,j}^{N} - u_{i,j} \mid 1 \leq i < j \leq n+1)$$
is a pointed Hopf algebra of dimension $N^{\frac{n(n+1)}{2}} \,\text{ord}(\Gamma)$ with
$ gr A_{\gamma} \simeq \mathfrak{B}(V) \# \ku \Gamma.$

(ii). Conversely, let $A$ be a finite-dimensional pointed Hopf
algebra   such that either

(a) $grA \simeq \mathfrak{B}(V) \# \ku \Gamma,$ and $N > 3,$ or

\smallbreak (b) the infinitesimal braiding of $A$ is of type
$A_{n}$ with $N>7$ and not divisible by 3.

\noindent Then $A$ is isomorphic to a Hopf algebra $A_{\gamma}$ in
(i).
\end{theorem}
\pf
(i). By Lemma \ref{centraluij}, the elements $u_{i,j}$ are central in $\widehat{U} := \widehat{\mathfrak{B}}(V) \# \ku \Gamma$, and $u_{i,j}=0$ if $g_{i,j}^{N} =1$ or $\chi_{i,j}^{N} \neq \varepsilon.$ Hence the residue classes of the elements $e^{l}g, l \in \mathbb{A}, 0 \leq l_{i,j} < N \text{ for all } 1 \leq i < j \leq n+1, g \in \Gamma,$ form a basis of $A_{\gamma}$ by Theorem \ref{basis}. By Theorem \ref{lifting-combinator}, the $u_{i,j}$ satisfy \eqref{combinquestion}. The ideal $I$ of  $\widehat{U}$ generated by all $e_{i,j}^{N} - u_{i,j}$ is a biideal, since
\begin{align}
\Delta(e_{i,j}^{N} - u_{i,j}) = &(e_{i,j}^{N} - u_{i,j}) \otimes 1 + g_{i,j}^{N} \otimes (e_{i,j}^{N} - u_{i,j}) \notag \\
&+ \sum_{i < p < j} C_{i,p}^{j}((e_{i,p}^{N} - u_{i,p}) g_{p,j}^{N} \otimes e_{p,j}^{N} + u_{i,p} g_{p,j}^{N} \otimes (e_{p,j}^{N} - u_{p,j})) \in I \otimes \widehat{U} + \widehat{U} \otimes I, \notag
\end{align}
by \eqref{coproduct} and \eqref{combinquestion}.

Since $A_{\gamma}$ is generated by group-like and skew-primitive elements,
and the group-like elements form a group, $A_{\gamma}$ is a Hopf algebra.

For all $1 \leq i \leq n$, let $a_{i} \in \text{gr}(A_{\gamma})(1)$ be
the residue class of $x_{i} \in (A_{\gamma})_{1}$. Define root vectors
$a_{i,j} \in \text{gr}(A_{\gamma}), 1 \leq i < j \leq N+1$ inductively
as in \eqref{2.2.1lift} and \eqref{2.2.2lift}. Then $a_{i,j}^{N} = 0$
in $\text{gr}(A_{\gamma})$ since $e_{i,j}^{N} =0$ in $A_{\gamma}$.
Therefore, by Theorem \ref{nichalg-typeA_n}, there is a surjective
Hopf algebra map
$$\mathfrak{B}(V) \# \ku \Gamma \to \text{gr}(A_{\gamma}) \text{ mapping } x_{i} \# g \text{ onto } a_{i}g, 1 \leq i \leq n, g \in \Gamma.$$
This map is an isomorphism, since $ \text{dim(gr}(A_{\gamma}))
= \text{dim}(A_{\gamma}) = N^{\frac{n(n+1)}{2}} \, |\Gamma|
= \text{dim}(\mathfrak{B}(V) \# \ku \Gamma)$ by Theorem \ref{Theorem2.2.8}.

(ii). As in Section \ref{lift-A_n}, we choose elements
$a_{i} \in (A_{1})_{g_{i}}^{\chi_{i}}$ such that
$\pi(a_{i}) = x_{i}, 1 \leq i \leq n$. By assumption resp. by
Lemma \ref{valenenA}, there is a Hopf algebra map
$$\phi : \widehat{\mathfrak{B}}(V) \# \ku \Gamma \to A, \phi(x_{i} \# g)
= a_{i}g, 1 \leq i \leq n, g \in \Gamma.$$
By Theorem \ref{degree1}, $A$ is generated in degree one, hence $\phi$
is surjective. We define the root vector $E_{i,j} \in A, 1 \leq i < j \leq n+1,$ by \eqref{2.2.1lift}, \eqref{2.2.2lift}. By Lemma \ref{vivenkgamma},
$E_{i,j}^{N} =:u_{i,j} \in \ku \Gamma$ for all $1\leq i < j \leq n+1.$ Then for all $g \in \Gamma$ and $i < j$, $ g E_{i,j}^{N} = \chi_{i,j}^{N}(g) E_{i,j}^{N} g,$ hence $g u_{i,j} =  \chi_{i,j}^{N}(g) u_{i,j} g,$ and $u_{i,j} = \chi_{i,j}^{N}(g) u_{i,j}.$ By \eqref{coproduct} and Theorem \ref{lifting-combinator} we therefore know that $u_{i,j} = u_{i,j}(\gamma)$ for all $i < j$, for some family $\gamma = (\gamma_{i,j})_{1 \leq i \leq J \leq n+1}$ of scalars in $\ku$ such that for all $i < j$, $\gamma_{i,j} = 0$ if $g_{i,j}^{N} = 1$ or $ \chi_{i,j}^{N} = \varepsilon.$ Hence $\phi$ indices a surjective Hopf algebra map $A_{\gamma} \to A$ which is an isomorphism since $\text{dim}(A_{\gamma}) = N^{\frac{n(n+1)}{2}} \, \text{ord}(\Gamma) = \text{dim}(A)$ by 1).
\epf
\begin{obs}
Up to isomorphism, $A_{\gamma}$ does not change if we replace each $x_{i}$ by a non-zero scalar multiple of itself. Hence in the definition of $A_{\gamma}$ we may always assume that
$$\gamma_{i,i+1} = 0 \text{ or } 1 \text{ for all } 1 \leq i \leq n.$$
\end{obs}

We close the paper with a very special case of Theorem \ref{classification}.
We obtain a large class of non-isomorphic
Hopf algebras which have exactly the same infinitesimal braiding as
$\mathfrak{u}_{q}^{\ge 0}(sl_{n})$. Here $q$ has order $N$, but the group is
$\prod_{i=1}^{n} \mathbb{Z}/(Nm_{i})$ and
not $(\mathbb{Z}/(N))^{n}$ as for $\mathfrak{u}_{q}^{\ge 0}(sl_{n + 1})$.

\begin{exa}
Let $N$ be $>2$, $q$ a root of unity of order $N$, and $m_{1}, \dots,m_{n}$ integers
$>1$ such that $m_{i} \neq m_{j}$ for all $i \neq j$. Let $\Gamma$ be the commutative
group generated by $g_{1},\dots,g_{n}$ with relations $g_{i}^{Nm_{i}} = 1$,
$1 \leq i \leq n.$ Define $\chi_{1} , \dots ,\chi_{n} \in \widehat{\Gamma}$ by
$$\chi_{j}(g_{i}) = q^{a_{ij}}, \text{ where } a_{ii} = 2 \text{ for all } i, a_{ij}
= -1 \text{ if } |i-j|=1, a_{ij} = 0 \text{ if } |i-j| \geq 2.$$
Then $\chi_{i,j}^{N} = \varepsilon$ and $g_{i,j}^{N} \neq 1$ for all $i < j.$
Thus for any family $\gamma = (\gamma_{i,j})_{1 \leq i < j \leq n+1}$ of
scalars in $\ku$, $A_{\gamma}$ in Theorem \ref{classification} has infinitesimal
braiding of type $A_{n}$.

Moreover, if $\gamma, \widetilde{\gamma}$ are arbitrary such families with
$\gamma_{i,i+1} = 1 = \widetilde{\gamma}_{i,i+1}$ for all $1 \leq i \leq n,$
then
$$A_{\gamma} \ncong A_{\widetilde{\gamma}}, \text{ if } \gamma \neq \widetilde{\gamma}.$$
\end{exa}
\pf
We let $\widetilde{x}_{i}$ and $\widetilde{e}_{i,j}$ denote the elements of $A_{\widetilde{\gamma}}$ corresponding to $x_{i}$ and $e_{i,j}$ in $A_{\gamma}$
as above, for all $i$ and $i < j$. Suppose $\phi : A_{\gamma} \to A_{\widetilde{\gamma}}$
is a Hopf algebra isomorphism. By Lemma \cite[Lemma 1.2]{AS3} there exist
non-zero scalars $\alpha_{1}, \dots, \alpha_{n} \in \ku$
and a permutation $\sigma \in \mathbb{S}_{n}$
such that $\phi(g_{i}) = g_{\sigma(i)}$ and $\phi(x_{i})
= \alpha_{i} \widetilde{x}_{\sigma(i)}$ for all $i$.
Since $\text{ord}(g_{i}) = m_{i}N \neq m_{j}N = \text{ord}(g_{j})$
for all $i \neq j$, $\sigma$ must be the identity, and $\phi$
induces the identity on $\Gamma$ by restriction. In particular,
$1 - g_{i}^{N} = \phi(x_{i}^{N}) = \alpha_{i}^{N} \widetilde{x}_{i}^{N}
= \alpha_{i}^{N} (1 -  g_{i}^{N})$, and $\alpha_{i}^{N} =1$ for all $i$.
Therefore we obtain for all $i < j$,
$$u_{i,j}(\gamma) = \phi(e_{i,j}^{N}) = \alpha_{i}^{N} \alpha_{i+1}^{N}
\dots \alpha_{j-1}^{N} \widetilde{e}_{i,j}^{N} = u_{i,j}(\widetilde{\gamma}),$$
and by Remark \ref{explicit}, 2), $\gamma = \widetilde{\gamma}.$ \epf

\end{document}